\documentclass{article}

\usepackage{lmodern}
\usepackage{amsthm,amsmath,amssymb,oldgerm,hyperref,mathrsfs,amsxtra, bm, accents}
\usepackage[toc,page]{appendix}
\usepackage[a4paper, margin=2.5 cm]{geometry}
\usepackage{mathtools}
\usepackage{dsfont}
\usepackage{bbm} 
\usepackage{theoremref}
\usepackage{amsopn}
\usepackage{aliascnt} 
\usepackage{enumerate}
\usepackage{stmaryrd} 
\usepackage{afterpage} 
\usepackage{extarrows} 
\usepackage{commath} 
\usepackage{pgf,tikz-cd}
\usepackage{mathrsfs}
\usetikzlibrary{arrows}
\usepackage{esint} 
\newcommand{\mynewtheorem}[2]{
\newaliascnt{#1}{dummy}
\newtheorem{#1}[#1]{#2}
\aliascntresetthe{#1}
\expandafter\def\csname #1autorefname\endcsname{#2}}
 
\theoremstyle{plain}
\mynewtheorem{theorem}{Theorem}
\mynewtheorem{proposition}{Proposition}
\mynewtheorem{corollary}{Corollary}
\mynewtheorem{lemma}{Lemma}

\newtheorem*{theorem*}{Theorem}
\theoremstyle{definition}
\mynewtheorem{definition}{Definition}
\mynewtheorem{example}{Example}
\theoremstyle{remark}
\mynewtheorem{remark}{Remark}
\newtheorem*{remark*}{Remark}

\def\equationautorefname~#1\null{equation~(#1)\null}

\numberwithin{equation}{section}
\DeclareMathOperator{\im}{Im}
\DeclareMathOperator{\re}{Re}
\DeclareMathOperator{\tr}{\mathrm{tr}}

\newcommand{\mf}{\mathfrak}
\newcommand{\mr}{\mathrm}
\newcommand{\mb}{\mathbb}
\newcommand{\ve}{\varepsilon}
\newcommand{\vr}{\varrho}
\newcommand{\ti}{\widetilde}
\newcommand{\ch}{\,\mathrm{ch}}
\newcommand{\sh}{\,\mathrm{sh}}
\newcommand{\Th}{\,\mathrm{th}}

\begin{document}

\title{Uniform sup-norm bounds on average for Siegel cusp forms}
\author{J\"urg Kramer and Antareep Mandal}
\date{\today}

\maketitle


\begin{abstract}
\noindent

Let $\Gamma\subsetneq \mathrm{Sp}_n(\mb{R})$ be an arithmetic subgroup of the symplectic group $\mathrm{Sp}_n(\mb{R})$
acting on the Siegel upper half-space $\mb{H}_n$ of degree $n$. Consider the $d$-dimensional space of Siegel cusp forms $\mathcal{S}_{\kappa}^n(\Gamma)$ of weight $\kappa$ for $\Gamma$ and let  $\{f_j\}_{1\leq j\leq d}$ be a basis of $\mathcal{S}_{\kappa}^n(\Gamma)$ orthonormal with respect to the Petersson inner product. In this paper we show using the heat kernel method that the sup-norm of the quantity $S_{\kappa}^{\Gamma}(Z):=\sum_{j=1}^{d}\det (Y)^{\kappa}\vert{f_j(Z)}\vert^2\,(Z\in\mb{H}_n)$ is bounded above by $c_{n,\Gamma} {\kappa}^{n(n+1)/2}$ when $M:=\Gamma\backslash\mb{H}_n$ is compact and by $c_{n,\Gamma} {\kappa}^{3n(n+1)/4}$ when $M$ is non-compact of finite volume, where $c_{n,\Gamma}$ denotes a positive real constant depending only on the degree $n$ and the group $\Gamma$. Furthermore, we show that this bound is uniform in the sense that if we fix a group $\Gamma_0$ and take $\Gamma$ to be a subgroup of $\Gamma_0$ of finite index, then the constant $c_{n,\Gamma}$ in these bounds depends only on the degree $n$ and the fixed group $\Gamma_0$. 
\end{abstract}


\section{Introduction}

Obtaining sup-norm bounds $\| \varphi \|_{\infty}$ of eigenfunctions $\varphi$ satisfying $\Delta_X \varphi + \lambda \phi =0$ for the Laplace--Beltrami operator $\Delta_X$ on a Riemannian manifold $X$ in terms of the eigenvalue $\lambda$ is a classical problem in spectral theory, for which local estimates exist \cite{hormander,seeger--sogge} that are essentially sharp in this level of generality. In arithmetic setting, these estimates are expected to be improved drastically\cite{sarnak-letter,sarnak-chaos}. Although major improvements over the classical general estimates have been obtained in such setting \cite{iwaniec--sarnak}, we are still far from the conjectured bound
\begin{align}\label{sarnak conjecture}
\| \varphi \|_{\infty}\ll_{\epsilon} \lambda^{\epsilon}\quad\quad(\epsilon>0)
\end{align}
for Hecke eigenforms. However, more interestingly, in this arithmetic setting, the sup-norm bound problem has been shown to have important connections to various fundamental queries in number theory such as the Lindel\"of hypothesis for the Riemann zeta function \cite{sarnak-chaos}, quantum ergodicity and entropy bounds \cite{key-4}, the subconvexity problem for $L$-functions \cite{key-6}, distribution of zeros of modular forms \cite{key-5} and the study of Arakelov invariants of arithmetic surfaces \cite{key-1,key-8,key-10}. This has created a sustained interest in the sup-norm bound problem in various number-theoretic aspects, one of which we address in this paper. 

Although we are far from obtaining \eqref{sarnak conjecture} for individual eigenforms, in the special setting of holomorphic cusp forms, in \cite{JK2}, optimal bounds have been obtained  on average over an orthonormal basis, without the assumption of strong arithmetic symmetries such as Hecke structure on the eigenforms. As this setup has a ready generalization in the case of Siegel modular forms, where the sup-norm bound has also recently been of some interest\cite{blomer-ikeda,blomer-spectral,cogdell,das-krishna}, we attempt here to extend this method and the results obtained in \cite{JK2} to the case of the Siegel upper half-space, which works to a large extent along with some significant non-trivialities that need special tools from the theory of harmonic analysis of semisimple Lie groups to get around. 


\subsection{Sup-norm bounds on $\mb{H}$}\label{Sup-norm bounds on H}
Let $\mb{H}:=\{z=x+i y \;|\; y>0\}$ denote the upper half-plane and $\Gamma\subsetneq \mathrm{SL}_2(\mb{R})$ denote a Fuchsian subgroup of the first kind. Let $\mathcal{V}_{\kappa}(\Gamma)$ denote the space of real analytic functions $\varphi\colon\mb{H}\rightarrow\mb{C}$ with the transformation behaviour 
\begin{align}\label{1-automorphy}
    \varphi(\gamma z)=\bigg(\frac{cz+d}{c\overline{z}+d}\bigg)^{\kappa/2}\varphi(z)\quad\quad\bigg(\gamma=\begin{pmatrix}a & b\\ c & d\end{pmatrix}\in\Gamma\bigg).
\end{align}
For $\varphi\in\mathcal{V}_{\kappa}(\Gamma)$, we set
\begin{align*}
\norm{\varphi}^2:=\int\limits_{\Gamma\backslash\mathbb{H}}\vert\varphi(z)\vert^2\;\frac{\dif x \wedge \dif y}{y^2},
\end{align*}
whenever it is defined.  Let $\mathcal{H}_{\kappa}(\Gamma):=\{\varphi\in \mathcal{V}_{\kappa}(\Gamma) \;|\; \norm{\varphi}<\infty\}$ denote the Hilbert space of square integrable functions in $\mathcal{V}_{\kappa}(\Gamma)$, equipped with the Petersson inner product. 
Let $\Delta^{(\kappa)}$ denote the Maa\ss{} Laplacian invariant with respect to the action of $\Gamma$ on $\varphi$ in \eqref{1-automorphy}. Then the operator $\Delta^{(\kappa)}$ acts on the smooth functions of $\mathcal{H}_{\kappa}(\Gamma)$ and extends to an essentially self-adjoint linear operator acting on a dense subspace of  $\mathcal{H}_{\kappa}(\Gamma)$. 

The eigenvalues for the Laplace equation $(\Delta^{(\kappa)}+\lambda)\varphi=0$ satisfy $\lambda\geq \kappa/2(1-\kappa/2)$ and in case $\lambda=\kappa/2(1-\kappa/2)$, the corresponding eigenfunction $\varphi$ in $\mathcal{H}_{\kappa}(\Gamma)$ can be shown to be of the form $\varphi(z)=y^{\kappa/2}f(z)$ with $f\in\mathcal{S}_{\kappa}(\Gamma)$, where $\mathcal{S}_{\kappa}(\Gamma)$ denotes the space of holomorphic cusp forms of weight $\kappa$ with respect to $\Gamma$.

Let $d=\dim(\mathcal{S}_{\kappa}(\Gamma))$ be the dimension of the space $\mathcal{S}_{\kappa}(\Gamma)$  and consider a  basis  $\{f_j\}_{1\leq j\leq d}$ of $\mathcal{S}_{\kappa}(\Gamma)$ orthonormal with respect to the Petersson inner product. Then, in order to obtain a sup-norm bound for the quantity
\begin{align*}
    S_{\kappa}^{\Gamma}(z):=\sum_{j=1}^{d} y^{\kappa} \vert{f_j(z)}\vert^2\quad(z\in\mb{H})
\end{align*}
in the weight aspect, using the spectral decomposition
of the weight-$\kappa$ heat kernel $K_t^{(\kappa,\Gamma)}$ and the $\mathbb{C}$-vector space isomorphism
\begin{align}\label{isomor}
\mathcal{S}_{\kappa}(\Gamma)\cong\ker\bigg(\Delta^{(\kappa)}+\frac{\kappa}{2}\bigg(1-\frac{\kappa}{2}\bigg)\mathrm{id}\bigg)
\end{align}
induced by the assignment $f\mapsto y^{\kappa/2}f $, one arrives at the important relation
\begin{align}\label{heatlimit}
    S_{\kappa}^{\Gamma}(z)=\lim_{t\rightarrow\infty}\exp\bigg(-\frac{\kappa}{2}\bigg(\frac{\kappa}{2}-1\bigg)\,t\bigg)\;K_t^{(\kappa,\Gamma)}(z,z)\quad\quad(\kappa\geq 2),
\end{align}
whence analyzing the heat kernel $K_t^{(\kappa,\Gamma)}$, in \cite{JK2} it is shown that
\begin{align*}
\sup_{z\in\mb{H}} S_{\kappa}^{\Gamma}(z)\leq
\begin{cases}
    c_{\Gamma}\, \kappa\quad\quad&(M\text{ compact}),\\
    c_{\Gamma}\, {\kappa}^{3/2}\quad\quad&(M\text{ non-compact of finite volume}),
\end{cases}
\end{align*}
where $c_{\Gamma}>0$ is a positive real number depending only on $\Gamma$. Furthermore, it is shown that this bound is uniform in the sense that if we fix a group $\Gamma_0\subsetneq \mathrm{SL}_2(\mb{R})$ and take $\Gamma$ to be a subgroup of $\Gamma_0$ of finite index, then the constant $c_{\Gamma}$ in these bounds depends only on the fixed group $\Gamma_0$.


\subsection{Sup-norm bounds on $\mb{H}_n$}
Let $\mathbb{H}_{n}:=\{Z=X+iY\in\mathbb{C}^{n\times n}\,\vert\,X,Y\in\mathrm{Sym}_{n}(\mathbb{R}):\,Y>0\}$ denote the Siegel upper half-space of degree $n$ and $\Gamma\subsetneq \mathrm{Sp}_n(\mb{R})$ be an arithmetic subgroup of the symplectic group $\mathrm{Sp}_n(\mb{R})$. Let $\mathcal{S}_{\kappa}^n(\Gamma)$ denote the space of all Siegel cusp forms of weight $\kappa$. Let $d=\dim(\mathcal{S}_{\kappa}^n(\Gamma))$ be the dimension of the space $\mathcal{S}_{\kappa}^n(\Gamma)$ and $\{f_j\}_{1\leq j\leq d}$ be an orthonormal  basis of $\mathcal{S}_{\kappa}^n(\Gamma)$ with respect to the Petersson inner product. We denote
\begin{align*}
    S_{\kappa}^{\Gamma}(Z):=\sum_{j=1}^{d}\det (Y)^{\kappa}\vert{f_j(Z)}\vert^2\quad(Z\in\mb{H}_n).
\end{align*}
In \cite{cogdell}, for $\Gamma=\Gamma_n=\mathrm{Sp}_n(\mb{Z})$, an asymptotic analysis of the Bergman kernel shows the bound
\begin{align}\label{cogdell bound}
    \sup\limits_{Z\in K} S_{\kappa}^{\Gamma}(Z)\asymp_{n} {\kappa}^{n(n+1)/2}\quad\quad(\kappa>2n)
\end{align}
in the weight aspect, where $K\subsetneq \mathscr{F}_n$ is any fixed compact subset of the standard fundamental domain $\mathscr{F}_n$ of $\mb{H}_n$ for $\Gamma_n$. 

In case of individual forms, on the basis of upper and lower bounds for certain specific kinds of Siegel cusp forms $F\in\mathcal{S}_{\kappa}^n(\Gamma)$, namely those coming from elliptic modular forms $f\in\mathcal{S}_{\kappa}(\Gamma)$ via Ikeda lifts, in the weight aspect, for $L^2$-normalized Siegel Hecke cusp forms $F$ for $\Gamma_n$ of (large) weight $\kappa$ and (fixed) genus $n$, Blomer in \cite{blomer-ikeda} conjectures
\begin{align*}
    \sup_{Z\in \mb{H}_n} \det(Y)^{\kappa}\vert F(Z) \vert^2\asymp_n {\kappa}^{n(n+1)/4}.
\end{align*}
 Combining this conjecture with Hashimoto's result \cite{hashimoto}
\begin{align*}
\dim_{\mathbb{C}} \mathcal{S}_{\kappa}^n\left(\Gamma\right)=2^{n(n-1)/2} \frac{\operatorname{vol}\left(\Gamma \backslash \mb{H}_{n}\right)}{(4 \pi)^{n(n+1) / 2}} {\kappa}^{n(n+1) / 2}+O({\kappa}^{n(n+1) / 2-1}),
\end{align*}
for $S_{\kappa}^{\Gamma}(Z)$ one conjectures
\begin{align}\label{blomer conjecture}
\sup_{Z\in\mb{H}_n} S_{\kappa}^{\Gamma}(Z)=O_{\Gamma}({\kappa}^{3n(n+1)/4}).
\end{align}
Recently, in \cite{das-krishna} it has been shown that 
\begin{align*}
\left.\begin{array}{l}
{\kappa}^{3n(n+1)/4} \\
{\kappa}^{3n(n+1)/4} \\
{\kappa}^{3n(n+1)/4}
\end{array}\right\} \ll_{n} \sup_{Z\in\mb{H}_n}\! S_{\kappa}^{\Gamma_n}(Z) \ll_{n, \epsilon}\left\{\begin{array}{ll}
{\kappa}^{3n(n+1)/4} & (n=1) \\
{\kappa}^{3n(n+1)/4 +\epsilon} & (n=2) \\
{\kappa}^{(5n-3)(n+1)/4 +\epsilon} & (n \geq 3)
\end{array}\right.,
\end{align*}
which establishes \eqref{blomer conjecture} for $n=1$ and $n=2$, but moves away from the optimal upper bound for $n>2$. 

 
\subsection{Statement of results}
The main result of this paper is the following theorem, which establishes the conjecture \eqref{blomer conjecture} for $n\geq 2$ by relating $S_{\kappa}^{\Gamma}(Z)$ with the heat kernel
$K_t^{(\kappa,\Gamma)}$ corresponding to the Siegel--Maa\ss{} Laplacian $\Delta^{(\kappa)}$ on $\Gamma\backslash\mb{H}_n$.

\begin{theorem}\label{main}
 Let $\Gamma\subsetneq \mathrm{Sp}_n(\mb{R})$ be an arithmetic subgroup and $\mathcal{S}_{\kappa}^n(\Gamma)$ denote the space of Siegel cusp forms of weight $\kappa$ on $\mb{H}_n$ with respect to $\Gamma$. Let $\{f_j\}_{1\leq j\leq d}$ be a basis of $\mathcal{S}_{\kappa}^n(\Gamma)$ orthonormal with respect to the Petersson inner product. Then, for all $n\geq 2$, we have
\begin{align*}
\sup_{Z\in\mb{H}_n}\! \sum_{j=1}^{d}\det (Y)^{\kappa}\vert{f_j(Z)}\vert^2\!\leq\!
\begin{cases}
    c_{n,\Gamma}\, {\kappa}^{n(n+1)/2}&(\Gamma\text{ cocompact}),\\
    c_{n,\Gamma}\, {\kappa}^{3n(n+1)/4}&(\Gamma\text{ cofinite}),
\end{cases}\quad(\kappa\geq n+1)
\end{align*}
where $c_{n,\Gamma}>0$ is a positive real number depending only on the degree $n$ of $\mb{H}_n$ and the group $\Gamma$.

Furthermore, this bound is uniform in the sense that if we fix a group $\Gamma_0\subsetneq \mathrm{Sp}_n(\mb{R})$ and take $\Gamma$ to be a subgroup of $\Gamma_0$ of finite index, then the constant $c_{n,\Gamma}$ in these bounds depends only on the degree $n$ and the fixed group $\Gamma_0$. 
\end{theorem}

This generalizes the theorems 4.2, 5.2, and 6.1 in \cite{JK2}. Furthermore, as $\kappa\geq n+1$, essentially the same arguments generalize the theorem 3.1 in \cite{JK1} to obtain for an orthonormal basis $\{f_j\}_{1\leq j\leq d}$ of $\mathcal{S}_{n+1}^n(\Gamma)$ and a positive real number $c_{n,\Gamma_0}>0$ depending only on $n$ and a fixed base space $M_0:=\Gamma_0\backslash \mb{H}_n$, the estimate
\begin{align*}
\frac{\dif \mu_B(Z)}{\dif \mu_S (Z)}=\sup_{Z\in\mb{H}_n} \sum_{j=1}^{d}\det (Y)^{n+1}\vert{f_j(Z)}\vert^2\ \leq c_{n,\Gamma_0},
\end{align*}
 where $\dif\mu_{B}$ denotes the volume form of the Bergman metric 
\begin{align*}
    \dif\mu_{B}(Z):=\sum_{j=1}^{d}\vert{f_j(Z)}\vert^2 \bigwedge\limits_{1\leq j\leq k\leq n} 
    \dif x_{j,k}\wedge\dif y_{j,k}
\end{align*}
and $\dif\mu_S$ denotes the volume form of the Siegel metric
\begin{align*}
\dif\mu_S(Z):=\frac{\bigwedge\limits_{1\leq j\leq k\leq n} \dif x_{j,k}\wedge \dif y_{j,k}}{\det(Y)^{n+1}} 
\end{align*}
on $\mb{H}_n$. 


\subsection{Strategy of the proof}

For the proofs, we follow the same general method developed in \cite{JK2}. Let  $\mathcal{V}_{\kappa}^n(\Gamma)$ denote the space of all real-analytic functions $\varphi\colon\mb{H}_n\rightarrow\mb{C}$, which have the transformation behaviour 
\begin{align}\label{n-automorphy}
\varphi(\gamma Z)&=\bigg(\frac{\det(C Z+D)}{\det(C\overline{ Z}+D)}\bigg)^{\kappa/2}\varphi( Z) \quad\quad\bigg(\gamma=\begin{pmatrix}A & B\\ C & D\end{pmatrix}\in\Gamma\bigg).
\end{align}
For $\varphi\in\mathcal{V}_{\kappa}^n(\Gamma)$, we set
\begin{align*}
\norm{\varphi}^2:=\int\limits_{\Gamma\backslash\mathbb{H}_n}\vert\varphi( Z)\vert^2\dif\mu_{n}(Z),
\end{align*}
whenever it is defined. Let $\mathcal{H}_{\kappa}^n(\Gamma):=\{\varphi\in \mathcal{V}_{\kappa}^n(\Gamma) \;|\; \norm{\varphi}<\infty\}$ denote the Hilbert space  of square integrable functions in $\mathcal{V}_{\kappa}^n(\Gamma)$, equipped with the Petersson inner product. Let $\Delta^{(\kappa)}$ denote the Siegel-Maa\ss{} Laplacian invariant with respect to the action of $\Gamma$ on $\varphi$ in \eqref{n-automorphy}. Then, the operator $\Delta^{(\kappa)}$ acts on the smooth functions of  $\mathcal{H}_{\kappa}^n(\Gamma)$ and extends to an essentially self-adjoint linear operator acting on a dense subspace of  $\mathcal{H}_{\kappa}^n(\Gamma)$. 

The eigenvalues for the Laplace equation $(\Delta^{(\kappa)}+\lambda)\,\varphi=0$ satisfy the inequality $\lambda\geq (n\kappa/4)\big((n+1)-\kappa\big)$, with the equality 
$\lambda= (n\kappa/4)\big((n+1)-\kappa\big)$ being attained if and only if $\varphi$ is of the form
$\varphi(Z)=\det(Y)^{\kappa/2}f(Z)$ for some Siegel cusp form $f\in \mathcal{S}_{\kappa}^n(\Gamma)$ of weight $\kappa$, i.e.,  the $\mb{C}$-vector space isomorphism
\begin{align}\label{rekernel}
\mathcal{S}_{\kappa}^n(\Gamma)\cong\ker\bigg(\Delta^{(\kappa)}+\frac{n\kappa}{4}((n+1)-\kappa)\mathrm{id}\bigg)
\end{align}
induced by the assignment $f\mapsto \det ( Y)^{\kappa/2}f$ holds (See \cite[corollary 5.4]{siegel--maass}).

Then, in a manner similar to \eqref{spectral decomposition}, we use the spectral decomposition of the heat kernel $K_t^{(\kappa,\Gamma)}$ corresponding to the Siegel--Maa\ss{} Laplacian $\Delta^{(\kappa)}$ on $\Gamma\backslash\mb{H}_n$ to generalize the relation in \eqref{heatlimit} to obtain
\begin{align}\label{limit-connection}
    S_{\kappa}^{\Gamma}(Z)=\lim\limits_{t\rightarrow\infty}\exp\big(-\frac{n \kappa}{4}(\kappa-(n+1))\,t\big) \,K_{t}^{(\kappa,\Gamma)}(Z,Z).
\end{align}
As, for $t>0$, both the function $\exp(-n\kappa(\kappa-(n+1))\,t/4)$ and  the heat kernel $K_{t}^{(\kappa,\Gamma)}$ are monotonically decreasing in $t$, from \eqref{limit-connection}, one also obtains the inequality
\begin{align}\label{limit-inequality}
    S_{\kappa}^{\Gamma}(Z)\leq \exp\big(-\frac{n\kappa}{4}(\kappa-(n+1))\,t\big) \,K_{t}^{(\kappa,\Gamma)}(Z)\quad\quad(t>0),
\end{align}
whence analyzing the heat kernel $K_t^{(\kappa,\Gamma)}$, we arrive at the results stated in \autoref{main}. 

The non-triviality in extending these results from $n=1$ to $n\geq 1$ lies in constructing the heat kernel $K_t^{(\kappa)}$ corresponding to the Siegel--Maa\ss{} Laplacian $\Delta^{(\kappa)}$ on $\mb{H}_n$ from which the heat kernel $K_t^{(\kappa,\Gamma)}$ on $\Gamma\backslash\mb{H}_n$ is obtained by periodization. We use a method of calculating spherical functions on real semisimple Lie groups by reducing them to the complex case developed by Flensted-Jensen in \cite{FLENSTED-JENSEN} to construct a spherical function for $\Delta^{(\kappa)}$ on $\mb{H}_n$ and then use the traditional method of obtaining the heat kernel from a spherical function to construct  $K_t^{(\kappa)}$. The spherical function and the ensuing heat kernel so constructed are not totally explicit, as they involve a change of variable that is somewhat implicit in nature. In the end, we got around this difficulty in the cocompact case by using a counting function argument to estimate the periodization sum in $K_t^{(\kappa,\Gamma)}$ with an integral of $K_t^{(\kappa)}$ over the radial coordinates, which allows us to change back from the implicit change of variables. In the cofinite case, we consider a limiting case of $K_t^{(\kappa)}$, which can be constructed explicitly. Thankfully, these special cases seem to suffice for our purpose. 


\subsection{Brief outline of the paper}
In section 2 we gather some basic preliminaries on the symplectic group, Siegel upper half-space and Siegel modular forms for later use in our calculations. In section 3 we gather the background on spherical functions and construction of the heat kernel on symmetric spaces. Then we use the Flensted-Jensen reduction method of calculating spherical functions on real semisimple groups via complex semisimple groups to construct the spherical function as well as the heat kernel on the Siegel upper half-space. All these calculations being for the Laplace--Beltrami operator $\Delta$, next we correct them for weight-$\kappa$ to obtain the heat kernel on Siegel upper half-space corresponding to the Siegel--Maass Laplacian $\Delta^{(\kappa)}$. Finally, in section 4, we analyze this weight-$\kappa$ heat kernel to obtain uniform sup-norm bounds for Siegel cusp forms on average over an orthonormal basis in both cocompact and cofinite cases.


\subsection{Acknowledgement}
The first author acknowledges support from the DFG Cluster of Excellence MATH+ and the second author acknowledges support from the Indian Institute of Technology, Madras. In parts, the material of this manuscript is contained in the doctoral dissertation \cite{thesis} of the second author completed under the supervision of the first author at the Humboldt-Universität zu Berlin, during which he acknowledges support from the university as well as from the Berlin Mathematical School. Both authors thank Valentin Blomer, Anilatmaja Aryasomayajula and Aprameyan Parthasarathy for inspiring discussions related to the material presented here. Both authors certify that they have no affiliations with or involvement in any organization or entity with any financial interest or non-financial interest in the subject matter or materials discussed in this manuscript.


\section{Background on Siegel modular forms}

In this section, we gather some basic preliminaries on the symplectic group, Siegel upper half-space and Siegel modular forms from some standard references such as \cite{Freitag}, \cite{Klingen} and \cite{Terras2}. 


\subsection{Siegel upper half-space}\label{SEC Siegel geometry}

For $n\in\mathbb{N}_{>0}$ and a commutative ring $R$, let $R^{m\times n}$ denote the set of $(m\times n)$-matrices with entries in $R$ and 
$\mathrm{Sym}_{n}(R)$ denote the set of symmetric matrices in $R^{n\times n}$.
The Siegel upper half-space 
$\mathbb{H}_{n}$ of degree $n$ is then defined by
\begin{align*}
\mathbb{H}_{n}:=\{Z=X+iY\in \mathbb{C}^{n\times n}\,\vert\,X,Y\in\mathrm{Sym}_{n}(\mathbb{R}):\,Y>0\}.
\end{align*}
The symplectic group $\mathrm{Sp}_{n}(\mathbb{R})$ of degree $n$ is defined by 
\begin{align*}
\mathrm{Sp}_{n}(\mathbb{R}):=\{g\in\mathbb{R}^{2n\times 2n}\,\vert\,\,g^{t}J_{n}g=J_{n}\},
\end{align*}
where $J_{n}\in\mathbb{R}^{2n\times 2n}$ is the skew-symmetric matrix 
\begin{align*}
J_{n}:=\begin{pmatrix}0&\mathbbm{1}_{n}\\-\mathbbm{1}_{n}&0\end{pmatrix}
\end{align*}
with $\mathbbm{1}_{n}$ denoting the identity matrix of $\mathbb{R}^{n\times n}$. The group $\mathrm{Sp}_{n}(\mathbb{R})$ acts by the 
symplectic action
\begin{align}
\label{sympaction}
\mathbb{H}_{n}\ni Z\mapsto g Z=(AZ+B)(CZ+D)^{-1}\qquad\big(g=\big(\begin{smallmatrix}A&B\\C&D\end{smallmatrix}\big)\in\mathrm
{Sp}_{n}(\mathbb{R})\big)
\end{align}
on $\mathbb{H}_{n}$. Under this action $\mathrm{Im}(Z)$ transforms as
\begin{align}
\label{imaginary transform}
\mathrm{Im}(g Z)=(CZ+D)^{-t}\,\mathrm{Im}(Z)(C\overline{Z}+D)^{-1},
\end{align}
which, on taking determinants on both sides, gives
\begin{align*}
\det(\mathrm{Im}(g Z))=\frac{\det(\mathrm{Im}(Z))}{\vert\det(CZ+D)\vert^{2}}.
\end{align*}
Similarly, taking matrix-differentials on both sides of the symplectic action~\eqref{sympaction}, it is easy to see that under this action, the matrix-differential 
form $\mathrm{d}Z$ transforms as
\begin{align}
\label{difZ calc2}
\mathrm{d}(g Z)=(CZ+D)^{-t}\,\mathrm{d}Z\,(CZ+D)^{-1}.
\end{align}
The arclength $\mathrm{d}s^{2}_{n}$ and the volume form $\mathrm{d}\mu_{n}$ on $\mathbb{H}_{n}$ in terms of $Z=(z_{j,k})_{1\leq j\leq k\leq n}\in 
\mathbb{H}_{n}$ are given by
\begin{equation}\label{siegel-metric}
\begin{aligned}
\mathrm{d}s^{2}_{n}(Z)&=\tr(Y^{-1}\,\mathrm{d}Z\,Y^{-1}\,\mathrm{d}\overline{Z})\quad\qquad(Z=X+iY), \\[1mm]
\mathrm{d}\mu_{n}(Z)&=\frac{\bigwedge\limits_{1\leq j\leq k\leq n}\mathrm{d}x_{j,k}\wedge\mathrm{d}y_{j,k}}{\det(Y)^{n+1}}\qquad(z_{j,k}=x_{j,k}+
iy_{j,k}).
\end{aligned}
\end{equation}
From equations~\eqref{imaginary transform} and~\eqref{difZ calc2} it is obvious that the arclength $\mathrm{d}s^{2}_{n}$ and the volume form 
$\mathrm{d}\mu_{n}$ on $\mathbb{H}_{n}$ given by the above equations are invariant under the symplectic action. Corresponding to this metric, 
we have the Laplace--Beltrami operator
\begin{align*}
    \Delta=\tr\Bigg( Y\bigg(\bigg(Y\dpd{}{ X}\bigg)^t\dpd{}{ X}+\bigg(Y\dpd{}{ Y}\bigg)^t\dpd{}{ Y}\bigg)\Bigg)
\end{align*}
on $\mathbb{H}_{n}$, called the Siegel Laplacian, which is also invariant under the symplectic 
action. 

The geodesic distance $s(Z,W)$ between the points $Z,\,W\in\mb{H}_n$ is given by
\begin{align*}
s(Z, W)=\sqrt{2}\Bigg(\sum_{j=1}^{n} \log ^{2} \bigg(\frac{1+\sqrt{\rho_{j}}}{1-\sqrt{\rho_{j}}}\bigg)\Bigg)^{1 / 2},
\end{align*}
where  $\rho_j\,(1\leq j\leq n)$ are the eigenvalues of the cross-ratio matrix
\begin{align}\label{siegel distance}
    \rho(Z, W)=(Z-W)(\overline{Z}-W)^{-1}(\overline{Z}-\overline{W})(Z-\overline{W})^{-1}\quad(Z,W\in\mb{H}_n).
\end{align}

\begin{remark}
Due to the action \eqref{sympaction} of $\mr{Sp}_n(\mb{R})$ on $\mb{H}_n$, the Siegel upper half-space can be viewed as a Riemannian globally symmetric space $G_0/K_0$, where $G_0=\mr{Sp}_n(\mb{R})$ and $K_0=\mr{Sp}_n(\mb{R})\cap O_{2n}(\mb{R})\cong U_n$ is the maximal compact subgroup of $\mr{Sp}_n(\mb{R})$ fixing the origin $i\mathbbm{1}_n\in \mb{H}_n$. An explicit correspondence between $G_0/K_0$ and $\mb{H}_n$ is given by $gK_0\mapsto g\,i\mathbbm{1}_n$. Under this identification, the action of $G_0$ by left translation on $G_0/K_0$ translates to the symplectic action \eqref{sympaction} of $G_0$ on $\mb{H}_n$. The space $G_0/K_0$ is equipped with a natural metric coming from the Killing form on the Lie algebra $\mf{g}_0=\mf{sp}_n(\mb{R})$ of $G_0$. Under this identification, this metric takes the form of the metric given by \eqref{siegel-metric} on $\mb{H}_n$. 
\end{remark}

\subsection{Arithmetic subgroups and fundamental domains}\label{SEC arithmetic subgroups}

A subgroup $\Gamma$ of $\mathrm{Sp}_n(\mathbb{R})$ is called \emph{discrete} if it acts discontinuously on $\mathbb{H}_n$, i.e., the orbit $\Gamma Z =\{\gamma Z \;|\; \gamma \in \Gamma\}$
has no accumulation point in $\mathbb{H}_n$, or equivalently, for any two compact sets $K_1,K_2\subset\mathbb{H}_n$, the set $\{\gamma\in\Gamma\;|\; \gamma(K_1)\cap K_2\neq \emptyset\}$ is finite. The most important example of a discrete subgroup of 
$\mathrm{Sp}_n(\mathbb{R})$ is the Siegel modular group $\Gamma_n:=\mathrm{Sp}_n(\mathbb{Z})$. 
\begin{definition}
A subgroup $\Gamma\subsetneq \mathrm{Sp}_n(\mathbb{R})$ is called an arithmetic subgroup if $\Gamma$ is commensurable to $\Gamma_n$, i.e., the group $\Gamma\cap\Gamma_n$ has finite index in both $\Gamma$ and $\Gamma_n$. 
\end{definition}
Because of their commensurability with the discrete subgroup $\Gamma_n$, arithmetic subgroups of $\mathrm{Sp}_n(\mathbb{R})$ are also discrete. 

Any arithmetic subgroup $\Gamma$ of $\mathrm{Sp}_n(\mathbb{R})$ has a fundamental domain, but it is not unique. A fundamental domain of the Siegel modular group $\Gamma_n$ can be explicitly constructed by means of reduction theory applied to the positive definite imaginary part $Y$ of $Z\in\mb{H}_n$. A vector $h^t=(h_1,h_2,\ldots,h_g)\in\mathbb{Z}^n$ is called primitive if for $1\leq k\leq n$, we have $\mathrm{gcd}(h_{k},\ldots,h_{n})=1$. 
\begin{definition}
A positive definite matrix $Y=(y_{j,k})_{1\leq j,k\leq n}\in\mr{P}_n$ is called Minkowski reduced if it satisfies $y_{k,k+1}\geq0\,(1\leq k\leq n-1)$ and for all primitive vectors $h\in\mathbb{Z}^n$, we have $h^t Y h\geq y_{k,k}\,(1\leq k\leq n)$. 
\end{definition}

\begin{proposition}
A Minkowski reduced positive definite matrix $Y$ satisfies the properties
\begin{itemize}
\item[(i)] $ y_{1,1}\leq y_{2,2}\leq \ldots \leq y_{n,n}$,
\item[(ii)] $\abs{2 y_{j,k}}\leq y_{j,j} \quad(1\leq j< k\leq n)$,
\item[(iii)] $y_{k,k+1}\geq0 \quad(1\leq k\leq n-1)$,
\item[(iv)] there exists a positive number $c_1(n)$ depending only on $n$ such that 
\begin{align*}
\det(Y)\leq \prod \limits_{j=1}^n y_{j,j}\leq c_1(n)\det(Y), 
\end{align*}
\item[(v)]  there exists a positive number $c_2(n)$ depending only on $n$ such that 
\begin{align*}
c_2(n)^{-1} Y<Y^{\mathrm{D}}<c_2(n) Y,
\end{align*}
where $Y^{\mathrm{D}}$ denotes the diagonal matrix made up of the diagonal elements $y_{1,1}, \ldots, y_{n,n}$ of the matrix $Y$, i.e.,
\begin{align*}
Y^{\mathrm{D}}=\begin{pmatrix}
{y_{1,1}} & {} & {0} \\
{} & {\ddots} & {} \\
{0} & {} & {y_{n,n}}
\end{pmatrix}.
\end{align*}
\end{itemize}
\end{proposition}
\begin{proof}
See \cite[Satz 2.5, Folgerung 2.6]{Freitag}, \cite[page 20]{Klingen}. 
\end{proof}

\begin{proposition}
 The set of points $Z=X+iY\in\mathbb{H}_n$ satisfying the following criteria forms a fundamental domain $\mathscr{F}_n$ of the Siegel modular group $\Gamma_n$:
\begin{itemize}
\item[(i)] $|\det(CZ+D)|\geq1$ for all $(\begin{smallmatrix}A&B\\C&D\end{smallmatrix})\in\Gamma_n$,
\item[(ii)]$Y=\im(Z)$ is Minkowski reduced,
\item[(iii)] for all $1\leq j,k\leq n$, the matrix $X=(x_{j,k})_{1\leq j,k\leq n}$ satisfies $|x_{j,k}|\leq 1/2$.
\end{itemize}
\end{proposition}
\begin{proof}
See \cite[Satz 2.9]{Freitag}.
\end{proof}

The fundamental domain $\mathscr{F}_n$ of $\Gamma_n$ is called the standard fundamental domain of $\Gamma_n$ and the matrices $Z\in\mathscr{F}_n$ are called Siegel reduced.

\begin{proposition}
If $Z=X+iY\in\mb{H}_n$ is also Siegel reduced, then $Y$ satisfies the properties
\begin{itemize}
    \item [(i)] $y_{1,1}\geq \sqrt{3}/2$.
    \item [(ii)] there exists a constant $c_3(n)>0$ depending only on $n$, such that $Y\geq c_3(n)\mathbbm{1}_n$.
\end{itemize}
\end{proposition}
\begin{proof}
See \cite[Hilfssatzs 2.11 ,2.12]{Freitag}.
\end{proof}

Using the fundamental domain $\mathscr{F}_n$ of the Siegel modular group $\Gamma_n$, we can construct fundamental regions of other arithmetic subgroups of $\Gamma$ of $\mathrm{Sp}_n(\mathbb{R})$. Consider the space of the left cosets $\{\Gamma\gamma\;|\;\gamma\in\Gamma_n\}$. Since $\Gamma\gamma_1=\Gamma\gamma_2$ if and only if $(\Gamma\cap\Gamma_n)\gamma_1=(\Gamma\cap\Gamma_n)\gamma_2$ and 
$[\Gamma_n:\Gamma\cap\Gamma_n]<\infty$, they have a finite system of representatives $\gamma_1,\gamma_2,\ldots,\gamma_m\,(m\in\mb{N}_{\geq 1})$. 
Then,
\begin{align}\label{arithmetic fundamental}
\mathscr{F}=\bigcup\limits_{j=1}^m\gamma_j\mathscr{F}_n
\end{align}
is a fundamental region of $\Gamma$. 


\subsection{Boundary of the Siegel upper half space}\label{SEC Boundary}

The Siegel upper half-space $\mathbb{H}_n$ can be realized as bounded domain 
$\mathbb{D}_n=\{\zeta\in \mathrm{Sym}_n(\mathbb{C})\,\vert\, \zeta \bar{\zeta}<\mathbbm{1}_n\}$
through the Cayley transformation $l\colon\mathbb{H}_n \xlongrightarrow{\sim}\mathbb{D}_n$ given by the assignment
\begin{equation}
Z\mapsto \zeta=(Z-i\mathbbm{1}_n)(Z+i\mathbbm{1}_n)^{-1},
\end{equation}
whose inverse $l^{-1}:\mathbb{D}_n\longrightarrow \mathbb{H}_n$ is given by the assignment
\begin{equation}
\zeta\mapsto Z=i(\mathbbm{1}_n+\zeta)(\mathbbm{1}_n-\zeta)^{-1}. 
\end{equation}

The topological closure of $\mathbb{D}_n$ is given by $\overline{\mathbb{D}}_n=\{\zeta\in \mathrm{Sym}_n(\mathbb{C})\,\vert\, \zeta\bar{\zeta}\leq\mathbbm{1}_n\}$. Through Cayley transformation, the symplectic action on $\mathbb{H}_n$ induces an analogous action of $\mathrm{Sp}_{n}(\mathbb{R})$ on $\mathbb{D}_n$ given by
\begin{align*}
\begin{pmatrix}A&B\\C&D\end{pmatrix}\zeta=\Big((A-iC)(\zeta+1)+i(B-iD)(\zeta-1)\Big)\Big((A+iC)(\zeta+1)+i(B+iD)(\zeta-1)\Big)^{-1},
\end{align*}
which extends to $\overline{\mathbb{D}}_n$ (see \cite[page 15]{Namikawa}). 

Two points $\zeta,\eta\in\overline{\mathbb{D}}_n$ are called equivalent if they can be connected by a finite number of holomorphic curves.
\begin{definition}
    A maximal subset in $\overline{\mathbb{D}}_n$ of mutually equivalent points is called a boundary component of $\mathbb{D}_n$.  
\end{definition}
The space $\overline{\mathbb{D}}_n$ is divided into a disjoint union of boundary components. Moreover, the symplectic action transforms one boundary component to another. Therefore, the division of $\overline{\mathbb{D}}_n$ into boundary components is invariant under the symplectic action.

For an integer $0\leq j \leq n$, let 
\begin{align*}
\mathbb{D}_n^j=\bigg\{\begin{pmatrix} \zeta_j & 0\\0 & \mathbbm{1}_{n-j}\end{pmatrix} \;\bigg\vert\; \zeta_j\in\mb{D}_j\bigg\}\cong \mb{D}_j
\end{align*}
Then for all $0\leq j\leq n$, $\mathbb{D}_n^j$ is a boundary component. In particular, $\mathbb{D}_n\cong \mathbb{D}_n^n$ itself is a boundary component. As
\begin{align*}
\overline{\mathbb{D}}_n=\bigcup\limits_{0\leq j\leq n} \mathrm{Sp}_{n}(\mathbb{R}) \mathbb{D}_n^j,
\end{align*} 
any boundary component $\mathbb{P}$ of $\mathbb{D}_n$ can be realized as $\mathbb{P}=g \mathbb{D}_n^j$
 for some $g\in\mathrm{Sp}_{n}(\mathbb{R})$ and $0\leq j\leq n$ (see \cite[page 17]{Namikawa}). Hence, we call $\mathbb{D}_n^j$ a standard boundary component.
 
As $\mathbb{D}_n^j$ is isomorphic to the bounded realization $\mathbb{D}_j$ of the degree $j$ Siegel upper half-space, we say $\mathbb{P}=g\mathbb{D}_n^j$ is a boundary component of degree $j$. If $j<n$, we call $\mathbb{P}$ a proper boundary component. For two boundary components $\mathbb{P}^{j}, \mathbb{P}^{k}$ of degrees $j,k$ respectively, we write $\mathbb{P}^{j}<\mathbb{P}^{k}$ if $\mathbb{P}^j\subset \overline{\mathbb{P}}_k$. In that case there exists $g\in\mathrm{Sp}_{n}(\mathbb{R})$ such that $g \mathbb{P}^{j}=\mathbb{D}_n^j$, $g \mathbb{P}^{k}=\mathbb{D}_n^k$ and $j\leq k$. 
\begin{remark}\label{boundary chain}
This result can be extended to the case of a chain of boundary components 
 \begin{align*}
 \mathbb{P}^0<\ldots<\mathbb{P}^j\ldots<\mathbb{P}^{n-1}
 \end{align*}    
where $\mathbb{P}^j$ is of degree $j \in \{0,1,\ldots,n-1\}$. Then in that case we have a $g\in\mathrm{Sp}_{n}(\mathbb{R})$ such that $g \mathbb{P}^j=\mathbb{D}_n^j\,(0\leq j <n)$.  
\end{remark}

\begin{definition}\label{parabolic general}
Let $\mathbb{P}$ be a boundary component of $\mathbb{D}_n$. Then the group $P(\mathbb{P})\subsetneq \mr{Sp}_n(\mb{R})$ defined by 
\begin{align*}
P(\mathbb{P})=\{g \in \mathrm{Sp}_{n}(\mathbb{R}) \,\vert\, g\mathbb{P}=\mathbb{P}\}
\end{align*}
is called the parabolic subgroup of $\mr{Sp}_n(\mb{R})$ associated to $\mathbb{P}$. 
\end{definition}
For the standard proper boundary components $\mathbb{P}=\mathbb{D}_n^j\,(0\leq j < n)$, the groups $P_j:=P(\mathbb{D}_n^j)$ has the structure (see \cite[page 21]{Namikawa})
\begin{align*}
P_j=\left\{\left.
\begin{pmatrix}
A^{\prime} & 0 & B^{\prime} & \star\\
\star & u & \star & \star\\
C^{\prime} & 0 & D^{\prime} & \star \\
0 & 0 & \star & u^{-t}
\end{pmatrix} \right\vert
\begin{pmatrix}
A^{\prime} & B^{\prime}\\
C^{\prime} & D^{\prime}
\end{pmatrix}\in \mr{Sp}_j(\mb{R}),\,u\in\mr{GL}_{n-j}(\mb{R})
\right\}.
\end{align*}
For a general boundary component $\mathbb{P}$ of $\overline{\mathbb{D}}_n$, realized as $\mathbb{P}=g \mathbb{D}_n^j$ for some $g\in\mathrm{Sp}_{n}(\mathbb{R})$ and some standard boundary component $\mathbb{D}_n^j$ $(0\leq j < n)$,the parabolic subgroup $P(\mathbb{P})$ associated to $\mathbb{P}$ can be obtained as $P(\mathbb{P})=gP_j\,g^{-1}$. 

\begin{definition}
A boundary component $\mathbb{P}$ of $\mathbb{D}_n$ is called rational if the parabolic subgroup $P(\mathbb{P})$ associated to it is defined over $\mb{Q}$. The set
\begin{align*}
\mb{D}_n^{\star}:=\bigsqcup_{\mathbb{P} \text{ rational}} \mathbb{P}\subsetneq \overline{\mathbb{D}}_n
\end{align*}
is called the rational closure of $\mb{D}_n$. 
\end{definition}

\begin{remark}\label{boundary scaling}
If $\mathbb{P}$ is a rational boundary component of $\mb{D}_n$, then there is a $\sigma\in\mr{Sp}_n(\mathbb{Z})$, such that $\sigma \mathbb{P}=\mathbb{D}_n^j$ for some standard boundary component $\mathbb{D}_n^j$ $(0\leq j < n)$. 
\end{remark}

The boundary components in the Siegel upper-half space $\mb{H}_n$ are obtained from the bounded realization $\mb{D}_n$ via the inverse Cayley transform. We denote the standard boundary components on $\mb{H}_n$ by $\mb{H}_n^j:=l^{-1}\mathbb{D}_n^j$ $(0\leq j< n)$. The rational closure $\mb{H}_n^{\star}$ of $\mb{H}_n$ is endowed with the cylindrical topology (see \cite[page 35]{Namikawa}). Under this topology, a sequence
\begin{align*}
Z(\nu):=\begin{pmatrix} Z_{1,1}(\nu) & Z_{1,2}(\nu)\\ Z_{1,2}(\nu) & Z_{2,2}(\nu)\end{pmatrix}\quad\quad (\nu\in\mb{N}_{>0},\,Z(\nu)\in\mb{H}_n,\,Z_{1,1}(\nu)\in \mb{H}_j)
\end{align*}
on $\mb{H}_n$ converges to a point $Z\in \mb{H}_n^j\cong\mb{H}_j$ in $\mb{H}_n^{\star}$ if and only if $Z_{1,1}(\nu)\rightarrow Z$ in $\mb{H}_j$ and 
$Y_{2,2}(\nu)-Y_{1,2}(\nu)Y_{2,2}(\nu)Y_{1,2}(\nu)\rightarrow \infty$ in $\mb{H}_{n-j}$. Under the assumption that $Y_{1,2}(\nu)$ is bounded, the latter condition reduces to $Y_{2,2}(\nu)\rightarrow \infty$.

In general, for any boundary component $\mb{P}$ of $\mb{H}_n$, one can show that there exists a one-parameter subgroup $w_{\mb{P}}\colon \mb{R}\rightarrow G$ such that
\begin{align*}
\lim_{t\rightarrow 0} w_{\mb{P}}(t)^{-1}O=O_{\mb{P}},
\end{align*}
where $O=i\mathbbm{1}_n$ is the base point of $\mb{H}_n$ and $O_{\mb{P}}$ is the base point of $\mb{P}$. For $\mb{P}=\mb{H}_n^j$, we denote $w_{\mb{P}}$ by $w_j$, which takes the form
\begin{align*}
w_j(t)=
\begin{pmatrix}
\mathbbm{1}_j & 0 & 0 & 0\\
0 & t\mathbbm{1}_{n-j} & 0 & 0\\
0 & 0 & \mathbbm{1}_j & 0\\
0 & 0 & 0 & t^{-1}\mathbbm{1}_{n-j}
\end{pmatrix}\qquad(t\in\mb{R}\setminus\{0\}). 
\end{align*}
It is easy to see that in the above sense, $w_j(t)^{-1}i\mathbbm{1}_n\rightarrow i\mathbbm{1}_j\in \mb{H}_n^j\cong\mb{H}_j$ as $t\rightarrow 0$. 

The parabolic subgroups $P(\mb{P})$ defined in \autoref{parabolic general} can be characterized in terms of $w_{\mb{P}}$ as
\begin{align*}
P(\mb{P})=\{g\in\mr{Sp}_n(\mb{R})\,\vert\,\lim_{t\rightarrow 0}w_{\mb{P}}(t)\,g \,w_{\mb{P}}(t)^{-1}<\infty\}. 
\end{align*}
We define
\begin{align*}
W(\mb{P})=\{g\in\mr{Sp}_n(\mb{R})\,\vert\,\lim_{t\rightarrow 0}w_{\mb{P}}(t)\,g \,w_{\mb{P}}(t)^{-1}=1\}
\end{align*}
\begin{remark}\label{arbitrary closeness}
Given $Z\in \mb{P}$, if for any sequence $Z(\nu)\,(\nu\in\mb{N}_{>0})$ in $\mb{H}_n$ such that $Z(\nu)\rightarrow Z$ in $\mb{H}_n^{\star}$, we have
$gZ(\nu)\rightarrow Z$ for some $g\in \mr{Sp}_n(\mb{R})$, then it is easy to see that $g\in W(\mb{P})$. 
\end{remark}
For $\mb{P}=\mb{H}_n^j$, we denote $W(\mb{P})$ by $W_j$, which can be shown to be (see \cite[page 21]{Namikawa})
\begin{align*}
W_j=\left\{ \left.
\begin{pmatrix}
\mathbbm{1}_j & 0 & 0 & Q\\
P^t & \mathbbm{1}_{n-j} & Q^t & B\\
0 & 0 & \mathbbm{1}_j & -P\\
0 & 0 & 0 & \mathbbm{1}_{n-j}
\end{pmatrix}
\right\vert
Q^t P+B=P^t Q+B^t
\right\}.
\end{align*}
Then, setting $P=L^t$, $Q=H^t$ and $B=LH^t+S_2$, we have
\begin{align}\label{W_j definition}
W_j=\bigg\{\begin{pmatrix}A & 0\\ 0 & A^{-t}\end{pmatrix}
\begin{pmatrix}\mathbbm{1}_n &  S\\ 0 & \mathbbm{1}_n\end{pmatrix}
\;\bigg\vert\; 
    A=\begin{pmatrix}\mathbbm{1}_j & 0\\ L & \mathbbm{1}_{n-j}\end{pmatrix},\,
    S=\begin{pmatrix}0 & H^t\\ H & S_2\end{pmatrix}\bigg\}\quad\quad(0\leq j< n), 
\end{align}
where $\,L,H\in\mb{R}^{(n-j)\times j}$ and $S_2\in\mb{R}^{(n-j)\times(n-j)},\,S_2=S_2^t$.

Next, for an arithmetic subgroup $\Gamma\subsetneq \mr{Sp}_n(\mb{R})$, consider the set $M:=\Gamma\backslash\mb{H}_n$.
\begin{theorem}\label{arithmetic boundary}
The quotient $M^{\star}:=\Gamma\backslash\mb{H}_n^{\star}$ endowed with the quotient topology, is a compact Hausdorff space. It contains $M$ as an open everywhere dense subset. $M^{\star}$ is the finite union of subspaces $M_j:=(\Gamma\cap P(\mb{P}_j))\backslash \mb{P}_j$, where $\mb{P}_j$ runs through a set of representatives of equivalence classes modulo $\Gamma$ of rational boundary components of $\mb{H}_n$. The closure of $M_j$ is the union of $M_j$ and the subspaces $M_k$ of $M_j$ of strictly smaller degree.
\end{theorem}
\begin{proof}
See \cite[Corollary 4.11]{Baily-Borel}. 
\end{proof}

The above compactification $M^{\star}$ of $M$ is called the Satake-Baily-Borel compactification of $M$. For $\Gamma=\Gamma_n=\mr{Sp}_n(\mb{Z})$, it takes the form
\begin{align}\label{modular boundary}
(\Gamma_n\backslash\mb{H}_n)^{\star}=\bigsqcup_{j=0}^n (\Gamma_n \cap P_j)\backslash \mb{H}_n^j\cong \bigsqcup_{j=0}^n \Gamma_j\backslash \mb{H}_j.
\end{align}

\begin{remark}\label{boundary equivalence}
By the remarks \ref{boundary chain} and \ref{boundary scaling}, the group $\Gamma_n$ acts transitively on the rational boundary components of $\mb{H}_n$.
Hence, for any arithmetic subgroup $\Gamma\subsetneq \mr{Sp}(n,\mb{R})$, we only need to consider $\mb{P}_j=\mb{H}_n^k\,(0\leq k < n)$ to fully describe the boundary of $M=\Gamma\backslash\mb{H}_n$. 
\end{remark}


\subsection{Siegel modular forms}
\begin{definition}
A function $f\colon\mathbb{H}_n\rightarrow \mathbb{C}$ is called a Siegel modular form of weight $\kappa$ and degree $n$ with respect to the Siegel modular group $\Gamma_n=\mr{Sp}_n(\mb{Z})$ if it satisfies the following conditions:
\begin{itemize}
\item[(i)] $f$ is holomorphic,
\item[(ii)] $f(\gamma Z)=\det(CZ+D)^{\kappa} f(Z)$ for all $\gamma=(\begin{smallmatrix}A&B\\C&D\end{smallmatrix})\in\Gamma$,
\item[(iii)] For every $Y_0>0$, the function $f$ is bounded in the region $Y\geq Y_0$. 
\end{itemize}
\end{definition}
We denote the space of all such functions by $\mathcal{M}_{\kappa}^n(\Gamma_n)$.
For all $S\in\mr{Sym}_n(\mathbb{Z})$, we have $(\begin{smallmatrix}\mathbbm{1}_n & S\\0 & \mathbbm{1}_n \end{smallmatrix})\in \Gamma_n$ . Then $f:\mathbb{H}_n\rightarrow \mathbb{C}$ is a holomorphic function satisfying $f(Z + S)=f(Z)$.
Therefore, $f$ has a Fourier expansion of the form
\begin{align*}
f(Z)=\sum_{\substack{T\in\mr{Sym}_n(\mb{Q})\\T\, \text{half-integral}}}a(T)\exp(2\pi i\:\tr(TZ)),
\end{align*} 
where $T=(t_{j,k})_{1\leq j,k\leq n}$ being half-integral implies that $t_{j,j},2t_{j,k}\in\mathbb{Z}\,(1\leq j\leq k\leq n)$. 
Also, since $U\in\mathrm{GL}(n,\mathbb{Z})$ implies $(\begin{smallmatrix} U^t&0\\0&U^{-1}\end{smallmatrix})\in\Gamma_n$, the function $f$ satisfies 
\begin{align*}
\det(U)^{\kappa} f(U^tZ U)=f(Z)\quad\quad(U\in\mr{GL}(n,\mb{Z})).
\end{align*} 
So, $U\in\mathrm{SL}(n,\mathbb{Z})$ implies that $f(U^tZU)=f(Z)$. One can show that a holomorphic function $f\colon\mathbb{H}_n\rightarrow\mathbb{C}$ satisfying these two transformation behaviours, i.e., $f(Z+S)=f(Z)$ for integral symmetric matrices $S$ and $f(U^tZU)=f(Z)$ for $U\in\mr{SL}(n,\mb{Z})$, under the assumption $n\geq2$, has a Fourier expansion of the form
\begin{align}\label{siegel fourier}
f(Z)=\sum_{\substack{T\in\mr{Sym}_n(\mb{Q}),\,T>0\\T\, \text{half-integral}}}a(T)\exp(2\pi i\:\tr(TZ)).
\end{align}
In particular, for some $Y_0>0$, the function $f$ is bounded in the region $Y\geq Y_0$. Thus, for $\Gamma=\Gamma_n$ and $n>1$, condition (iii) follows from conditions (i) and (ii). This is the so-called Koecher's principle. 

\begin{definition}
Let $f\colon\mathbb{H}_n\rightarrow \mathbb{C}$ be a function so that the limit 
\begin{align*}
\lim_{t\rightarrow \infty} f\begin{pmatrix}Z & 0\\ 0 & it\end{pmatrix}\quad\quad(Z\in\mb{H}_{n-1})
\end{align*}
exists. Then we obtain another function $\Phi(f)\colon \mb{H}_{n-1}\rightarrow \mb{C}$ defined by
\begin{align*}
\Phi(f)(Z):=\lim_{t\rightarrow \infty} f\begin{pmatrix}Z & 0\\ 0 & it\end{pmatrix}\quad\quad(Z\in\mb{H}_{n-1}).
\end{align*}
This operator $\Phi\colon \mathcal{M}_{\kappa}^n(\Gamma_n)\rightarrow \mathcal{M}_{\kappa}^{n-1}(\Gamma_{n-1})$ is called the Siegel $\Phi$-operator. 
\end{definition}
\begin{definition}\label{Siegel cusp form definition}
A Siegel modular form $f\in\mathcal{M}_{\kappa}^n(\Gamma_n)$ is called a Siegel cusp form if $\Phi(f)=0$. We denote the space of Siegel cusp forms by $\mathcal{S}_{\kappa}^n(\Gamma_n)$.
\end{definition}
\begin{proposition}\label{cusp form fourier coefficient}
A Siegel modular form $f\in\mathcal{M}_{\kappa}^n(\Gamma_n)$ is a Siegel cusp form if and only if in the Fourier expansion \eqref{siegel fourier}, $a(T)\neq 0$ implies that $T$ is positive definite.
\end{proposition}
\begin{proof}
See \cite[Hilfssatz 3.9]{Freitag}.
\end{proof}
\begin{proposition}
Let $f\in \mathcal{S}_{\kappa}^n(\Gamma_n)$ and let $c>0$. Then there exist positive numbers $c_1$ and $c_2$ such that 
\begin{align*}
\vert f(Z) \vert \leq c_1 \exp(-c_2 \tr(Y))
\end{align*}
for all $Z\in\mb{H}_n$, for which $Y$ is Minkowski reduced and $Y\geq c\mathbbm{1}_n$. 
\end{proposition}
\begin{proof}
See \cite[page 57]{Klingen}
\end{proof}

Next we define Siegel modular forms for arithmetic subgroups.
\begin{definition}
Let $\Gamma\subset\mathrm{Sp}_{n}(\mathbb{R})$ be an arithmetic subgroup and $\gamma_
{j}\in\mathrm{Sp}_{n}(\mathbb{Z})$ ($j=1,\ldots,h$) denote a set of representatives for the left cosets of $\Gamma\cap\mathrm{Sp}_{n}(\mathbb{Z})$
in $\mathrm{Sp}_{n}(\mathbb{Z})$. Then, a \emph{Siegel modular form of weight $\kappa$ and degree $n$ for $\Gamma$} is a function $f\colon\mathbb
{H}_{n}\longrightarrow\mathbb{C}$ satisfying the following conditions:
\begin{itemize}
\item[(i)] 
$f$ is holomorphic;
\item[(ii)] 
$f(\gamma Z)=\det(CZ+D)^{\kappa}f(Z)$ for all $\gamma=\big(\begin{smallmatrix}A&B\\C&D\end{smallmatrix}\big)\in\Gamma$;
\item[(iii)] 
given $Y_{0}\in\mathrm{Sym}_{n}(\mathbb{R})$ with $Y_{0}>0$, the quantities $\det(C_{j}Z+D_{j})^{-\kappa}f(\gamma_{j} Z)$ are bounded in the region 
$\{Z=X+iY\in\mathbb{H}_{n}\,\vert\,Y\geq Y_{0}\}$ for the set of representatives $\gamma_{j}=\big(\begin{smallmatrix}A_{j}&B_{j}\\C_{j}&D_{j}\end
{smallmatrix}\big)\in\mathrm{Sp}_{n}(\mathbb{Z})$ ($j=1,\ldots,h$).
\end{itemize}
\end{definition}
We denote the space of all such functions by $\mathcal{M}_{\kappa}^n(\Gamma)$. Just as in \autoref{Siegel cusp form definition}, a Siegel modular form $f\in \mathcal{M}_{\kappa}^n(\Gamma)$ with respect to the arithmetic subgroup $\Gamma$ is called a Siegel cusp forms with respect to $\Gamma$ if $\Phi(f)=0$. We denote the space of all such functions by $\mathcal{S}_{\kappa}^n(\Gamma)$.

For an arithmetic subgroup $\Gamma\subsetneq \mathrm{Sp}_n(\mathbb{R})$, define 
\begin{align*}
t(\Gamma):=\bigg\{S\in \mr{Sym}_n(\mb{R}) \;\bigg|\; \begin{pmatrix}\mathbbm{1}_n & S\\0 & \mathbbm{1}_n \end{pmatrix}\in\Gamma \bigg\}.
\end{align*}
Then $t(\Gamma)$ is commensurable to $t(\Gamma_n)$. Also, since $t(\Gamma_n)$ is commutative, $[t(\Gamma_n):t(\Gamma_n\cap\Gamma)]<\infty$ implies that there is an $\ell\in\mathbb{N}_{\geq 1}$ such that $\ell t(\Gamma)\subseteq t(\Gamma_n)$. Hence, $f\in \mathcal{M}_{\kappa}^n(\Gamma)$ satisfies 
$f(Z+\ell S)=f(Z)\,(S\in \mr{Sym}_n(\mb{R}))$. Therefore, the function $f_{\ell}$ defined by $f_{\ell}(Z)=f(\ell Z)$ satisfies $f_{\ell}(Z+S)=f_{\ell}(Z),\, (S\in \mr{Sym}_n(\mb{R}))$ and hence has a Fourier expansion 
\begin{align*}
f_{\ell}(Z)=f(\ell Z)=\sum_{\substack{T\in\mr{Sym}_n(\mb{Q}),\,T>0\\T\, \text{half-integral}}}a(T)\exp(2\pi i\:\tr(TZ)),
\end{align*}
whence, replacing $Z$ by $Z/\ell$, we have a Fourier expansion of $f$ of the form
\begin{align}\label{eq:fourier}
f(Z)=\sum_{\substack{T\in\mr{Sym}_n(\mb{Q}),\,T>0\\T\, \text{half-integral}}} a(T)\exp\bigg(\frac{2\pi i}{\ell}\:\tr(TZ)\bigg).
\end{align}
Just like in \autoref{cusp form fourier coefficient}, for a Siegel cusp form $f\in \mathcal{S}_{\kappa}^n(\Gamma)$ for which, the Fourier coefficients $a(T)$ are 0 unless $T$ is positive definite.  

\begin{proposition}
Let $f\in \mathcal{S}_{\kappa}^n(\Gamma)$ be a Siegel cusp form. Then, for the function
\begin{align*}
\varphi(Z):=\det(Y)^{\kappa/2} f(Z),
\end{align*}
$\vert\varphi(Z)\vert$ has a maximum in $\mb{H}_n$. 
\end{proposition}
\begin{proof}
See \cite[Bemerkung 6.10]{Freitag}. 
\end{proof}

Both $\mathcal{M}_{\kappa}^{n}(\Gamma)$ and $\mathcal{S}_{\kappa}^{n}(\Gamma)$ form finite-dimensional vector spaces over $\mathbb{C}$. The space $\mathcal{S}_{\kappa}^n(\Gamma)$, with the \emph{Petersson inner product} given by
\begin{align*}
\langle f, g\rangle:=\int_{\Gamma \backslash \mathbb{H}_{n}}(\det Y)^{\kappa} f(Z) \overline{g}(Z) \dif \mu_{n}(Z)\quad(f,g\in S_{\kappa}^{n}(\Gamma)),
\end{align*}
becomes a Hermitian inner product space.

\subsection{Siegel Maa\ss{} forms and $\Delta^{(\kappa)}$}\label{SEC Siegel Maass forms}

In order to derive sup-norm bounds for cusp forms $f\in  \mathcal{S}_{\kappa}^n(\Gamma)$, one introduces the function
\begin{align*}
\varphi(Z):=\det(Y)^{\kappa/2}f(Z)\quad (Z=X+iY\in\mb{H}_n,\,f\in  \mathcal{S}_{\kappa}^n(\Gamma))
\end{align*}
with transformation behaviour
\begin{equation}\label{ycusp}
\begin{aligned}
\varphi(\gamma Z)&=\det(\im(\gamma Z))^{\kappa/2}f(\gamma Z)=\bigg(\frac{\det(C Z+D)}{\det(C\overline{ Z}+D)}\bigg)^{\kappa/2}\varphi( Z),
\end{aligned}
\end{equation}
for all $\gamma=(\begin{smallmatrix}A&B\\C&D\end{smallmatrix})\in\Gamma$. We begin by defining a space $\mathcal{V}_{\kappa}^{n}(\Gamma)$ of real-analytic functions on $\mb{H}_n$ that transforms like \eqref{ycusp} with appropriate growth conditions.

\begin{definition}
\label{siegel-maass definition}
Let $\Gamma\subset\mathrm{Sp}_{n}(\mathbb{R})$ be a subgroup commensurable with $\mathrm{Sp}_{n}(\mathbb{Z})$, i.e., the intersection $\Gamma\cap
\mathrm{Sp}_{n}(\mathbb{Z})$ is a finite index subgroup of $\Gamma$ as well as of $\mathrm{Sp}_{n}(\mathbb{Z})$. We let $\gamma_{j}\in\mathrm{Sp}_{n}
(\mathbb{Z})$ ($j=1,\ldots,h$) denote a set of representatives for the left cosets of $\Gamma\cap\mathrm{Sp}_{n}(\mathbb{Z})$ in $\mathrm{Sp}_{n}(\mathbb
{Z})$. We then let $\mathcal{V}_{\kappa}^{n}(\Gamma)$ denote the space of all functions $\varphi\colon\mathbb{H}_{n}\rightarrow\mathbb{C}$ 
satisfying the following conditions:
\begin{itemize}
\item[(i)]
$\varphi$ is real-analytic;
\item[(ii)]
$\varphi(\gamma Z)=\det(CZ+D)^{\kappa/2}\det(C\overline{Z}+D)^{-\kappa/2}\varphi(Z)$ for all $\gamma=\big(\begin{smallmatrix}A&B\\C&D\end{smallmatrix}\big)\in
\Gamma$;
\item[(iii)]
given $Y_{0}\in\mathrm{Sym}_{n}(\mathbb{R})$ with $Y_{0}>0$, there exist $M\in\mathbb{R}_{>0}$ and $N\in\mathbb{N}$ such that the inequalities
\begin{align*}
\vert\det(C_{j}Z+D_{j})^{-\kappa/2}\det(C_{j}\overline{Z}+D_{j})^{\kappa/2}\varphi(\gamma_{j} Z)\vert\leq M\tr(Y)^{N} 
\end{align*}
hold in the region $\{Z=X+iY\in\mathbb{H}_{n}\,\vert\,Y\geq Y_{0}\}$ for the set of representatives $\gamma_{j}=\big(\begin{smallmatrix}A_{j}&B_{j}\\C_{j}&
D_{j}\end{smallmatrix}\big)\in\mathrm{Sp}_{n}(\mathbb{Z})$ ($j=1,\ldots,h$).
\end{itemize}
\end{definition}

\begin{remark}
\label{boundedness}
For $\varphi\in\mathcal{V}_{\kappa}^{n}(\Gamma)$, we set
\begin{align*}
\Vert{\varphi}\Vert^{2}:=\int\limits_{\Gamma\backslash\mathbb{H}_{n}}\vert\varphi(Z)\vert^{2}\,\mathrm{d}\mu_{n}(Z),
\end{align*}
whenever it is defined. In this way we obtain the Hilbert space
\begin{align*}
\mathcal{H}_{\kappa}^{n}(\Gamma):=\big\{\varphi\in\mathcal{V}_{\kappa}^{n}(\Gamma)\,\big\vert\,\Vert\varphi\Vert<\infty\big\}
\end{align*}
equipped with the inner product
\begin{align*}
\langle\varphi,\psi\rangle=\int\limits_{\Gamma\backslash\mathbb{H}_{n}}\varphi(Z)\overline{\psi}(Z)\,\mathrm{d}\mu_{n}(Z)\qquad(\varphi,
\psi\in\mathcal{H}_{\kappa}^{n}(\Gamma)).
\end{align*}
We note that in order to enable $\Vert\varphi\Vert<\infty$, the exponent $N\in\mathbb{N}$ in part (iii) of Definition~\ref{siegel-maass definition} has to be~$0$. 
\end{remark}

To compensate for not being holomorphic, the functions of the form $\varphi(Z)=\det (Y)^{\kappa/2}f(Z)\,(f\in \mathcal{S}_{\kappa}^n(\Gamma))$ satisfy the property of being eigenfunctions of a certain differential operator introduced by Maa\ss{} that is invariant under the transformation behaviour \eqref{ycusp}. 

Given $Z=X+iY\in\mathbb{H}_{n}$, we start by introducing the following symmetric $(n\times n)$-matrices of partial derivatives:
\begin{align*}
\mathrm{(i)}&\quad\bigg(\dpd{}{X}\bigg)_{j,k}:=\frac{1+\delta_{j,k}}{2}\dpd{}{x_{j,k}}, \\
\mathrm{(ii)}&\quad\bigg(\dpd{}{Y}\bigg)_{j,k}:=\frac{1+\delta_{j,k}}{2}\dpd{}{y_{j,k}}, \\
\mathrm{(iii)}&\quad\dpd{}{Z}:=\frac{1}{2}\bigg(\dpd{}{X}-i\dpd{}{Y}\bigg), \\
\mathrm{(iv)}&\quad\dpd{}{\overline{Z}}:=\frac{1}{2}\bigg(\dpd{}{X}+i\dpd{}{Y}\bigg),
\end{align*}
where $\delta_{j,k}$ is the Kornecker delta symbol.

\begin{definition}
Given a positive integer $\kappa$, the differential operator $\Delta^{(\kappa)}$ given by
\begin{align*}
\Delta^{(\kappa)}=\tr\bigg(Y\bigg(\bigg(Y\dpd{}{X}\bigg)^{t}\dpd{}{X}+\bigg(Y\dpd{}{Y}\bigg)^{t}\dpd{}{Y}\bigg)-i\kappa Y\dpd{}{X}\bigg)
\end{align*}
acting on smooth complex valued functions on~$\mathbb{H}_{n}$ is called the Siegel--Maa\ss{} Laplacian of weight $\kappa$.
\end{definition}
By its invariance under the transformation behaviour \eqref{ycusp}, the operator $\Delta^{(\kappa)}$ acts on the Hilbert space $\mathcal{H}_{\kappa}^{n}(\Gamma)$ (see \cite[Remark 4.6]{siegel--maass}).

\begin{definition}
Let $\Gamma\subset\mathrm{Sp}_{n}(\mathbb{R})$ be a subgroup commensurable with $\mathrm{Sp}_{n}(\mathbb{Z})$. The elements of the Hilbert space 
$\mathcal{H}_{\kappa}^{n}(\Gamma)$ are called \emph{automorphic forms of weight $\kappa$ and degree $n$ for $\Gamma$}. Moreover, if 
$\varphi\in\mathcal{H}_{\kappa}^{n}(\Gamma)$ is an eigenform of $\Delta^{(\kappa)}$, it is called a \emph{Siegel--Maa\ss{} form of weight $\kappa$ and degree $n$ for $\Gamma$}.
\end{definition}

\begin{theorem}
\label{kernel connection}
Let $\Gamma\subset\mathrm{Sp}_{n}(\mathbb{R})$ be a subgroup commensurable with $\mathrm{Sp}_{n}(\mathbb{Z})$ and let $\varphi\in\mathcal{H}_{\kappa}^{n}
(\Gamma)$ be a Siegel--Maa\ss{} form of weight $\kappa$ and degree $n$ for $\Gamma$. Then, if $\varphi$ is an eigenform of $\Delta^{(\kappa)}$ with eigenvalue $\lambda$, 
we have $\lambda\in\mathbb{R}$ and 
\begin{align*}
\lambda\geq \frac{n\kappa}{4}(n-\kappa+1),
\end{align*}
with equality attained if and only if the function $\varphi$ is of the form $\varphi(Z)=\det(Y)^{\kappa/2}f(Z)$ for some Siegel cusp form $f\in\mathcal{S}_{\kappa}^{n}(\Gamma)$ 
of weight $\kappa$ and degree $n$ for $\Gamma$. In other words, there is an isomorphism
\begin{align*}
\mathcal{S}_{\kappa}^{n}(\Gamma)\cong\ker\bigg(\Delta^{(\kappa)}+\frac{n\kappa}{4}(n-\kappa+1)\mathrm{id}\bigg)
\end{align*}
of $\mathbb{C}$-vector spaces, induced by the assignment $f\mapsto\det(Y)^{\kappa/2}f$.  
\end{theorem}
\begin{proof}
See \cite[Corollary 5.4]{siegel--maass}
\end{proof}


\section{Construction of the heat kernel}

To use \eqref{limit-inequality} to obtain sup-norm bounds for the quantity $S_{\kappa}^{\Gamma}(Z)$, we need to obtain a somewhat explicit form for the heat kernel $K_{t}^{(\kappa)}$ corresponding to the Siegel--Maa\ss{} Laplacian $\Delta^{(\kappa)}$ on $\mb{H}_n$. In the theory of harmonic analysis on symmetric spaces, there is a standard way of obtaining the heat kernel $K_t$ corresponding to the Laplace--Beltrami operator $\Delta$ from the spherical function on the given symmetric space, from which, one can use a weight-correction technique to obtain the heat kernel $K_{t}^{(\kappa)}$ corresponding to $\Delta^{(\kappa)}$. Thus, the problem of obtaining a somewhat explicit form for $K_{t}^{(\kappa)}$ on $\mb{H}_n$ translates to the problem of obtaining a somewhat explicit form for the spherical function $\phi_{\lambda}$ on $\mathrm{Sp}_n(\mb{R})$. It is difficult to do it directly. Instead, we wield a technique developed by Flensted-Jensen in \cite{FLENSTED-JENSEN} of obtaining the spherical function $\phi_{\lambda}$ on a real semisimple group by reducing it to obtaining the spherical function $\Phi_{\Lambda}$ on the corresponding complex semisimple group, which is much simpler.

In the first two subsections, we briefly recall the general theory of spherical functions on a real semisimple group via the Flensted-Jensen reduction. The general reference for these subsections are \cite{Helgason2} and \cite{FLENSTED-JENSEN}. In the third subsection, we recall from \cite{Gangolli1} and \cite{Gangolli2}, the general procedure of construction of the heat kernel on a real semisimple group by spherical inversion. Then in the next two subsections we implement this procedure for the special case of the symplectic group. Finally, in the last subsection, we apply a weight-correction procedure to obtain the heat kernel on $\mb{H}_n$ corresponding to $\Delta^{(\kappa)}$. 

But first, we need to fix some basic notations for this section. 


\subsection{Notation}\label{SEC notation}
 
Let $\mf{g}$ be a complex semisimple Lie algebra with a Cartan decomposition $\mf{g}=\mf{u}+\mf{p}$ corresponding to a Cartan involution $\theta$ on $\mf{g}$.
Let $\mf{g}_{0}$ be a non-compact real form of $\mf{g}$. The Cartan involution $\theta$ on $\mf{g}$ restricts to the Cartan involution $\theta_0$ on $\mf{g}_0$.
Let $\mf{g}_{0}=\mf{k}_{0}+\mf{p}_{0}$ be the Cartan decomposition of $\mf{g}_0$ corresponding to the Cartan involution $\theta_0$ on $\mf{g}_0$. Then $\mf{u}=\mf{k}_{0}+i \mf{p}_{0}$ is a compact real form of $\mf{g}$ and $\mf{k}=\mf{k}_{0}+i\mf{k}_{0}$ is a complex subalgebra of $\mf{g}$. Denote by $\mf{a}$ and $\mf{a}_0$ the maximal abelian subspaces of $\mf{p}$ and 
$\mf{p}_0$, respectively. We note here that $\mf{k}_0$, $\mf{p}_0$ and $\mf{a}_0$ are related to $\mf{k}$, $\mf{u}$, $\mf{p}$ and $\mf{a}$ via $\mf{k}_0=\mf{k}\cap\mf{g}_0=\mf{u}\cap\mf{g}_0$, $\mf{p}_0=\mf{p}\cap\mf{g}_0$ and $\mf{a}_0=\mf{a}\cap\mf{g}_0$. 

The Killing form $B_{0}$ on $\mf{g}_{0}$ is just the restriction of the complex Killing form $B^{\prime}$ of $\mf{g}$, whereas the Killing form $B$ of $\mf{g}$ as a real Lie algebra, is $2 B^{\prime}$. This means that the Euclidean structures on $\mf{a}_{0}$, induced by $B_{0}$ and $B$ are different. Denote by $\langle\cdot,\cdot\rangle_{0}$ and $\vert\vert\cdot\vert\vert_0$ the scalar product and norm induced by $B_{0}$ on $\mf{a}_{0}$ as well as $\langle\cdot,\cdot\rangle$ and $\vert\vert\cdot\vert\vert$ for the scalar product and norm induced by $B$ on $\mf{p}$. So in particular
\begin{align*}
\|H\|_{0}^{2}=\frac{1}{2}\|H\|^{2} \quad \text { for all } \quad H \in \mf{a}_{0}.
\end{align*}
By the Killing form identification of $\mf{a}_{0}$ and $\mf{a}$ with their duals, $\mf{a}_{0}\spcheck$ is embedded in $\mf{a}\spcheck$. The Euclidean structures on the spaces $\mf{a}$ and $\mf{p}$ induce Euclidean structures on the dual spaces $\mf{a}\spcheck$ and $\mf{p}\spcheck$, respectively, by duality. Denote by $\langle\cdot,\cdot\rangle_{0}$ and $\vert\vert\cdot\vert\vert_0$ the induced scalar product and norm on $\mf{a}_{0}\spcheck$ as well as by $\langle\cdot,\cdot\rangle$ and $\vert\vert\cdot\vert\vert$ the induced scalar product and norm on $\mf{p}\spcheck$. So in particular
\begin{align*}
\|\lambda\|_{0}^{2}=\frac{1}{2}\|\lambda\|^{2} \quad \text { for all } \quad \lambda \in \mf{a}_{0}\spcheck.
\end{align*}

Let $\Delta$ be the root system of the pair $(\mf{g},\mf{a})$, by which we mean that $\Delta$ is the set of restricted roots for the real Lie algebra $\mf{g}$ with respect to the maximal abelian subalgebra $\mf{a}$. Then each root space $\mf{g}^{\alpha}\,(\alpha \in \Delta)$ has dimension $m_{\alpha}=2$. Let $\Delta_{0}$ be the restricted root system of the pair $(\mf{g}_0,\mf{a}_0)$. 
Then $\Delta_{0}=\{\alpha\vert_{\mf{a}_{0}}\; \vert\, \alpha \in \Delta, \alpha\vert_{\mf{a}_{0}} \neq 0\}$. Let $W$ and $W_0$ be the Weyl groups corresponding to the restricted root systems $\Delta$ and $\Delta_0$.

Let $\Delta^{+}$ and $\Delta_{0}^{+}$ be choices of positive restricted roots in the restricted root systems $\Delta(\mf{g},\mf{a})$ and $\Delta_0(\mf{g}_0,\mf{a}_0)$, respectively. Let $\mf{a}^{+}$ and $\mf{a}_{0}^{+}$ be the corresponding choices of positive Weyl chambers in $\mf{a}$ and $\mf{a}_0$, respectively. Let $\rho$ and $\rho_0$ denote the half-sums (with multiplicity) 
\begin{align*}
\rho=\frac{1}{2} \sum_{\alpha \in \Delta^{+}} m_{\alpha}\, \alpha=\sum_{\alpha \in \Delta^{+}} \alpha \quad\text{and}\quad
\rho_{0}=\frac{1}{2} \sum_{\lambda \in \Delta_{0}^{+}} m_{\alpha} \alpha
\end{align*}
of the positive restricted roots for $(\mf{g}, \mf{a})$ and $(\mf{g}_{0}, \mf{a}_{0}),$ respectively. Similarly, let $\pi$ and $\pi_0$ denote the products of the indivisible positive restricted roots
\begin{align}\label{pi-def}
\pi (\lambda)=\prod_{\alpha \in \Delta^{+}}\langle\alpha, \lambda\rangle\quad(\lambda\in \mf{a}\spcheck) \quad\text{and}\quad \pi_{0}(\lambda)=\prod_{\alpha \in \Delta_{0}^{+}} \langle\alpha, \lambda \rangle \quad (\lambda\in\mf{a}_0\spcheck),
\end{align}
respectively.

Denoting $\mf{n}=\sum_{\alpha\in\Delta^{+}} \mf{g}^{\alpha}$ and $\mf{n}_{0}=\mf{n} \cap \mf{g}_{0}$, the algebras $\mf{g}$ and $\mf{g}_0$ have Iwasawa decompositions $\mf{g}=\mf{u}+\mf{a}+\mf{n}$ and $\mf{g}_{0}=\mf{k}_{0}+\mf{a}_{0}+\mf{n}_{0}$, respectively.

Let $G$ be a Lie group with Lie algebra $\mf{g}$, and let $K$, $U$, $A$, $N$, $G_0$, $K_0$, $A_0$ and $N_0$ be the analytic subgroups corresponding to $\mf{k}$, $\mf{u}$, $\mf{a}$, $\mf{n}$, $\mf{g}_0$, $\mf{k}_0$, $\mf{a}_0$ and $\mf{n}_0$. Corresponding to the algebra level Iwasawa decompositions $\mf{g}=\mf{u}+\mf{a}+\mf{n}$ and $\mf{g}_{0}=\mf{k}_{0}+\mf{a}_{0}+\mf{n}_{0}$ of $\mf{g}$ and $\mf{g}_0$, respectively, the groups $G$ and $G_0$ have the group level Iwasawa decompositions $G=UAN$ and $G_0=K_0 A_0 N_0$, so that the mapping $(u,a,n)\mapsto uan$ is a diffeomorphism of $U\times A \times N$ onto $G$ and $K_0\times A_0 \times N_0$ onto $G_0$. Let for $g\in G$, $H(g)\in \mf{a}$ be determined by $g\in U\exp(H(g))N$. If $g\in G_0$, then we have $H(g)\in\mf{a}_0$. 

The group $G_0$ also has the polar decomposition $G_0=K_0 A_0 K_0$, by which,
for each $g\in G_0$, there is an $a\in A_0$ such that $g\in K_0\,a K_0$. If, for a particular choice of a positive Weyl chamber $\mf{a}_0^{+}$, we restrict ourselves to $\overline{A_0^{+}}=\exp(\overline{\mf{a}_0^{+}})$, then for each $g\in G_0$, the choice of $a\in \overline{A_0^{+}}$ such that $g\in K_0\,a K_0$, is unique. It can be shown that the set $C^{\infty}(K_0\backslash G_0/K_0)$ of $K_0$-bi-invariant $C^{\infty}$-functions on $G_0$, via restriction to $A_0$, is in bijective correspondence with $C_{W_0}^{\infty}(A_0)$, the set of $W_0$-invariant $C^{\infty}$-functions on $A_0$. 

Similarly, for the complex group $G$, we have the polar decomposition $G=U \overline{A^{+}} U$. Furthermore, $G$ also has the Mostow decomposition $G=U \overline{A_0^{+}} K$, by which, for each $g\in G$, there is a unique $a\in \overline{A_0^{+}}$ such that $g\in UaK$. The set $C^{\infty}(U\backslash G/K)$ is in bijective correspondence, via restriction to $A_0$, with the set $C_{W_0}^{\infty}(A_0)$ of $W_0$-invariant $C^{\infty}$-functions on $A_0$ (see \cite[Theorem 4.1]{FLENSTED-JENSEN}).


\subsection{Spherical functions on $G_0/K_0$}\label{SEC spherical general}

Consider the Riemannian globally symmetric space $G_0/K_0$. Let $\pi\colon G_0\rightarrow G_0/K_0$ denote the natural mapping of $G_0$ onto $G_0/K_0$ and $o\in G_0/K_0$ denote the point $o=\pi(e)$, where $e\in G_0$ is the neutral element of $G_0$. If $f$ is any function on $G_0/K_0$, let $\widetilde{f}$ denote the function $\widetilde{f}=f \circ \pi$ on $G$. Let $D(G_0)$ denote the set of all left-invariant differential operators on $G_0$, $D_{K_0}(G_0)\subsetneq D(G_0)$ the subspace of $D(G_0)$ containing left-invariant differential operators on $G_0$ which are also right-invariant under $K_0$ and $D(G_0/K_0)$ the algebra of differential operators on $G_0/K_0$ invariant under all left translations of $G_0/K_0$ by $G_0$.
 
\begin{definition} 
A complex-valued function $\phi\in C^{\infty}(G_0/K_0)$ on $G_0/K_0$ is called a spherical function on $G_0/K_0$ if it satisfies the following properties:
\begin{enumerate}[(i)]
\item $\phi(o)=1$,
\item $D \phi=\lambda_{D}\, \phi$ for each $D \in D(G_0 / K_0)$, where $\lambda_{D}$ is a complex number,
\item $\phi(k_0 g K_0)=\phi(gK_0)$ for all $g\in G_0$ and $k_0 \in K_0$.
\end{enumerate}
The function $\widetilde{\phi}=\phi \circ \pi$ on $G_0$ is called a spherical function on $G_0$ if and only if $\phi$ is a spherical function on $G_0/K_0$. 
\end{definition}

From the above definition it is easy to see that a spherical function $\widetilde{\phi}$ on $G_0$ is characterized by the following properties:
\begin{enumerate}[(i)]
    \item $\widetilde{\phi}(e)=1$,
    \item $D \widetilde{\phi}=\lambda_{D}\, \widetilde\phi$ for each $D \in D_{K_0}(G_0)$, where $\lambda_{D}$ is a complex number,
    \item  $\widetilde{\phi}(k_0 g k_0^{\prime})=\widetilde{\phi}(g)$ for all $g \in G_0$ and all $k_0, k_0^{\prime} \in K_0$. 
\end{enumerate}
As noted in the last subsection, due to the bi-invariance of $\widetilde{\phi}$ under $K_0$, it suffices to know $\widetilde{\phi}\in C^{\infty}(K_0\backslash G_0 /K_0)$ on the Weyl chamber $A_0^{+}=\exp (\mf{a}_0^{+})$.

\begin{remark}
As the notion of spherical functions on the group $G_0$ is equivalent to that on the symmetric space $G_0/K_0$, for convenience, we denote the spherical functions on both $G_0$ and $G_0/K_0$ by $\phi$.  
\end{remark}

For a symmetric space  $G_0/K_0$ of non-compact type, Harish-Chandra \cite{Harish1} gave the following characterization of spherical functions on $G_0/K_0$ in terms of an integral.

\begin{theorem}\label{harish-integral}
Let $G_0$ be a connected semisimple Lie group with finite centre and $K_0$ a maximal compact subgroup of $G_0$. Then, as $\lambda$ runs through $(\mf{a}^{\mb{C}}_0)\spcheck$, the functions
\begin{align}\label{harish-spherical 2}
\phi_{\lambda}(g)=\int_{K_0} \exp\big((i \lambda-\rho_0)(H(g k_0))\big) \dif \mu(k_0) \quad (g \in G_0),
\end{align}
where $\dif\mu(k_0)$ denotes the Haar measure on $K_0$, exhaust the class of spherical functions on $G_0$. Moreover, two such functions
$\phi_{\mu}$ and $\phi_{\lambda}$ are identical if and only if $\mu=\sigma\lambda$ for some $\sigma$ in the Weyl group $W_0$.
\end{theorem}

\begin{proof}
See \cite[Chapter IV, Theorem 4.3]{Helgason2}.
\end{proof}

\begin{lemma}\label{helgason lemma}
Let $G_0$ be a connected semisimple Lie group with finite centre and $K_0$ a maximal compact subgroup of $G_0$. Then, for $a\in A_0^+$, we have
\begin{align*}
    \phi_0(a)\geq \exp(-\rho_0(\log(a))).
\end{align*}
\end{lemma}

\begin{proof}
 Given the positive Weyl chamber $\mf{a}_0^+$, let $^+\mf{a}_0$ denote the dual cone defined by
 \begin{align*}
     ^+\mf{a}_0:=\{H\in\mf{a}_0\,\vert\,B(H,H^{\prime})>0,\,\forall\, H^{\prime}\in \mf{a}_0^+\},
 \end{align*}
 and let $\overline{^+\mf{a}_0}$ denote its closure. Then, by \cite[Chapter IV, Lemma 6.5]{Helgason2}, for $a\in A_0^+$, we have
 \begin{align*}
     \log(a)-H(a k_0)\in \overline{^+\mf{a}}\quad(k_0\in K_0),
 \end{align*}
 which implies that
 \begin{align}\label{helgaq}
     \rho_0(\log(a))\geq\rho_0(H(ak_0))\quad(k_0\in K_0). 
 \end{align}
 Then, by Harish-Chandra's characterization of $\phi_{\lambda}\,(\lambda\in\mf{a}_0\spcheck)$ in terms of the integral
 \begin{align*}
\phi_{\lambda}(g)=\int_{K_0} \exp\big((i \lambda-\rho_0)(H(g k_0))\big) \dif \mu(k_0) \quad (g \in G_0)
\end{align*}
given in equation \eqref{harish-spherical 2}, we have

\begin{align*}
    \phi_{0}(a)& =\int_{K_0} \exp\big(-\rho_0(H(a k_0))\big) \dif \mu(k_0)\\
    & \geq \exp\big(-\rho_0(\log(a))\big)\int_{K} \dif \mu(k_0)
     =\exp\big(-\rho_0(\log(a))\big)  
\quad (a\in A_0^+),
\end{align*}
thereby proving the lemma. 
\end{proof}

\begin{remark}\label{remark c-function}
In \cite{Harish1}, Harish-Chandra gave a series expansion of $\phi_{\lambda}$ with leading coefficients $c(\sigma \lambda)\,(\sigma\in W_0)$. This function, called \emph{Harish-Chandra's $c$-function}, features prominently in the theory of spherical transforms and 
was explicitly determined by Gindikin and Karpelevic as a meromorphic function on $(\mf{a}^{\mb{C}}_0)\spcheck$.
In particular, for $G_0=\mr{Sp}_n(\mb{R})$, corresponding to the vector $\lambda=\lambda_1e_1+\ldots+\lambda_ne_n\in\mf{a}\spcheck\cong\mb{R}^n$ (see \autoref{SEC Spherical function on Hn}), Bhanu Murti \cite{bhanu} showed that
\begin{align}\label{bhanu}
\begin{aligned}
\abs{c(\lambda)}^{-2} =&\frac{1}{\pi^{n^{2} / 2}} \prod_{1\leq j\leq n}\frac{\lambda_{j}}{2} \operatorname{th} \Big(\frac{\lambda_{j}}{2}\pi\Big) \prod_{1\leq j<k \leq n}\frac{\lambda_{j}+\lambda_{k}}{2} \operatorname{th} \Big(\frac{\lambda_{j}+\lambda_{k}}{2} \pi\Big) \times\\ 
&\quad\times\prod_{1 \leq j<k \leq n} \frac{\lambda_{j}-\lambda_{k}}{2} \operatorname{th} \Big(\frac{\lambda_{j}-\lambda_{k}}{2} \pi\Big).
\end{aligned}
\end{align}
\end{remark}

\begin{definition}
Let $f$ be a smooth function on $G_0$ which is bi-invariant under $K_0$. The function $\widehat{f}\colon \mf{a}_0\spcheck\rightarrow\mb{C}$ defined by
\begin{align}\label{spherical transform}
\widehat{f}(\lambda)=\int_{G} f(g) \phi_{-\lambda}(g) \dif \mu(g)\quad(\lambda\in \mf{a}_0\spcheck),
\end{align}
is called the \emph{spherical transform} of $f$ at $\lambda\in\mf{a}_0\spcheck$. 
\end{definition}

The next theorem states the crucial inversion formula for the spherical transform. 
\begin{theorem}\label{spherical inversion theorem}
For $g=k_1\exp(H(g))k_2\in G_0\;(k_1,k_2\in K_0,\,H(g)\in\mf{a}_0)$, define $\vert g\vert:=B(H(g),H(g))$. Then the spherical $L^2$-Schwartz space $\mathscr{C}(K_0\backslash G_0/K_0)$ is the space of all functions $f\in C^{\infty}(K_0\backslash G_0 /K_0)$ such that for all $N\in\mb{N}_{\geq 1}$ and $D\in D(G_0)$,
\begin{align*}
   \sup\limits_{g\in G_0}\,(1+\vert g\vert)^N \vert Df(g)\vert\, \phi_0(g)^{-1}<\infty.
\end{align*}
Let $\mathscr{S}(\mf{a}_0\spcheck)$ denote the usual Schwartz space on $\mf{a}_0\spcheck$ of rapidly decreasing smooth functions and $\mathscr{S}_{W_0}(\mf{a}_0\spcheck)$ be the subspace of $W_0$-invariant elements. Further, let $\mathscr{H}^{R} (\mf{a}_0\spcheck)\,(R\in\mb{R}_{>0})$ denote the set of functions $f$ on $\mf{a}_0\spcheck$ satisfying the criterion that for each $N\in\mb{N}_{\geq 1}$, there exists a constant $C_N\in\mb{R}_{>0}$, for which the function $f$ satisfies the condition
\begin{align*}
|f(\lambda)| \leq C_{N}(1+|\lambda|)^{-N} \exp(R\vert\operatorname{Im} (\lambda)\vert) \quad (\lambda \in \mf{a}_0\spcheck)
\end{align*}
and let $\mathscr{H}(\mf{a}_0\spcheck)=\bigcup_{R>0} \mathscr{H}^{R}(\mf{a}_0\spcheck)$. Let $\mathscr{H}_{W_0}(\mf{a}_0\spcheck)$ and $\mathscr{H}_{W_0}^{R}(\mf{a}_0\spcheck)$ denote the respective subspaces of $W_0$-invariant elements. Then the following assertions hold:
\begin{enumerate}[(i)]
\item The spherical transform given by the assignment $f \mapsto\widehat{f}$ induces a bijection of $\mathscr{C}(K_0\backslash G_0/K_0)$ onto $\mathscr{S}_{W_0}(\mf{a}_0\spcheck)$.

\item Restriction of the domain of the above transform to $C_{c}^{\infty}(K_0\backslash G_0/K_0)\subsetneq \mathscr{C}(K_0\backslash G_0/K_0)$ restricts the bijection onto the subspace $\mathscr{H}_{W_0}(\mf{a}_0\spcheck)\subsetneq\mathscr{S}_{W_0}(\mf{a}_0\spcheck)$. 
\item For $f\in \mathscr{C}(K_0\backslash G_0/K_0)$ and $\widehat{f}\in \mathscr{S}_{W_0}(\mf{a}_0\spcheck)$, we have the formula, called the inverse spherical transform, given by
\begin{align*}
 f(g)=\int_{\mf{a}_0\spcheck} \widehat{f}(\lambda) \phi_{\lambda}(g)\abs{{c}(\lambda)}^{-2} \dif \lambda \quad (g \in G_0),
\end{align*}
where $\dif\lambda$ denotes the Euclidean measure on $\mf{a}_0\spcheck/W_0$.
\item As $\mathscr{C}(K_0\backslash G_0/K_0)$ is dense in $L^2(K_0\backslash G_0/K_0)$ and its image $\mathscr{S}_{W_0}(\mf{a}_0\spcheck)$ is dense in \\ $L^{2}(\mathfrak{a}_0\spcheck / W_0,\,\abs{{c}(\lambda)}^{-2} \dif \lambda)$, the spherical transform given by the assignment $f \mapsto \widehat{f}$\\  extends by continuity to an isometry  of $L^2(K_0\backslash G_0/K_0)$ onto $L^{2}(\mf{a}_0\spcheck / W_0,\,\abs{{c}(\lambda)}^{-2} \dif \lambda)$, thereby giving the equality
\begin{align*}
 \int_{G_0}|f(g)|^{2} \dif \mu(g)=\int_{\mf{a}_0\spcheck}|\widehat{f}(\lambda)|^{2}|c(\lambda)|^{-2} \dif \lambda.
 \end{align*}
\end{enumerate}
\end{theorem}
\begin{proof}
See \cite[Chapter IV, Section 7]{Helgason2} and \cite[Sections 5 and 6]{Gangolli3}. 
\end{proof}

One should note here that everything said above concerning spherical functions
on real semisimple groups $G_0$ with respect to $K_0$ applies equally to spherical functions on complex semisimple groups $G$ with respect to $U$, which is just the special case where $\mf{g}$ has a complex structure. However, using Harishchandra's series expansion of spherical function in \autoref{remark c-function}, in case of a complex Lie group $G$, the spherical function on $G$ takes the following much simpler form.
\begin{theorem}\label{complex spherical theorem}
Let $G$ be a complex Lie group. Then the spherical function of $G$ corresponding to 
$\lambda\in \mf{a}\spcheck$ is given by
\begin{align*}
\phi_{\lambda}(a)=\frac{\pi(\rho)}{\pi(i \lambda)} \frac{\sum_{\sigma\in W}
\det (\sigma) \exp(i \sigma \lambda(\log (a)))}{\sum_{\sigma \in W}\det (\sigma)\, \exp(\sigma \rho(\log (a)))} \quad (a \in A^+),
\end{align*}
where $\pi(\lambda)=\prod\limits_{\alpha \in \Delta^{+}}\langle\alpha, \lambda\rangle $. Moreover, the $c$-function in this case is given by
\begin{align*}
c(\lambda)=\pi(\rho) / \pi(i \lambda).
\end{align*}
\end{theorem}
\begin{proof}
See \cite[Chapter IV, Section 5]{Helgason2}
\end{proof}


\subsection{Flensted-Jensen reduction}\label{SEC Flensted-Jensen reduction}

Now consider $\mf{g}$ as a Lie algebra over $\mb{R}$. Let $D_{R}(K \backslash G)$ denote the set of right-invariant differential operators on the coset space $K \backslash G=\{K g\,\vert\, g \in G\}$.

Let
\begin{align*}
C^{\infty}(K \backslash G / U)=\{\phi \in C^{\infty} (G)\,\vert\, \phi(k g u) = \phi(g)\}.
\end{align*}
The main result in \cite{FLENSTED-JENSEN} is the following theorem that enables us to lift many questions related to the analysis of spherical functions on a real group $G_{0}$, to analogous questions concerning the spherical functions on the corresponding complex group $G$. 
\begin{theorem}\label{eta correspondence}
There is a one-to-one correspondence induced by $\phi \mapsto \phi^{\eta}$ between the set of spherical functions $\phi$ on $G_{0} / K_{0}$ and the set of functions $\psi=\phi^{\eta}$ on $G$ satisfying
\begin{enumerate}[(i)]
    \item $\psi(e)=1$,
    \item $D \psi=\lambda_{D} \psi$ for all $D \in D_{R}(K \backslash G)$, where $\lambda_{D}$ is a complex number,
    \item $\psi \in C^{\infty}(K \backslash G / U)$,
\end{enumerate}
such that
\begin{align*}
\phi(g \Theta(g)^{-1})=\phi^{\eta}(g) \quad (g \in G_{0}),
\end{align*}
where $\Theta\colon G_0\rightarrow G_0$ is the involutive automorphism of $G_0$ such that $(\dif\Theta)_e=\theta$.
\end{theorem}
\begin{proof}
See \cite[Section 5]{FLENSTED-JENSEN}. 
\end{proof}

This allows us to identify $C^{\infty}(K_0 \backslash G / K_0)$ with $C^{\infty}(K \backslash G / U)$ and write $\phi$ instead of $\phi^{\eta}$. Let $\phi_{\lambda}\,(\lambda \in \mf{a}_{0}\spcheck)$ denote the spherical functions on $G_{0} / K_{0}$ and $\Phi_{\Lambda}$ $(\Lambda \in \mf{a}\spcheck)$ the spherical functions on $G / U$. If $\lambda \in \mf{a}_{0}\spcheck$, define $\Lambda \in \mf{a}\spcheck$ by
\begin{align*}
\Lambda+i \rho=2(\lambda+i \rho_{0}).
\end{align*}
Then we have
\begin{align*}
\Phi_{\Lambda}(g)=\int_{U} \phi_{\lambda}(u g) \dif \mu(u) \quad (g \in G),
\end{align*}
where $\dif \mu(u)$ is the normalized Haar measure on $U$. Under this setup, the following theorem enables us to calculate spherical functions on non-compact real Lie groups from spherical functions on the corresponding complex Lie group via an integral transform.

\begin{theorem}\label{Flensted reduction}
Let $\mf{g}_{0}$ be a normal real form of the complex Lie algebra $\mf{g}$. Assume that the Haar measure $\dif \mu(k)$ on $K$ is normalized such that on compact groups the total mass is $1$ and on non-compact, $d$-dimensional spaces the measure is $(2\pi)^{-d/2}$ times the volume element so that the Euclidean Fourier transform is an isometry. Then, the spherical functions $\phi_{\lambda}$ on $G_0/K_0$ and $\Phi_{\Lambda}$ on $G/U$ are related by the equation
\begin{align}\label{Flensted formula}
\phi_{\lambda}(g \Theta(g)^{-1})=|c(\lambda)|^{2}|\pi_{0}(\lambda)|^{2} \int_{K} \Phi_{2 \lambda}(k g) \dif \mu(k) \quad (\lambda \in \mathfrak{a}_{0}\spcheck).
\end{align}
In particular, we have
\begin{align*}
|c(\lambda)|^{-2}=|\pi_{0}(\lambda)|^{2} \int_{K} \Phi_{2 \lambda}(k) \dif \mu(k) \quad (\lambda \in \mathfrak{a}_{0}\spcheck).
\end{align*}
\end{theorem}
\begin{proof}
See \cite[Section 7]{FLENSTED-JENSEN}.
\end{proof}


\subsection{Heat kernel on $G_0/K_0$}\label{SEC Heat kernel on $G_0/K_0$}

Let $\Delta_X$ be the Laplace-Beltrami operator on $X=G_0/K_0$ corresponding to the natural metric on $X$ defined by the Killing form $B$ on $\mf{g}_0$. Then $\Delta_X$ can be shown to be descending from the Casimir element $\omega\in U(\mf{g}_0)$ of the universal enveloping algebra $U(\mf{g}_0)$, which, subject to a choice of  basis $\{X_j\}_{1\leq j\leq n}$ of $\mf{g}_0$, can be defined as the sum
\begin{align*}
    \omega=\sum_{j=1}^n X_j^*X_j\,,
\end{align*}
where $\{X_j^*\}_{1\leq j\leq n}$ is the dual  basis with respect to the Killing form $B$ on $\mf{g}_0$ (see \cite[p.~331]{Helgason2}). The spherical function $\phi_{\lambda}\,(\lambda\in\mf{a}_0\spcheck)$ is then an eigenfunction of $\Delta_X$ with eigenvalue $\lambda_{\omega}=-(\langle \rho_0,\rho_0\rangle_0+\langle\lambda,\lambda\rangle_0)$ (see \cite[p.~427] {Helgason2}), i.e.,
\begin{align*}
    \Delta_X\phi_{\lambda}(x)= \lambda_{\omega}\,\phi_{\lambda}(x)\quad(x\in G_0/K_0,\,\lambda\in\mf{a}_0\spcheck).
\end{align*}
\begin{definition}
The heat kernel on $G_0/K_0$ is the fundamental solution $K_t(x)\in L^2(G_0 /K_0)$ for each $t>0$ to the heat equation
\begin{equation}\label{Heat eq}
\begin{aligned}
    & \frac{\partial u_t(x)}{\partial t}=\Delta_X u_t(x) &\,& (t> 0,\,x\in G_0/K_0)\\
    & u_0(x)=f(x) &\,& (f\in C_c^{\infty})
\end{aligned}
\end{equation}
in the sense that for any $f\in C_c^{\infty}$, its convolution $u_t=f\ast K_t$ is a solution to the above equation satisfying $\norm{f\ast K_t}_2\rightarrow 0$ as $t\rightarrow 0$. 
\end{definition}

In \cite{Gangolli1, Gangolli2}, Gangolli, using spherical transform, constructs a function $K_t$ that has the standard properties of the fundamental solution of the heat equation on $G_0/K_0$.
\begin{theorem}
Let $\mathscr{C}(K_0\backslash G_0/K_0)$ be the $L^2$-Schwartz space defined in \autoref{spherical inversion theorem}. The function $K_t\colon G_0/K_0\rightarrow\mb{R}\;(t>0)$ defined by
\begin{align}\label{heat kernel formula}
K_t(x)=\int_{\mf{a}\spcheck}\exp(\lambda_{\omega}\,t)\,\phi_{\lambda}(x)\,\vert c(\lambda)\vert^{-2}\,\dif\lambda.
\end{align}
satisfy the following properties:
\begin{enumerate}[(a)]
    \item $K_t\in\mathscr{C}(K_0\backslash G_0/K_0)$ for each $t>0$. 
    \item $\widehat{K_t}(\lambda)=\exp(\lambda_{\omega}\,t)$ for all $\lambda\in\mf{a}\spcheck$.
    \item $K_t\ast K_s= K_{t+s}$ for all $t,s>0$.
    \item For any $f\in C_c^{\infty},\,f\ast K_t$ is a solution to the equation $\partial/\partial t=\Delta_X$ and $\norm{f\ast K_t}_2\rightarrow 0$ as $t\rightarrow 0$.
\end{enumerate}
\end{theorem}

\begin{proof}
One obtains (b) by taking a spherical transform of the heat equation \eqref{Heat eq}. Then \eqref{heat kernel formula} is obtained by taking an inverse spherical transform of (b). For further details on the proof,
see \cite[Proposition 3.1]{Gangolli1} and \cite[Theorem 1]{Gangolli2}. 
\end{proof}


\subsection{Spherical function on $\mb{H}_n$}\label{SEC Spherical function on Hn}

In this section, we obtain the spherical function on $\mathrm{Sp}_n(\mb{R})/U_n(\mb{R})\cong\mb{H}_n$ by using the general procedure for obtaining spherical functions on Riemannian symmetric spaces via the Flensted-Jensen reduction established in \autoref{SEC Flensted-Jensen reduction}. But first we need to specialize the general notation in \autoref{SEC notation} for the symplectic group in order to have a more explicit structure for these groups and algebras that in turn is essential for obtaining a more explicit structure for the spherical function and the heat kernel in this particular case.  

The Lie algebras shall as usual be denoted by gothic letters. Let
\begin{align*}
\mathfrak{g}_0=\bigg\{\begin{pmatrix} A & B \\ C & -A^t\end{pmatrix}\;\bigg|\; A,B,C\in\mathbb{R}^{n\times n},\, B=B^t,\, C=C^t\bigg\}
\end{align*}
denote the real symplectic algebra $\mf{sp}_n(\mb{R})$, while $\mf{g}$ shall denote the complex symplectic algebra $\mf{sp}_n(\mb{C})=\mf{g}_0+i\mf{g}_0$. On $\mf{g}$, we have the Cartan involution $\theta(X)=-\overline{X}^t\,(X\in\mf{g})$, which restricts to $\theta_0(X)=-X^t\,(X\in\mf{g}_0)$ on $\mf{g}_0$. Accordingly, we have the Cartan decomposition $\mf{g}=\mf{u}+\mf{p}$ of $\mf{g}$ into the $(+1)$-eigenspace
\begin{align*}
\mf{u}=\bigg\{\begin{pmatrix} A & B \\ -\overline{B} & \overline{A}\end{pmatrix}\;\bigg|\; A,B\in\mathbb{C}^{n\times n},\, A=-\overline{A}^t, B=B^t\bigg\}.
\end{align*}
and the $(-1)$-eigenspace
\begin{align*}
\mathfrak{p}=\bigg\{\begin{pmatrix} A & B \\ \overline{B} & -A\end{pmatrix}\;\bigg|\; A,B\in\mathbb{C}^{n\times n},\, A=\overline{A}^t,\, B=B^t\bigg\}.
\end{align*}
of $\theta$. Similarly, we have the Cartan decomposition $\mf{g}_0=\mf{k}_0+\mf{p}_0$ of $\mf{g}_0$ into the $(+1)$-eigenspace
\begin{align*}
\mathfrak{k}_0&=\bigg\{\begin{pmatrix} A & B \\ -B & A\end{pmatrix}\;\bigg|\; A,B\in\mathbb{R}^{n\times n},  A=-A^t, B=B^t \bigg\},
\end{align*}
and the $(-1)$-eigenspace 
\begin{align*}
\mathfrak{p}_0=\bigg\{\begin{pmatrix} A & B \\ B & -A\end{pmatrix}\;\bigg|\; A,B\in\mathbb{R}^{n\times n},\, B=B^t,\, A=A^t\bigg\}
\end{align*}
of $\theta_0$. 

\begin{remark}\label{algebra identification}
We note here that $\mf{u}$ is the symplectic unitary algebra $\mf{u}=\mf{sp}_n(\mb{C})\cap \mf{u}_{2n}$ while $\mf{k}_0$ is the symplectic real orthogonal algebra $\mf{k}_0=\mf{sp}_n(\mb{R})\cap \mf{0}_{2n}(\mb{R})$. The subalgebra $\mf{k}=\mf{k}_0+i\mf{k}_0$ is then given by the symplectic complex orthogonal algebra $\mf{k}=\mf{sp}_n(\mb{C})\cap \mf{0}_{2n}(\mb{C})$.
\end{remark}

The maximal abelian subspaces $\mf{a}$ and $\mf{a}_0$ of $\mf{p}$ and $\mf{p}_0$, respectively, are given by the diagonal elements in $\mf{p}$ and $\mf{p}_0$, respectively. As the elements of $\mf{p}=\{X\in\mf{g}\,\vert\,X=\overline{X}^t\}$ are Hermitian and the elements of $\mf{p}_0=\{X\in\mf{g}_0\,\vert\,X=X^t\}$ are symmetric, the diagonal entries in both $\mf{p}$ and $\mf{p}_0$ are real. Therefore, we have
\begin{align*}
   \mf{a}=\mf{a}_0=\bigg\{ r=\begin{pmatrix}R & 0 \\ 0 & -R\end{pmatrix} \; \bigg| \; R= \bigg(\begin{smallmatrix}
{r_{1}} &{} & {0} \\
{} & {\ddots}& \\
{0} &{}& {r_{n}}
\end{smallmatrix}\bigg),\;r_j\in\mb{R},\;1\leq j\leq n\bigg\}.
\end{align*}.
As $\mf{a}=\mf{a}_0$, we drop the distinction and from here on denote both by $\mf{a}$.  

The  basis of the dual space $\mf{a}\spcheck$ shall be denoted by $\{e_1,\,e_2,\,\ldots,e_n\}$ such that $e_j(r)=r_j$ $(r\in\mf{a},\,1\leq j\leq n)$. The generic element of $\mf{a}\spcheck$ shall be denoted by $\lambda=\lambda_1\,e_1+\ldots+ \lambda_n\,e_n\;(\lambda_j\in\mb{R},\;1\leq j\leq n)$.The Killing form on $\mf{g}$ is given by
\begin{align*}
    B(X,Y)=2(n+1)\tr(XY).
\end{align*}
Let $E_{j,k}\,(1\leq j,k\leq 2n)$ denote the $(2n\times 2n)$-matrix with entry $1$ where the $j$-th row and the $k$-th column meet, all other entries being $0$. Using this Killing form, we can assign to each  basis vector $e_j\,(1\leq j\leq n)$ in $\mf{a}\spcheck$ an element 
\begin{align*}
H_j=\frac{1}{4(n+1)}(E_{j,j}-E_{n+j,n+j})\in\mf{g}\quad\quad (1\leq j\leq n),
\end{align*}
so that $B(H_j,r)=e_j(r)$. This assignment induces a scalar product on $\mf{a}\spcheck$ defined by
\begin{align*}
    \langle e_{j},e_{k}\rangle:=B(H_{j},H_{k})=\frac{\delta_{j,k}}{4(n+1)}.
\end{align*}

The roots of $\mf{g}$ corresponding to $\mf{a}$ are given by
\begin{align*}
\Delta=\{\pm 2e_{j}\;|\; 1\leq j\leq n\}\,\cup\,\{\pm e_{j} \pm e_{k}\;|\; 1\leq j< k\leq n\},
\end{align*}
with each root space $\mf{g}^{\alpha}\,(\alpha\in\Delta)$ having real dimension (and hence root multiplicity) $m_{\alpha}=2$. As $\mf{a}=\mf{a}_0$, for roots of $\mf{g}_0$ corresponding to $\mf{a}_0$, we have $\Delta_0=\Delta$. However, in case of the real algebra $\mf{g}_0$, each root space $\mf{g}_0^{\alpha}\,(\alpha\in\Delta_0)$ has real dimension 1. The Weyl group $W=W_0$ consists of the permutations $\sigma\colon \mf{a}\rightarrow \mf{a}$ of elements $r\in\mf{a}$, i.e., 
\begin{align*}
\sigma(r)=\sigma\begin{pmatrix} R & 0 \\ 0 & -R \end{pmatrix}=\begin{pmatrix} \sigma(R) & 0 \\ 0 & -\sigma(R) \end{pmatrix}\quad\text{where}\quad \sigma(R)=\sigma \bigg(\begin{smallmatrix}
{r_{1}} &{} & {0} \\
{} & {\ddots}& \\
{0} &{}& {r_{n}}
\end{smallmatrix}\bigg) = \Bigg(\begin{smallmatrix}
{\pm r_{\tau(1)}} &{} & {0} \\
{} & {\ddots}& \\
{0} &{}& {\pm r_{\tau(n)}}
\end{smallmatrix}\Bigg)\quad(\tau\in S_n).
\end{align*}
Consequently, we have $W=\mb{Z}/n\mb{Z}\times S_n$.

The canonical choice of positive roots in $\Delta$ is given by
\begin{align}\label{positive roots}
\Delta^+ =\{2e_{j}\;|\; 1\leq j\leq n\}\,\cup\,\{e_{j}+e_{k}\;|\; 1\leq j<k\leq n\}\,\cup\, \{e_{j}-e_{k}\;|\; 1\leq j< k\leq n\}.
\end{align}
Thus, the half-root sum $\rho_0=1/2\sum_{\alpha\in\Delta^+}m_{\alpha}\alpha$ in the real algebra $\mf{g}_0$ is given by
\begin{align*}
    \rho_0=n\,e_1+(n-1)\,e_2+\ldots+(n-j+1)\,e_j+\ldots+2\,e_{n-1}+ e_n
\end{align*}
and the half-root sum $\rho$ in the complex algebra $\mf{g}$ is given by $\rho=2\,\rho_0$. Corresponding to the choice \eqref{positive roots} of positive roots, the positive Weyl chamber $\mf{a}^{+}$ of $\mf{a}$ is given by
\begin{align*}
\mf{a}^{+}=\bigg\{ r=\begin{pmatrix}R & 0 \\ 0 & -R\end{pmatrix} \; \bigg| \; R= \bigg(\begin{smallmatrix}
{r_{1}} &{} & {0} \\
{} & {\ddots}& \\
{0} &{}& {r_{n}}
\end{smallmatrix}\bigg),\;r_1\geq\ldots \geq r_j\geq \ldots r_n\geq 0\bigg\}
\end{align*}
and the nilpotent algebra $\mf{n}=\sum_{\alpha\in\Delta^{+}} \mf{g}^{\alpha}$ is given by
\begin{align*}
\mathfrak{n}&=\bigg\{\begin{pmatrix} P & Q \\ 0 & -P^t\end{pmatrix}\;\bigg|\; P,Q\in\mathbb{C}^{n\times n},\, P \text{ strictly upper-triangular},\, Q \text{ symmetric}\bigg\},
\end{align*}
wherefrom $\mf{n}_0=\sum_{\alpha\in\Delta_0^{+}} \mf{g}_0^{\alpha}$ can be obtained by restricting to $\mf{g}_0$.

Groups shall as usual continue to be denoted by capital Roman letters. In particular, $G_0$ shall denote the real symplectic group $\mathrm{Sp}_n(\mb{R})$, while $G$ shall denote the complex symplectic group $\mathrm{Sp}_n(\mb{C})$. By \autoref{algebra identification}, the subgroup $K_0=\exp(\mf{k}_0)\subsetneq G_0$ is given by the real orthogonal subgroup $\mr{Sp}_n(\mb{R})\cap O_{2n}(\mb{R})$ of $\mr{Sp}_n(\mb{R})$, while $K=\exp(\mf{k})\subsetneq G$ is given by the complex orthogonal group $\mr{Sp}_n(\mb{C})\cap O_{2n}(\mb{C})$ of $\mr{Sp}_n(\mb{C})$. Group elements, i.e., the matrices in the matrix groups shall continue to be denoted by small Roman letters. The scalar entries of the matrices shall also be denoted by small letters, while matrix-blocks in matrices written in a block-matrix format shall be denoted by capital letters.

Both real and complex symplectic orthogonal matrices have the same structure
\begin{align*}
    k=\begin{pmatrix} A & B\\ -B & A\end{pmatrix}\quad\quad (AA^t+BB^t=\mathbbm{1}_n\,,\,AB^t=BA^t),
\end{align*}
but while for $k\in K_0$ this implies that the matrix $A+iB$ is an $(n\times n)$-unitary matrix, no such implication is possible in case of complex orthogonal symplectic matrices $k\in K\setminus K_0$. However, any complex orthogonal matrix $k\in K$ can be represented as $k=k_0\; k_h$, where $k_0\in K_0$ is a real orthogonal matrix and $k_h$ is a Hermitian orthogonal matrix (see \cite[Theorem 1]{Gantmacher}).
Therefore, a general $k\in K$ can be represented by
\begin{align*}
    k=\begin{pmatrix} A_0 & B_0\\ -B_0 & A_0\end{pmatrix}\begin{pmatrix} A & B\\ -B & A\end{pmatrix}
\end{align*}
such that $A_0+iB_0$ is $(n\times n)$-unitary and $A+iB$ is $(n\times n)$-Hermitian.

The group $A=\exp(\mf{a})$ is given by the group of real diagonal symplectic matrices
\begin{align}\label{diagonal notation}
   A=\bigg\{ \exp(r)=\begin{pmatrix}\exp(R) & 0 \\ 0 & \exp(-R) \end{pmatrix} \; \bigg| \; R= \bigg(\begin{smallmatrix}
{r_{1}} &{} & {0} \\
{} & {\ddots}& \\
{0} &{}& {r_{n}}
\end{smallmatrix}\bigg),\;r_j\in\mb{R},\;1\leq j\leq n\bigg\}.
\end{align}
By \autoref{algebra identification}, the group $U=\exp(\mf{u})\subsetneq G$ is 
given by the unitary subgroup $U=\mr{Sp}_n(\mb{C})\cap U_{2n}=Sp(n)$ of $\mr{Sp}_n(\mb{C})$, whose elements can be shown to have the structure
\begin{align*}
U=\bigg\{\begin{pmatrix} A & B \\ -\overline{B} & \overline{A}\end{pmatrix}\;\bigg|\; A,B\in\mathbb{C}^{n\times n}, A\overline{A}^t+B\overline{B}^t=\mathbbm{1}_n,\,AB^t=BA^t\bigg\}.
\end{align*}
The group $N=\exp(\mf{n})\subsetneq G$ is given by
\begin{align*}
N=\bigg\{\begin{pmatrix} P & Q \\ 0 & P^{-t}\end{pmatrix}\;\bigg|\; P,Q\in\mathbb{C}^{n\times n},\, P \text{ unit upper-triangular},\, PQ^t=QP^t\bigg\},
\end{align*}
wherefrom $N_0=\exp(\mf{n}_0)\subsetneq G_0$ can be obtained by restricting to $G_0$. 

The Haar measures of the groups shall be denoted by $\dif\mu(x)$, while the Euclidean measures shall be denoted by $\dif x$.

This prepares the setup needed to compute the spherical function on $\mathrm{Sp}_n(\mb{C})$ corresponding to 
$\lambda\in \mf{a}\spcheck$ using the formula
\begin{align}\label{2complex-spherical}
\Phi_{\lambda}(\exp(r))=\frac{\pi(\rho)}{\pi(i \lambda)}\, \frac{\sum_{\sigma\in W}
\det(\sigma) \exp(i \sigma \lambda(r))}{\sum_{\sigma\in W}
\det(\sigma) \exp( \sigma \rho (r))}
\end{align}
in \autoref{complex spherical}. The quantity
\begin{align*}
    \pi(\rho)=\prod_{\alpha\in\Delta^+}\langle 2n\,e_1+\ldots+2(n-j+1)\,e_j+\ldots+2e_n\,,\,\alpha\rangle
\end{align*}
will come out to be a positive real constant depending only on $n$. Since we are not interested in the exact nature of this dependence and keeping track of these constants soon get quite tedious, we shall club all such constants which are not crucial to our calculation under the generic symbol $c_n$, which should be interpreted as a positive real constant depending only on $n$. 

The quantity $\pi(i\lambda)$ is of the form
\begin{align*}
    \pi(i\lambda)=i^{n^2}\prod_{1\leq j\leq n} \frac{2\lambda_j}{4(n+1)}\,
    \prod_{1\leq j<k\leq n} \frac{\lambda_j+\lambda_k}{4(n+1)}\,\prod_{1\leq j<k\leq n} \frac{\lambda_j-\lambda_k}{4(n+1)}.
\end{align*}
The above homogeneous polynomial of degree $n^2$ plays a crucial role in our analysis. Let us formally denote it by
\begin{align*}
    \varepsilon(\lambda):=\varepsilon(\lambda_1,\,\lambda_2,\,\ldots,\,\lambda_n)
    =\prod_{1\leq j\leq n} \lambda_j\,
    \prod_{1\leq j<k\leq n} (\lambda_j+\lambda_k)\,\prod_{1\leq j<k\leq n} (\lambda_j-\lambda_k).
\end{align*}
 Under transpositions $\sigma_{j,k}\in W$ given by the assignment 
 \begin{align*}
 (\lambda_1,..,\lambda_j,..,\lambda_k,..,\lambda_n)\mapsto(\lambda_1,..,\lambda_k,..,\lambda_j,..,\lambda_n)
 \end{align*}
 and sign-changes $\sigma_{j}\in W$ given by the assignment
 \begin{align*}
     (\lambda_1,..,\lambda_j,..,\lambda_n)\mapsto(\lambda_1,..,-\lambda_j,..,\lambda_n), 
\end{align*}
 we have $\sigma_{j,k}((\lambda_j+\lambda_k)(\lambda_j-\lambda_k))=-(\lambda_j+\lambda_k)(\lambda_j-\lambda_k)$ and $\sigma_j((\lambda_j+\lambda_k)(\lambda_j-\lambda_k))=(\lambda_j+\lambda_k)(\lambda_j-\lambda_k)$, respectively. Thus, in both the cases we have $\varepsilon(\sigma_{j,k}(\lambda))=-\varepsilon(\lambda)$ and $\varepsilon(\sigma_j (\lambda))=-\varepsilon(\lambda)$. Now, since the Weyl group $W$ is generated by these transpositions and sign-changes, the polynomial $\varepsilon$ for any formal variable $\lambda=(\lambda_1,\ldots,\lambda_n)$ satisfies the property
 \begin{align}\label{symmetric weyl}
     \varepsilon(\sigma (\lambda))=\det(\sigma)\varepsilon(\sigma(\lambda))\quad\quad(\sigma\in W).
 \end{align}

Next we consider the denominator $\sum_{\sigma\in W}\det(\sigma) \exp( \sigma \rho (r))$ in the right-hand side of the \autoref{2complex-spherical}. As $\rho=\sum_{\alpha\in\Delta^+} \alpha$, this denominator can be expressed in the form
\begin{align*}
    \sum_{\sigma\in W}\det(\sigma) \exp( \sigma \rho (r))=\prod_{\alpha\in\Delta^+}(\exp(\alpha(r))-\exp(-\alpha(r))),
\end{align*}
whence we obtain
\begin{align*}
    \sum_{\sigma\in W}\det(\sigma) \exp(\sigma\rho (r))=2^{n^2}\prod_{1\leq j\leq n} \sh(2r_j)
    \prod_{1\leq j<k\leq n} \sh(r_j+r_k)\prod_{1\leq j<k\leq n} \sh(r_j-r_k).
\end{align*}
The product of the $\sh'$s on the right-hand side of the formula would also play an important role in our later analysis. We formally denote it as
\begin{align*}
    \delta(r):=\prod_{1\leq j\leq n} \sh(r_j) \prod_{1\leq j<k\leq n} \sh\Big(\frac{r_j+r_k}{2}\Big)\prod_{1\leq j<k\leq n} \sh\Big(\frac{r_j-r_k}{2}\Big).
\end{align*}
This gives us the formula
\begin{align}\label{complex spherical}
    \Phi_{\lambda}(\exp(r))=\frac{c_n}{i^{n^2}}\,\frac{\sum_{\sigma \in W}
\det(\sigma) \exp(i \sigma \lambda(r))}{\varepsilon(\lambda)\,\delta(2r)}
\end{align}
for the spherical function on $\mathrm{Sp}_n(\mb{C})$ corresponding to 
$\lambda\in \mf{a}\spcheck$ at $\exp(r)\in A$. An interesting limiting case of this formula is to determine the spherical function on $\mathrm{Sp}_n(\mb{C})$ corresponding to $\lambda=0$ at $\exp(r)\in A$, which we consider in the next proposition.
\begin{proposition}\label{0-complex spherical}
The spherical function on $\mathrm{Sp}_n(\mb{C})$ corresponding to $\lambda=0$ at $\exp(r)\in A$ is given by the formula
\begin{align*}
    \Phi_{0}(\exp(r))=c_n \frac{\varepsilon(r)}{\delta(2r)}=c_n \prod_{1\leq j\leq n} \frac{r_j}{\sh(2r_j)}\,
    \prod_{1\leq j<k\leq n} \frac{r_j+r_k}{\sh(r_j+r_k)}\,\prod_{1\leq j<k\leq n} \frac{r_j-r_k}{\sh(r_j-r_k)},
\end{align*}
where $c_n$ denotes a constant depending only on $n$. 
\end{proposition}
\begin{proof}
We begin by defining a polynomial differential operator 
\begin{align*}
\ve(\partial/\partial\lambda):=\ve(\partial/\partial\lambda_1,\ldots,\partial/\partial\lambda_j,\ldots,\partial/\partial\lambda_n).
\end{align*}
Since both the numerator and the denominator in the formula
\begin{align*}
    \Phi_{\lambda}(\exp(r))=\frac{c_n}{i^{n^2}}\,\frac{\sum_{\sigma \in W}
\det(\sigma) \exp(i \sigma \lambda(r))}{\varepsilon(\lambda)\,\delta(2r)}
\end{align*}
are analytic at $\lambda=0$, the limit at $\lambda\rightarrow 0$ can be calculated as
\begin{align*}
    \Phi_{0}(\exp(r))&=\lim_{\lambda\rightarrow 0}\frac{c_n}{i^{n^2}}\,\frac{\sum_{\sigma \in W}
\det(\sigma) \exp(i \sigma \lambda(r))}{\varepsilon(\lambda)\,\delta(2r)}\\
                &=\frac{c_n}{i^{n^2}\,\delta(2r)}\lim_{\lambda\rightarrow 0}
                \frac{\ve(\partial/\partial\lambda)\Big(\sum_{\sigma \in W}
\det(\sigma) \exp(i \sigma\lambda(r))\Big)}{\ve(\partial/\partial\lambda)\varepsilon(\lambda)},
\end{align*}
provided both the derivatives converge at $\lambda=0$. 

Now, since $\partial/\partial\lambda_j(\exp(i\sigma\lambda(r)))=\sigma(i\,r_j)\exp(i\sigma\lambda(r)$, it is easy to see that 
\begin{align*}
\ve(\partial/\partial\lambda)\exp(i \sigma \lambda(r))=\ve(i\sigma(r))\exp(i \sigma \lambda(r)),
\end{align*}
which, by the property (\ref{symmetric weyl}) of $\ve$, becomes
\begin{align}\label{epsilon derivative}
\ve(\partial/\partial\lambda)\exp(i \sigma \lambda(r))=i^{n^2}\det(\sigma)\,\ve(r)\exp(i \sigma \lambda(r)).
\end{align}
Thus, the derivative of the numerator by $\ve(\partial/\partial\lambda)$, in the limit $\lambda\rightarrow 0$, converges to $i^{n^2}\det(\sigma)^2\ve(r)\vert W\vert=i^{n^2}\ve(r) 2^n n!$. 

For a monomial $\lambda_1^{\alpha_1}\cdots\lambda_j^{\alpha_j}\cdots \lambda_n^{\alpha_n}\,(\alpha_j\in\mb{N}_{\geq 0},\,\sum_{j=1}^n \alpha_j=n^2)$ and a differential operator 
\begin{align*}
(\partial/\partial\lambda_1)^{\beta_1}\ldots (\partial/\partial\lambda_j)^{\beta_j}\ldots (\partial/\partial\lambda_n)^{\beta_n}\quad(\beta_j\in\mb{N}_{\geq 0},\,\sum_{j=1}^n \beta_j=\sum_{j=1}^n \alpha_j=n^2), 
\end{align*}
we have 
\begin{align*}
    & \Bigg(\dpd[\beta_1]{}{\lambda_1}\Bigg)\cdots\Bigg(\dpd[\beta_j]{}{\lambda_j}\Bigg)\cdots\Bigg(\dpd[\beta_n]{}{\lambda_n}\Bigg) \big(\lambda_1^{\alpha_1}\cdots\lambda_j^{\alpha_j}\cdots \lambda_n^{\alpha_n}\big)\\
    &\quad=\begin{cases}\displaystyle
    \prod_{j=1}^n \frac{\alpha_j(\alpha_j+1)}{2} & \text{ for }\alpha_j=\beta_j\quad \forall 1\leq j\leq n,\\
    0 & \text{ otherwise}.
    \end{cases}
\end{align*}
Therefore it is easy to see that $\ve(\partial/\partial\lambda)\ve(\lambda)$ is a constant depending only on $n$ and hence we have the requisite limit stated in the proposition. 
\end{proof}

Next we apply \autoref{Flensted reduction} to calculate the spherical function $\phi_{\lambda}$ on $G_0=\mathrm{Sp}_n(\mb{R})$ corresponding to $\lambda\in \mf{a}_0\spcheck$ by reducing it to the complex case. 
\begin{theorem}\label{n-spherical theorem}
The spherical function $\phi_{\lambda}$ on $G_0=\mathrm{Sp}_n(\mb{R})$ corresponding to
\begin{align*}
\lambda=\lambda_1\,e_1+\ldots+\lambda_j\,e_j+\ldots+\lambda_n\,e_n\in\mf{a}\spcheck
\end{align*}
at $\exp(r)=\big(\begin{smallmatrix} \exp(R) & 0\\0 & \exp(-R)\end{smallmatrix}\big)\in A$ with 
\begin{align*}
R=\begin{pmatrix}
{r_{1}} &{} & {0} \\
{} & {\ddots}& \\
{0} &{}& {r_{n}}
\end{pmatrix}\quad\quad(r_j\in\mb{R}_{\geq 0},\;1\leq j\leq n),
\end{align*}
is given by
\begin{align}\label{n-spherical function}
    \phi_{\lambda}(\exp(r))=& \frac{c_n}{i^{n^2}\,\tau(\lambda)}
    \int\limits_{k\in K}
    \frac{\sum_{\sigma\in W}
\det(\sigma) \exp\big(i\sigma\lambda(\vr(r,k))\big)}{\delta(\vr(r,k))}\dif\mu(k),
\end{align}
where $\vr(r,k)$ is the diagonal matrix $\vr(r,k)=\big(\begin{smallmatrix}P(r,k) & 0\\0 & -P(r,k)\end{smallmatrix}\big)$ with 
\begin{align*}
P(r,k)= \begin{pmatrix}
{\vr_{1}(r,k)} &{} & {0} \\
{} & {\ddots}& \\
{0} &{}& {\vr_{n}(r,k)}
\end{pmatrix}\quad (\vr_j(r,k)\in\mb{R},\;1\leq j\leq n)
\end{align*}
related to $r=\big(\begin{smallmatrix}R & 0\\0 & -R\end{smallmatrix}\big)$ via the matrix equality $k \exp(r)\overline{k}^t=u\exp(\vr)\overline{u}^t$ with $k\in K$ and $u\in U$. The functions $\tau$ and $\delta$ are given by
\begin{align*}
    \tau(\lambda)&=\prod_{1\leq j\leq n} \operatorname{th} \Big(\frac{\lambda_{j}}{2}\pi\Big) \prod_{1\leq j<k \leq n} \operatorname{th} \Big(\frac{\lambda_{j}+\lambda_{k}}{2} \pi\Big) \prod_{1 \leq j<k \leq n}  \operatorname{th} \Big(\frac{\lambda_{j}-\lambda_{k}}{2} \pi\Big) ,\\
    \delta(\vr)&=\prod_{1\leq j\leq n} \sh(\vr_j)\prod_{1\leq j<k\leq n} \sh\Big(\frac{\vr_j+\vr_k}{2}\Big)\prod_{1\leq j<k\leq n} \sh\Big(\frac{\vr_j-\vr_k}{2}\Big),
\end{align*}
while $c_n$ is a positive real constant depending only on $n$. 
\end{theorem}
\begin{proof}
We begin by putting $x=\exp(r/2)$ in \autoref{Flensted reduction}, so that $x\theta(x)^{-1}=xx^t=\exp(r)$, and we have
\begin{align}\label{local flensted}
\phi_{\lambda}(\exp(r))=|c(\lambda)|^{2}|\pi_{0}(\lambda)|^{2} \int\limits_{k\in K} \Phi_{2 \lambda}(k \exp(r/2)) \dif \mu(k) \quad (\lambda \in \mathfrak{a}_{0}\spcheck),
\end{align}
where $\Phi_{2\lambda}$ is the spherical function on the complex group $G=\mathrm{Sp}_n(\mb{C})$ corresponding to $2\lambda\in\mf{a}\spcheck$.
As $k \exp(r/2)\in G$, from the $G=U \overline{A^{+}} U$ decomposition of the complex group, we obtain $\vr(r,k)\in \overline{A^{+}}$ and $u,v\in U$, so that $g=k\exp(r/2)=u\exp(\vr(r,k)/2)v\in G$. Therefore, we have
\begin{align}\label{connecting formula}
    g\overline{g}^t= 
    k \exp(r)\overline{k}^t=u\exp(\vr(r,k))\overline{u}^t.
\end{align}
Furthermore, since the spherical function $\Phi_{2\lambda}\in C^{\infty}(U\backslash G/U)$ on $G$ is bi-invariant under $U$, we have
\begin{align*}
    \Phi_{2 \lambda}(k \exp(r/2))=\Phi_{2 \lambda}(u\exp(\vr(r,k)/2)v)=\Phi_{2 \lambda}(\exp(\vr(r,k)/2)).
\end{align*}
Now, we plug in the formula for $\Phi_{2 \lambda}(\exp(\vr(r,k)/2))$ in \autoref{local flensted} using the formula (\ref{complex spherical}) to obtain
\begin{align}\label{local integral}
   \phi_{\lambda}(\exp(r))=c_n\frac{|c(\lambda)|^{2}|\pi_{0}(\lambda)|^{2}}{i^{n^2}\ve(2\lambda)} \int\limits_{k\in K} \,\frac{\sum_{\sigma\in W}
\det(\sigma) \exp(i \sigma 2\lambda(\vr(r,k)/2))}{\delta(\vr(r,k))} \dif\mu(k).
\end{align}
Now $2\lambda(\vr(r,k)/2)=\lambda(\vr(r,k))$, so the numerator inside the integral becomes
\begin{align*}
\sum_{\sigma\in W} \det(\sigma) \exp(i \sigma \lambda(\vr(r,k))). 
\end{align*}
Next $\ve(2\lambda)=2^{n^2}\ve(\lambda)$ and the positive real constant $2^{n^2}$ gets assumed within the generic $c_n$. Also, from the formula   
\begin{align*}
\abs{c(\lambda)}^{-2} =&\frac{1}{\pi^{n^{2} / 2}}\prod_{1\leq j\leq n}\frac{\lambda_{j}}{2} \operatorname{th} \Big(\frac{\lambda_{j}}{2}\pi\Big) \prod_{1\leq j<k \leq n}\frac{\lambda_{j}+\lambda_{k}}{2} \operatorname{th} \Big(\frac{\lambda_{j}+\lambda_{k}}{2} \pi\Big) \times\\ 
&\quad\times \prod_{1 \leq j<k \leq n} \frac{\lambda_{j}-\lambda_{k}}{2} \operatorname{th} \Big(\frac{\lambda_{j}-\lambda_{k}}{2} \pi\Big)
\end{align*}
for the Harish-Chandra $c$-function on $\mathrm{Sp}_n(\mb{R})$ due to Bhanu Murti mentioned in \autoref{SEC spherical general}, it clearly follows, that in terms of the special functions $\ve$ and $\tau$ introduced to make our calculations less cumbersome, the above formula can be simply written as
\begin{align}\label{c-reduction}
    \abs{c(\lambda)}^{-2}=c_n \,\ve(\lambda)\tau(\lambda). 
\end{align}
Thus, noting that
\begin{align*}
    |\pi_{0}(\lambda)|^{2}=\bigg\vert \prod_{1\leq j\leq n} \frac{2\lambda_j}{4(n+1)}\,
    \prod_{1\leq j<k\leq n} \frac{\lambda_j+\lambda_k}{4(n+1)}\,\prod_{1\leq j<k\leq n} \frac{\lambda_j-\lambda_k}{4(n+1)} \bigg\vert^2
    =c_n\,\ve(\lambda)^2,
\end{align*}
we see that $|c(\lambda)|^{2}|\pi_{0}(\lambda)|^{2}=c_n(1/\tau(\lambda))$, which brings the integral in the right-hand side of \autoref{local integral} to \begin{align*}
    \phi_{\lambda}(\exp(r))=& \frac{c_n}{i^{n^2}\,\tau(\lambda)}
    \int\limits_{k\in K}
    \frac{\sum_{\sigma\in W}
\det(\sigma) \exp(i \sigma \lambda(\vr(r,k)))}{\delta(\vr(r,k))}\dif\mu(k),
\end{align*}
thereby proving the theorem.
\end{proof}

\begin{remark}
Since the real symplectic matrix $\exp(r)=(\begin{smallmatrix} \exp(R) & 0\\0 & \exp(-R)\end{smallmatrix})\in A\subsetneq G_0=\mathrm{Sp}_n(\mb{R})$ maps the point $i\mathbbm{1}_n\in\mb{H}_n$ by symplectic action to $Z=i\,\exp(2R)\in\mb{H}_n$ and $\phi_{\lambda}(\exp(r))\in  C^{\infty}(K_0\backslash G_0/K_0) $  is a radial function on $\mb{H}_n$, the formula for $\phi_{\lambda}$ in the above theorem also gives the formula for spherical function on $\mb{H}_n$ at a point $Z=k_0\, i\exp(2R)\in\mb{H}_n\,(k_0\in K_0)$ and hence is also denoted by $\phi_{\lambda}(2R)$. 
\end{remark}

\subsection{Heat kernel on $\mb{H}_n$}\label{SEC Heat kernel on Hn}

In this section, we obtain the heat kernel on $\mb{H}_n$ by following the general procedure established in \autoref{SEC Flensted-Jensen reduction}. We continue with the notation and the basic setup fixed in \autoref{SEC Spherical function on Hn}.

We begin by computing the eigenvalue  $\lambda_{\omega}=-(\langle \rho_0,\rho_0\rangle_0+\langle\lambda,\lambda\rangle_0)$ for the Casimir operator $\omega$ on $G_0=\mathrm{Sp}_n(\mb{R})$. For the  basis vectors $e_j\,(1\leq j\leq n)$ in $\mf{a}_0\spcheck$, the inner product induced by the Killing form of $\mf{g}_0$ on $\mf{a}_0\spcheck$ takes the form
\begin{align*}
    \langle e_j,e_k\rangle=\frac{\delta_{j,k}}{4(n+1)}\quad\quad(1\leq j,k\leq n).
\end{align*}
Then, for the half-root sum
\begin{align*}
    \rho_0=n\,e_1+(n-1)\,e_2+\ldots+(n-j+1)\,e_j+\ldots+2\,e_{n-1}+e_n
\end{align*}
in $\mf{g}_0$, the inner product $\langle \rho_0,\rho_0\rangle_0$ turns out to be
\begin{align*}
    \langle \rho_0,\rho_0\rangle_0=\frac{1^2+2^2+\ldots+n^2}{4(n+1)}. 
\end{align*}
Similarly, for $\lambda=\lambda_1e_1+\lambda_2e_2+\ldots+\lambda_ne_n\in\mf{a}_0\spcheck$, we have
\begin{align*}
    \langle \lambda,\lambda\rangle_0=\frac{\lambda_1^2+\lambda_2^2+\ldots+\lambda_n^2}{4(n+1)}. 
\end{align*}
the Casimir operator $\omega$ on $G_0=\mathrm{Sp}_n(\mb{R})$ descends on the Siegel upper half-space $\mb{H}_n:=\{Z=X+iY\;|\; X,Y\in\mb{R}^{n\times n}, X=X^t, Y=Y^t, Y>0\}$ to the operator
\begin{align*}
    \Delta=\frac{1}{(n+1)}\tr\Bigg(Y\bigg(\bigg(Y\dpd{}{ X}\bigg)^t\dpd{}{ X}+\bigg(Y\dpd{}{ Y}\bigg)^t\dpd{}{ Y}\bigg)\Bigg).
\end{align*}
Traditionally, this factor of $1/(n+1)$ is ignored and the Laplace--Beltrami operator on $\mb{H}_n$ is written as
\begin{align*}
    \Delta=\tr\Bigg(Y\bigg(\bigg(Y\dpd{}{ X}\bigg)^t\dpd{}{ X}+\bigg(Y\dpd{}{ Y}\bigg)^t\dpd{}{ Y}\bigg)\Bigg),
\end{align*}
due to which, we correct the value of $\lambda_{\omega}$ calculated above by multiplying it with a factor of $(n+1)$, thereby setting
\begin{align*}
   \lambda_{\omega}= -\frac{\sum_{j=1}^n j^2+\sum_{j=1}^n \lambda_j^2}{4}.
\end{align*}
Now we are ready to compute the heat kernel $K_t$ on $\mb{H}_n$ using \autoref{n-spherical theorem} and the formula (\ref{heat kernel formula}), which is the subject of the next theorem.
\begin{theorem}\label{heat kernel theorem}
The heat kernel $K_t$ at a point $Z=k_0\,i\exp(2R)\in\mb{H}_n$ on the Siegel upper half-space with $k_0\in K_0$ and 
\begin{align*}
R=\begin{pmatrix}
{r_{1}} &{} & {0} \\
{} & {\ddots}& \\
{0} &{}& {r_{n}}
\end{pmatrix}\quad\quad(r_j\in\mb{R}_{\geq 0},\;1\leq j\leq n),
\end{align*} 
is given by
\begin{align*}
   K_t(2R)=c_n\;\frac{\exp{\big(-\sum_{j=1}^n j^2t/4\big)}}{t^{n^2+n/2}} \int\limits_{k\in K} \frac{\ve(\vr(r,k))\exp\Big(-\sum_{j=1}^n\,\vr_j(r,k)^2/t\Big)}{\delta(\vr(r,k))}\dif\mu(k),
\end{align*}
where $\vr(r,k)$ is the diagonal matrix $\vr(r,k)=\big(\begin{smallmatrix}P(r,k) & 0\\0 & -P(r,k)\end{smallmatrix}\big)$ with 
\begin{align*}
P(r,k)= \begin{pmatrix}
{\vr_{1}(r,k)} &{} & {0} \\
{} & {\ddots}& \\
{0} &{}& {\vr_{n}(r,k)}
\end{pmatrix}\quad(\vr_j(r,k)\in\mb{R},\;1\leq j\leq n)
\end{align*}
related to $r=\big(\begin{smallmatrix}R & 0\\0 & -R\end{smallmatrix}\big)$ via the matrix equality $k\exp(r)\overline{k}^t=u\exp(\vr)\overline{u}^t$ with $k\in K$ and $u\in U$. The functions $\ve$ and $\delta$ are given by
\begin{align*}
        \varepsilon(\vr)&=\prod_{1\leq j\leq n} \vr_j\,
    \prod_{1\leq j<k\leq n} (\vr_j+\vr_k)\,\prod_{1\leq j<k\leq n} (\vr_j-\vr_k) ,\\
    \delta(\vr)&=\prod_{1\leq j\leq n} \sh(\vr_j) \prod_{1\leq j<k\leq n} \sh\Big(\frac{\vr_j+\vr_k}{2}\Big)\prod_{1\leq j<k\leq n} \sh\Big(\frac{\vr_j-\vr_k}{2}\Big),
\end{align*}
while $c_n$ is a positive real constant depending only on $n$. 
\end{theorem}
\begin{proof}
In \autoref{n-spherical theorem}, we had calculated the spherical function $\phi_{\lambda}$ corresponding to
\begin{align*}
\lambda=\lambda_1\,e_1+\ldots+\lambda_j\,e_j+\ldots+\lambda_n\,e_n\in\mf{a}\spcheck
\end{align*}
at $Z=k_0\,i\exp(2R)$ as
\begin{align*}
    \phi_{\lambda}(2R)=& \frac{c_n}{i^{n^2}\,\tau(\lambda)}
    \int\limits_{k\in K}
    \frac{\sum_{\sigma\in W}
\det(\sigma) \exp(i \sigma \lambda(\vr(r,k)))}{\delta(\vr(r,k))}\dif\mu(k),
\end{align*}
Therefore, using the formula (\ref{heat kernel formula}), we have
\begin{align}\label{local heat kernel}
    K_t(2R)=\frac{c_n}{i^{n^2}} \exp{\Big(-\sum_{j=1}^n j^2t/4\Big)} \int\limits_{k\in K} \frac{I(\vr(r,k))}{\delta(\vr(r,k))}\dif\mu(k)
\end{align}
where the function $I(\vr(r,k))$ given by the integral
\begin{align*}
    I(\vr(r,k))=\sum_{\sigma\in W}\det(\sigma)  \int\limits_{\lambda\in\mf{a}\spcheck} \frac{\vert c(\lambda)\vert^{-2}}{\tau(\lambda)}\exp{\Big(-\sum_{j=1}^n \lambda_j^2 t/4}+i \sigma \lambda(\vr(r,k))\Big)\dif\lambda.
\end{align*}
As we noted in \autoref{c-reduction}, the quantity $\vert c(\lambda)\vert^{-2}/\tau(\lambda)$ is just the polynomial $c_n \ve(\lambda)$, where $c_n$ is a positive real constant depending only on $n$ . So our integral simply becomes 
\begin{align*}
    I(\vr(r,k))=c_n \sum_{\sigma\in W}\det(\sigma)  \int\limits_{\lambda\in\mf{a}\spcheck} \ve(\lambda)\exp{\Big(-\sum_{j=1}^n \lambda_j^2 t/4}+i \sigma \lambda(\vr(r,k))\Big)\dif\lambda.
\end{align*}
Also, as in \autoref{0-complex spherical}, we had noted in \autoref{epsilon derivative} that for the polynomial differential operator $\ve(\partial/\partial\lambda):=\ve(\partial/\partial\lambda_1,\ldots,\partial/\partial\lambda_j,\ldots,\partial/\partial\lambda_n)$ we have
\begin{align*}
\ve(\partial/\partial\lambda)\exp\big(i \sigma \lambda(\vr(r,k))\big)=i^{n^2}\det(\sigma)\,\ve(\vr(r,k))\exp\big(i \sigma \lambda(\vr(r,k))\big),
\end{align*}
the integral $I(\vr(r,k))$ reduces to calculating the derivative by $\ve(\partial/\partial\lambda)$ of the integral
\begin{align*}
    I_0(\vr(r,k))=\sum_{\sigma\in W} \int\limits_{\lambda\in\mf{a}\spcheck} \exp{\Big(-\sum_{j=1}^n \lambda_j^2 t/4}+i \sigma \lambda(\vr(r,k))\Big)\dif\lambda
\end{align*}
as we have $c_n\,\ve(\partial/\partial\lambda)I_0(\vr(r,k))=i^{n^2}I(\vr(r,k))$. This last integral splits into integrals over the individual $\lambda_j$-s $(1\leq j\leq n)$ as
\begin{align*}
    I_0(\vr(r,k))=\sum_{\sigma\in W}\,\prod\limits_{j=1}^n\; \int\limits_{\lambda_j=-\infty}^{\infty} \exp\big(-\lambda_j^2 t/4+i\lambda_j \sigma(\vr(r,k))\big)\dif\lambda_j.
\end{align*}
But as we saw before, these individual integrals over $\lambda_j$-s are simply
\begin{align*}
   \int\limits_{\lambda_j=-\infty}^{\infty} \exp\big(-\lambda_j^2 t/4+i\lambda_j \sigma(\vr(r,k))\big)\dif\lambda= \frac{2\sqrt{\pi}\exp\big(-\sigma(\vr(r,k))^2/t\big)}{\sqrt{t}}.
\end{align*}
Therefore, their product over $1\leq j\leq n$ becomes
\begin{align*}
    \prod\limits_{j=1}^n\; \int\limits_{\lambda_j=-\infty}^{\infty} \exp\big(-\lambda_j^2 t/4+i\lambda_j \sigma(\vr(r,k))\big)\dif\lambda_j=
    \frac{(2\sqrt{\pi})^n}{t^{n/2}}\exp\big(-\sum_{j=1}^n\,\sigma(\vr(r,k))^2/t\big).
\end{align*}
However, as we have $\sum_{j=1}^n\,\sigma(\vr(r,k))^2=\sum_{j=1}^n\,\vr_j(r,k)^2$, the integral $I_0(\vr(r,k))$ evaluates to give
\begin{align*}
    I_0(\vr(r,k))=\frac{c_n}{t^{n/2}} \exp\Big(-\sum_{j=1}^n\,\vr_j(r,k)^2/t\Big).
\end{align*}
Therefore, we have
\begin{align*}
    I(\vr(r,k))=\frac{c_n}{i^{n^2}} \,\ve(\partial/\partial\lambda)I_0(\vr(r,k))
    =\frac{c_n}{i^{n^2}t^{n/2}} \,\ve(-\vr(r,k)/t)\exp\Big(-\sum_{j=1}^n\,\vr_j(r,k)^2/t\Big).
\end{align*}
Now as $\ve$ is a homogeneous polynomial of order $n^2$, we have
\begin{align*}
    I(\vr(r,k))=\frac{c_n\,i^{n^2}}{t^{n^2+n/2}}\, \ve(\vr(r,k))\exp\Big(-\sum_{j=1}^n\,\vr_j(r,k)^2/t\Big). 
\end{align*}
Putting this back to \autoref{local heat kernel}, we have the theorem. 
\end{proof}


\subsection{weight-$\kappa$ correction}

We continue with the notation in subsections \ref{SEC Spherical function on Hn} and \ref{SEC Heat kernel on Hn}. 
Given a function $f\colon G_0/K_0\rightarrow \mb{C}$ and $g\in G_0$ define the function $f^g\colon G_0/K_0\rightarrow \mb{C}$ 
by $f^{g}(x):= f(g^{-1}x)\,(x\in X=G_0/K_0)$.

As the spherical function on $X$ is supposed to be invariant under the left action $f\mapsto f^g\,(g\in G_0)$ of the elements of $K_0$ on the functions $f\colon X\rightarrow\mb{C}$, it is constructed by having an eigenfunction $u$ of the invariant differential operators $D\in D(G_0/K_0)$ acted upon by elements $k_0\in K_0$ and then integrating over $K_0$ to produce an eigenfunction $\phi(x)=\int_{k_0\in K_0} u^{k_0}(x)\dif\mu(k_0)$ of $D\in D(G_0/K_0)$ that is invariant under the action $f\mapsto f^{k_0}\,(k_0\in K_0)$ of the elements of $K_0$. This is called the \emph{method of images} and this is basically how one obtains Harish-Chandra's characterization of the spherical function on the symmetric space $G_0/K_0$ as the integral \eqref{harish-spherical 2}. 

However, the action $f\mapsto f^g\,(g\in G_0)$ of the elements of $G_0$ on the functions $f\colon X\rightarrow\mb{C}$ of $X$ that we have considered in this process is the one that is normally considered in case of group actions, i.e., the action $f\mapsto f^g\,(g\in G_0)$ of the elements of $G_0$ on the functions $f\colon X\rightarrow\mb{C}$, where $f^g\colon X\rightarrow \mb{C}$ is given by the assignment $x\mapsto f(g^{-1} x)$. One can instead introduce a weight factor, i.e., a function $j\colon G_0\times X\rightarrow \mb{C}$ satisfying
\begin{align}\label{cocycle condition}
    j(g_1\,g_2,x)=j(g_1,g_2 x)\,j(g_2, x)
\end{align}
and consider the action $f\mapsto f^g_j\,(h\in G_0)$ of the elements of $G_0$ on the functions $f\colon X\rightarrow\mb{C}$, where $f^g_j\colon X\rightarrow \mb{C}$ is given instead by the assignment $x\mapsto j(g^{-1},x) f(g^{-1}x)$. Then, to compute a spherical function $\phi_j$, i.e., an eigenfunction of the invariant differential operators $D\in D_j(G_0/K_0)$ that is invariant under this action of the elements of $K_0$, we must have an eigenfunction $u_j$ of $D\in D_j(G_0/K_0)$ acted upon by this action of the elements of $K_0$ and then integrate over $K_0$ to produce
\begin{align*}
    \phi_j(x)=\int_{k_0\in K_0} j({k_0}^{-1},x)\, u_j({k_0}^{-1}x) \dif\mu(k_0). 
\end{align*}

In \autoref{SEC Siegel Maass forms}, we had considered one such weighted action of symplectic matrices due to Maa\ss{} and obtained a Laplacian invariant under this action. As we eventually want to construct the heat kernel for this weight-$\kappa$ Siegel--Maa\ss{} Laplacian $\Delta^{(\kappa)}$, in this section we adapt the computation of the spherical function and the heat kernel on $\mb{H}_n$ in subsections \ref{SEC Spherical function on Hn} and \ref{SEC Heat kernel on Hn} for the weight-$\kappa$ case.

As introduced in (\ref{ycusp}), the weight-$\kappa$ action of a real symplectic matrix $g\in G_0=\mathrm{Sp}_n(\mb{R})$ on functions $f\colon\mb{H}_n\rightarrow\mb{C}$ on $\mb{H}_n$ is given by
\begin{align}\label{weight-k action}
    f^{\,g^{-1}}(Z)=j_{\kappa}(g,Z)f(g\,Z)\quad(g\in G_0,\,Z\in\mb{H}_n),
\end{align}
where the weight-factor $j_{\kappa}(g,Z)$ is given by
\begin{align*}
    j_{\kappa}(g,Z)=\bigg(\frac{\det(C\overline{Z}+D)}{\det(CZ+D)}\bigg)^{\kappa/2}\quad \bigg(g=\begin{pmatrix}A&B\\C&D\end{pmatrix}\in \mathrm{Sp}_n(\mb{R}),\,Z\in\mb{H}_n\bigg).
\end{align*}
It is easy to check that the weight-factor $j_{\kappa}$ satisfies the property \eqref{cocycle condition}.

The functions $f\colon\mb{H}_n\rightarrow\mb{C}$ on $\mb{H}_n$ can be lifted to functions $\widetilde{f}\colon \mathrm{Sp}_n(\mb{R}) \rightarrow\mb{C}$ defined by
\begin{align*}
    \widetilde{f}(g):= j_{\kappa}(g,i\mathbbm{1}_n)\,f(g\, i\mathbbm{1}_n).
\end{align*}
There is a one-to-one correspondence between the functions on $f\colon\mb{H}_n\rightarrow\mb{C}$ on $\mb{H}_n$ that satisfy 
\begin{align*}
f(Z)=j_{\kappa}(g^{\prime},Z)f(g^{\prime} Z)\quad(Z=g\, i\mathbbm{1}_n\in\mb{H}_n)
\end{align*}
for some $g^{\prime}\in G_0$ and the functions $\widetilde{f}\colon G_0\rightarrow\mb{C}$ on $G_0$ that satisfy 
\begin{center}
\begin{enumerate}
\item[(i)] $\widetilde{f}(g^{\prime}\,g)=\widetilde{f}(g)$ for all $g\in G_0$
\item[(ii)] $\widetilde{f}(g\,k_0)=j_{\kappa}(k_0,i\mathbbm{1}_n)\widetilde{f}(g)$ for all $g\in G_0$ and $k_0\in K_0$.
\end{enumerate}
\end{center}
Therefore, to compute the weight-$\kappa$ spherical function on $\mb{H}_n$, we need to integrate over the action 
\begin{align*}
    \widetilde{f}^{\,k_0^{-1}}(g)=j_{\kappa}(k_0,i\mathbbm{1}_n)^{-1}\widetilde{f}(k_0\,g)
\end{align*}
of $K_0$ on $G_0$, which takes the explicit form
\begin{align*}
    \widetilde{f}^{\,k_0^{-1}}(g)&=\bigg(\frac{\det(A+iB)}{\det(A-iB)}\bigg)^{\kappa/2}\widetilde{f}(k_0\,g)\quad \bigg(k_0=\begin{pmatrix}A&B\\-B&A\end{pmatrix}\in K_0,\,g\in \mathrm{Sp}_n(\mb{R})\bigg)
\end{align*}
when we write $k_0$ in the familiar block-matrix form for symplectic matrices.

However, this lifts the weight-$\kappa$ action of $G_0$ on $\mb{H}_n$ to $G_0$, while as our computation of the spherical function, made by reducing it to the complex case, takes place in the complex group $G$, we need to determine this action in $G$. As in the complex reduction method, we consider the Lie algebra $\mf{g}=\mf{g}_0+i\mf{g}_0$ of $G$ as a real Lie algebra, its elements are canonically embedded in the space of $(4n\times 4n)$-real matrices as
\begin{align*}
    X\mapsto \begin{pmatrix}\re(X)& \im(X)\\-\im(X)&\re(X)\end{pmatrix}\quad\quad(X\in\mf{g}).
\end{align*}
Therefore, the elements of $G$ are also canonically embedded in the space of $(4n\times 4n)$-real matrices as
\begin{align*}
    g=\exp(X)\mapsto \begin{pmatrix}\re(g)& \im(g)\\-\im(g)&\re(g)\end{pmatrix}\quad\quad(g\in G).
\end{align*}
Accordingly, the element $i\mathbbm{1}_n\in\mb{H}_n$, under this canonical embedding of $(n\times n)$-complex matrices into $(2n\times 2n)$-real matrices takes the form
\begin{align*}
    i\mathbbm{1}_n\mapsto \begin{pmatrix}0& \mathbbm{1}_n\\-\mathbbm{1}_n& 0\end{pmatrix}=J_{n}.
\end{align*}
This gives us the weight-$\kappa$ action of $K$ on the functions $\widetilde{f}\colon G\rightarrow\mb{C}$ on $G$ as
\begin{align*}
    \widetilde{f}^{\,k^{-1}}(g)&=j_{\kappa}(k,J_n)\widetilde{f}(k\,g)\quad\quad (k\in K,\,g\in \mathrm{Sp}_n(\mb{C})),
\end{align*}
where the weight-factor $j_{\kappa}$ is of the form
\begin{align*}
    j_{\kappa}(k,J_n)=\bigg(\frac{\det(\re(k)+\im(k)J_n)}{\det(\re(k)-\im(k)J_n)}\bigg)^{\kappa/2}
\end{align*}
Now writing $k\in K$ as $k=k_0\,k_h$, where $k_0$ is real orthogonal and $k_h$ is Hermitian orthogonal, by the property \eqref{cocycle condition} of $j_{\kappa}$, we have 
$j_{\kappa}(k,J_n)=j_{\kappa}(k_0,k_h J_n)j_{\kappa}(k_h,J_n)$. Since $k_0$ is real orthogonal, $\im(k_0)=0$ so that $j_{\kappa}(k_0,k_hJ_n)=1$, thereby giving $j_{\kappa}(k,J_n)=j_{\kappa}(k_h,J_n)$. 

To calculate $j_{\kappa}(k_h,J_n)$ more explicitly, we need to write $k_h$ in the block-diagonal form
\begin{align*}
    k_h=\begin{pmatrix}A&B\\-B&A\end{pmatrix}\quad\quad(AA^t+BB^t=\mathbbm{1}_n,\,AB^t=BA^t,\,A=\overline{A}^t,\,B=-\overline{B}^t). 
\end{align*}
The matrix $h:=A+iB$ is obviously $(n\times n)$-Hermitian. The orthogonality condition $AA^t+BB^t=\mathbbm{1}_n,\,AB^t=BA^t$ can be restated as
\begin{align*}
    (A+iB)(A^t-iB^t)=\mathbbm{1}_n, 
\end{align*}
so that, we have $A-iB=h^{-t}$. As $h$ is Hermitian, this also implies that $A-iB=\overline{h}^{-1}$. Then to calculate $\det(\re(k)+\im(k)J_n)$ explicitly, we have
\begin{align*}
   & \det(\re(k)+\im(k)J_n)\\
   &\quad =\det\bigg[\begin{pmatrix}\frac{1}{2}(A+\overline{A})&\frac{1}{2}(B+\overline{B})\\-\frac{1}{2}(B+\overline{B}) & \frac{1}{2}(A+\overline{A}) \end{pmatrix}+\begin{pmatrix}-\frac{i}{2}(A-\overline{A})&-\frac{i}{2}(B-\overline{B})\\\frac{i}{2}(B-\overline{B}) & -\frac{i}{2}(A-\overline{A}) \end{pmatrix}\begin{pmatrix}0& \mathbbm{1}_n\\-\mathbbm{1}_n& 0\end{pmatrix}\bigg]\\
   &\quad=\det\bigg[\begin{pmatrix}\frac{1}{2}((A+iB)+(\overline{A}-i\overline{B})) & -\frac{i}{2}((A+iB)-(\overline{A}-i\overline{B})) \\ 
   \frac{i}{2}((A+iB)-(\overline{A}-i\overline{B})) & \frac{1}{2}((A+iB)+(\overline{A}-i\overline{B})) \end{pmatrix}\bigg],
\end{align*}
whence using the relations
\begin{align}\label{hermitians}
    A+iB=h,\quad\quad A-iB=\overline{h}^{-1},\quad\quad\overline{A}+i\overline{B}=h^{-1},\quad\quad\overline{A}-i\overline{B}=\overline{h},
\end{align}
it follows that
\begin{align*}
    \det(\re(k)+\im(k)J_n)=&\det\bigg[\begin{pmatrix}\frac{1}{2}(h+\overline{h})&
    -\frac{i}{2}(h-\overline{h})\\ \frac{i}{2}(h-\overline{h}) & \frac{1}{2}(h+\overline{h}) \end{pmatrix}\bigg]\\
    =&\det\bigg[\begin{pmatrix} \re(h) & \im(h)\\ -\im(h)&\re(h))\end{pmatrix}\bigg].
\end{align*}
Now, since for any two $(n\times n)$ real matrices $X,\,Y$, we have
\begin{align*}
    \det\bigg[\begin{pmatrix} X & Y \\ -Y & X\end{pmatrix}\bigg]=\det(X)\det(X+YX^{-1}Y)=\vert \det(X+iY) \vert^2,
\end{align*}
we have here $\det(\re(k)+\im(k)J_n)=\vert\det(h)\vert^2$ and similarly, $\det(\re(k)-\im(k)J_n)=\vert\det(h^{-1})\vert^2$, thereby giving
\begin{align*}
    j_{\kappa}(k,J_n)=j_{\kappa}(k_h,J_n)=\det(h)^{2\kappa}\quad\quad
    \bigg(k_h=\begin{pmatrix}A&B\\-B&A\end{pmatrix},\,h=A+iB,\,h \text{ hermitian}\bigg).
\end{align*}
To obtain the weight-$\kappa$ spherical function on $\mb{H}_n$, we need only to multiply the integrand in \autoref{n-spherical function} in \autoref{n-spherical theorem} with this weight-factor corresponding to $k\in K$. We restate this result as a theorem for future reference.

\begin{theorem}\label{k-spherical theorem}
The spherical function $\phi^{(\kappa)}_{\lambda}$ on the Siegel upper half-space $\mb{H}_n$ for the weight-$\kappa$ action 
\begin{align*}
    f^{g^{-1}}(Z)=\bigg(\frac{\det(C\overline{Z}+D)}{\det(CZ+D)}\bigg)^{\kappa/2} f(g\,Z)\quad \bigg(g=\begin{pmatrix}A&B\\C&D\end{pmatrix}\in \mathrm{Sp}_n(\mb{R}),\,Z\in\mb{H}_n\bigg),
\end{align*}
of the symplectic matrices $g\in \mathrm{Sp}_n(\mb{R})$ on the functions $f\colon\mb{H}_n\rightarrow\mb{C}$, corresponding to
\begin{align*}
\lambda=\lambda_1\,e_1+\ldots+\lambda_j\,e_j+\ldots+\lambda_n\,e_n\in\mf{a}\spcheck
\end{align*}
at $Z=k_0\,i \exp(2R)$ with $k_0\in K_0$ and
\begin{align*}
R=\begin{pmatrix}
{r_{1}} &{} & {0} \\
{} & {\ddots}& \\
{0} &{}& {r_{n}}
\end{pmatrix}\quad\quad(r_j\in\mb{R}_{\geq 0},\;1\leq j\leq n), 
\end{align*}
is given by
\begin{align}\label{n-spherical function k}
    \phi_{\lambda}^{(\kappa)}(2R)=& \frac{c_n}{i^{n^2}\,\tau(\lambda)}
    \int\limits_{k\in K}
    \frac{\sum_{\sigma\in W}
\det(\sigma) \exp(i \sigma \lambda(\vr(r,k)))}{\delta(\vr(r,k))}\det(h(k))^{2\kappa}\dif\mu(k),
\end{align}
where $h(k)=A+iB$ is the Hermitian matrix obtained from the decomposition of $k\in K$ into real orthogonal $k_0\in K_0$ and Hermitian orthogonal 
\begin{align*}
    k_h=\begin{pmatrix}A&B\\-B&A\end{pmatrix}\quad\quad(AA^t+BB^t=\mathbbm{1}_n,\,AB^t=BA^t,\,A=\overline{A}^t,\,B=-\overline{B}^t)
\end{align*}
as $k=k_0\,k_h$ and $\vr(r,k)$ is the diagonal matrix $\vr(r,k)=\big(\begin{smallmatrix}P(r,k) & 0\\0 & -P(r,k)\end{smallmatrix}\big)$ with 
\begin{align*}
P(r,k)= \begin{pmatrix}
{\vr_{1}(r,k)} &{} & {0} \\
{} & {\ddots}& \\
{0} &{}& {\vr_{n}(r,k)}
\end{pmatrix}\;(\vr_j(r,k)\in\mb{R},\;1\leq j\leq n)
\end{align*}
related to $r=\big(\begin{smallmatrix}R & 0\\0 & -R\end{smallmatrix}\big)$ via the matrix equality $k\exp(r)\overline{k}^t=u\exp(\vr)\overline{u}^t$ with $k\in K$ and $u\in U$. The functions $\tau$ and $\delta$ are given by
\begin{align*}
    \tau(\lambda)&=\prod_{1\leq j\leq n} \operatorname{th} \Big(\frac{\lambda_{j}}{2}\pi\Big)\prod_{1\leq j<k \leq n} \operatorname{th} \Big(\frac{\lambda_{j}+\lambda_{k}}{2} \pi\Big) \prod_{1 \leq j<k \leq n}  \operatorname{th} \Big(\frac{\lambda_{j}-\lambda_{k}}{2} \pi\Big) ,\\
    \delta(\vr)&=\prod_{1\leq j\leq n} \sh(\vr_j) \prod_{1\leq j<k\leq n} \sh\Big(\frac{\vr_j+\vr_k}{2}\Big)\prod_{1\leq j<k\leq n} \sh\Big(\frac{\vr_j-\vr_k}{2}\Big),
\end{align*}
while $c_n$ is a positive real constant depending only on $n$. 
\end{theorem}
The Siegel-Maa\ss{} Laplacian $\Delta^{(\kappa)}$ is invariant under the weight-$\kappa$ action \eqref{weight-k action} of the symplectic group. This is due to the fact that the Casimir operator $\omega\in U\mf{g}_0$ descends under this action to the $\Delta^{(\kappa)}$. The only part in our computation of the heat kernel where the action of the group $G_0$ on functions on $G_0/K_0$ played a role was in the computation of the spherical function in \autoref{SEC Flensted-Jensen reduction}, which was done by integrating over the action of the complex orthogonal group $K$ on the spherical function for the complex group $G$. Therefore, to construct the heat kernel for the Laplacian $\Delta^{(\kappa)}$ on $\mb{H}_n\,(n\in\mb{N}_{\geq 1})$, we only need to adapt the formula for the spherical function by suitably altering the action of the group $K$ on the spherical function $\Phi_{\Lambda}$ for the complex group $G$, which was done in \autoref{k-spherical theorem} by multiplying the integrand in \autoref{n-spherical function} in \autoref{n-spherical theorem} with a weight-factor $\det(h(k))^{-2\kappa}$, where $h(k)=A+iB$ is the Hermitian matrix obtained from the decomposition of $k\in K$ into real orthogonal $k_0\in K_0$ and Hermitian orthogonal 
\begin{align*}
    k_h=\begin{pmatrix}A&B\\-B&A\end{pmatrix}\quad\quad(AA^t+BB^t=\mathbbm{1}_n,\,AB^t=BA^t,\,A=\overline{A}^t,\,B=-\overline{B}^t)
\end{align*}
as $k=k_0\,k_h$. To obtain a more explicit bound for the heat kernel on $\mb{H}_n$ corresponding to $\Delta^{(\kappa)}$, we need a bound for this factor $\det(h(k))$ in terms of the diagonal matrices $\vr$ and $r$, which is what we undertake next.

\begin{lemma}\label{horn-johnson lemma}
Let $A$ be a $(n\times n)$-Hermitian matrix. Let eigenvalues of $A$ be labeled according to increasing size:
\begin{align*}
\lambda_{\min}(A)=\lambda_1(A)\leq \ldots \leq\lambda_n(A)=\lambda_{\max}(A)
\end{align*}
Let $r$ be an integer with $1\leq r\leq n$, and let $A_r$ denote any $(r\times r)$-principal submatrix of $A$ obtained by deleting $n-r$ rows and the corresponding $n-r$ columns from $A$. For each integer $k$ such that $1\leq k\leq r$, we have
\begin{align*}
\lambda_k(A)\leq \lambda_k(A_r)\leq \lambda_{k+n-r}(A)
\end{align*}
\end{lemma}
\begin{proof}
See \cite[p. 189, Theorem 4.3.15]{Horn-Johnson}
\end{proof}

\begin{theorem}\label{heat kernel correction theorem}
Let $k\in K$ be a complex symplectic orthogonal matrix and $h(k)=A+iB$ be the Hermitian matrix obtained from the decomposition of $k$ into real orthogonal $k_0\in K_0$ and Hermitian orthogonal 
\begin{align*}
    k_h=\begin{pmatrix}A&B\\-B&A\end{pmatrix}\quad\quad(AA^t+BB^t=\mathbbm{1}_n,\,AB^t=BA^t,\,A=\overline{A}^t,\,B=-\overline{B}^t)
\end{align*}
as $k=k_0\,k_h$. Let $R$ be the diagonal matrix \begin{align*}
R=\begin{pmatrix}
{r_{1}} &{} & {0} \\
{} & {\ddots}& \\
{0} &{}& {r_{n}}
\end{pmatrix}\quad\quad(r_j\in\mb{R}_{\geq 0},\;1\leq j\leq n) 
\end{align*}
and $r=\big(\begin{smallmatrix}R & 0\\0 & -R\end{smallmatrix}\big)$. Let $u\exp(\vr)\overline{u}^t=k \exp(r)\overline{k}^t$ be the eigendecomposition of the Hermitian matrix $k \exp(r)\overline{k}^t$ with $u\in U$ unitary symplectic and $\vr(r,k)$ the diagonal matrix $\vr(r,k)=\big(\begin{smallmatrix}P(r,k) & 0\\0 & -P(r,k)\end{smallmatrix}\big)$ with 
\begin{align*}
P(r,k)= \begin{pmatrix}
{\vr_{1}(r,k)} &{} & {0} \\
{} & {\ddots}& \\
{0} &{}& {\vr_{n}(r,k)}
\end{pmatrix}\quad(\vr_j(r,k)\in\mb{R},\;1\leq j\leq n).
\end{align*} 
Then $\det(h(k))$ is bounded above by
\begin{align*}
    \det(h(k))\leq \frac{\exp(\sum_{j=1}^n \vert\vr_j\vert)}{\prod\limits_{j=1}^n \ch(r_j)}. 
\end{align*}
\end{theorem}

\begin{proof}
Let $l$ be the $(2n\times 2n)$-matrix
\begin{align}\label{ldef}
    l=\frac{1-i}{2}
    \begin{pmatrix}
    \mathbbm{1}_n & -i\mathbbm{1}_n\\
    \mathbbm{1}_n & i\mathbbm{1}_n
    \end{pmatrix}.
\end{align}
It is easy to check that $l$ is a symplectic unitary matrix, whose inverse is given by
\begin{align*}
    l^{-1}=\overline{l}^t=\frac{1+i}{2}
    \begin{pmatrix}
    \mathbbm{1}_n & \mathbbm{1}_n\\
    i\mathbbm{1}_n & -i\mathbbm{1}_n
    \end{pmatrix}.
\end{align*}
Also let the symplectic real orthogonal $k_0$ and Hermitian orthogonal $k_h$ be of the forms
\begin{align*}
    k_0=\begin{pmatrix}A_0&B_0\\-B_0&A_0\end{pmatrix}\quad\quad\text{and}\quad\quad k_h=\begin{pmatrix}A&B\\-B&A\end{pmatrix},
\end{align*}
respectively. Then multiplying the matrix $k\exp(r)\overline{k}^t=k_0\,k_h\,\exp(r)\,k_h\,k_0^t$ from the left by $l$ and from the right by $l^{-1}$, and writing the product as
\begin{align*}
l(k\exp(r)\overline{k}^t)l^{-1}=(l\,k_0\,l^{-1})(l\,k_h\,l^{-1})(l\,\exp(r)\,l^{-1})(l\,k_h\,l^{-1})(l\,k_0^t\, l^{-1}), 
\end{align*}
in the block decomposed form, we have
\begin{align*}
  l(k\exp(r)\overline{k}^t)l^{-1}=
  \big(\begin{smallmatrix}A_0+ iB_0 & 0\\0& A_0-iB_0\end{smallmatrix}\big)
  \big(\begin{smallmatrix}A+ iB & 0\\0& A-iB\end{smallmatrix}\big)
  \big(\begin{smallmatrix}\ch(R) & \sh(R)\\\sh(R)& \ch(R)\end{smallmatrix}\big)
  \big(\begin{smallmatrix}A+ iB & 0\\0& A-iB\end{smallmatrix}\big)
  \big(\begin{smallmatrix}A_0^t-iB_0^t & 0\\0& A_0^t+iB_0^t\end{smallmatrix}\big).
\end{align*}
Since $k_0\in K_0$ is real orthogonal, we know that the matrix $w:=A_0+iB_0$ is unitary. By the hypothesis of the theorem, $h=A+iB$ is Hermitian. In that case, in \autoref{hermitians}, we noted that $A-iB=\overline{h}^{-1}=h^{-t}$. With these notations, the above matrix equation becomes
\begin{align*}
   l(k\exp(r)\overline{k}^t)l^{-1}= \begin{pmatrix}w\,h\,\ch(R)\,h\,\overline{w}^t & w\,h\,\sh(R)\,h^{-t}\,w^t\\\overline{w}h^{-t}\sh(R)\,h\,\overline{w}^t & \overline{w}\,h^{-t}\,\ch(R)\,h^{-t}\,w^t\end{pmatrix}.
\end{align*}
Note, that the determinant of the $(1,1)$-block of the above matrix is
$\det(h)^2 \det(\ch(R))$. 

Now coming to the other side of the matrix equation $u\exp(\vr)\overline{u}^t=k\exp(r)\overline{k}^t$, as $l$ is symplectic unitary, the matrix $s=l\,u$ is also unitary. Writing $s$ in the block decomposed form
\begin{align*}
   s=\begin{pmatrix} A & B \\ -\overline{B} & \overline{A}\end{pmatrix},
\end{align*}
we write the matrix $l(u\exp(\vr)\overline{u}^t)l^{-1}=s\exp(\vr)\overline{s}^t$ as
\begin{align*}
    s \exp(\vr)\overline{s}^t &=
    \begin{pmatrix} A & B \\ -\overline{B} & \overline{A}\end{pmatrix}
    \begin{pmatrix}\exp(P) & 0 \\ 0 & \exp(-P)\end{pmatrix}\!\!
    \begin{pmatrix} \overline{A}^t & -B^t \\ 
    \overline{B}^t & A^t\end{pmatrix}\\
    &=\begin{pmatrix} A\exp(P)\overline{A}^t+B\exp(-P)\overline{B}^t & -A\exp(P)B^t+B\exp(-P)A^t \\ 
    -\overline{B}\exp(P)\overline{A}^t+\overline{A}\exp(-P)\overline{B}^t & \overline{B}\exp(P)B^t+\overline{A}\exp(-P)A^t\end{pmatrix}
\end{align*}
Comparing the determinant of the (1,1)-block of this matrix with that of  $l(k\exp(r)\overline{k}^t)l^{-1}$, we have
\begin{align*}
    \det(h)^2 \det(\ch(R))=\det(A\exp(P)\overline{A}^t+B\exp(-P)\overline{B}^t). 
\end{align*}
Let us denote by $m$ the $(2n\times 2n)$-Hermitian matrix
$s\exp(\vr)\overline{s}^t$, and by $M$ its $(n\times n)$-principal submatrix
\begin{align*}
M:=A\exp(P)\overline{A}^t+B\exp(-P)\overline{B}^t. 
\end{align*} 
Now, $m$ being a $(2n\times 2n)$-Hermitian matrix with eigenvalues $\exp(\pm\vr_1),\ldots,\exp(\pm\vr_n)$ and $M$ being the $(n\times n)$-principal submatrix of $m$, by \autoref{horn-johnson lemma}, we have
\begin{align*}
\lambda_k(m)\leq \lambda_k(M)\leq \lambda_{n+k}(m)\quad\quad(1\leq k\leq n),
\end{align*}
which implies
\begin{align*}
\lambda_1(m)\cdots \lambda_n(m)\leq \det(M) \leq \lambda_{n+1}(m)\cdots\lambda_{2n}(m).
\end{align*}

The $n$ largest eigenvalues of $m$ are $\exp(\vert\vr_1\vert) ,\ldots,\exp(\vert\vr_n\vert)$, and
The $n$ smallest eigenvalues of $m$ are $\exp(-\vert\vr_1\vert),\ldots,\exp(-\vert\vr_n\vert)$.
Therefore, we have
\begin{align*}
\lambda_1(A)\cdots \lambda_n(A)=\exp\Big(-\sum_{j=1}^n \vert\vr_j\vert\Big) \quad\text{and}\quad
\lambda_{n+1}(A)\cdots\lambda_{2n}(A)=\exp\Big(\sum_{j=1}^n \vert\vr_j\vert\Big),
\end{align*}
from where it follows that
\begin{align*}
\exp\Big(-\sum_{j=1}^n \vert\vr_j\vert\Big)\leq\det(M)\leq \exp\Big(\sum_{j=1}^n \vert\vr_j\vert\Big),
\end{align*}
thereby proving the requisite determinant-inequality.
\end{proof}

This of course provides a very useful upper bound on the heat kernel corresponding to the Siegel--Maa\ss{} Laplacian $\Delta^{(\kappa)}$, which we state as the next theorem.
\begin{theorem}\label{heat corollary}
Let $K_t^{(\kappa)}$ denote the heat kernel at a point $Z=k_0\, i\exp(2R)\in\mb{H}_n$ corresponding to the Siegel--Maa\ss{} Laplacian $\Delta^{(\kappa)}$ of weight $\kappa$  on the Siegel upper half-space with $k_0\in K_0$ and 
\begin{align*}
R=\begin{pmatrix}
{r_{1}} &{} & {0} \\
{} & {\ddots}& \\
{0} &{}& {r_{n}}
\end{pmatrix}\quad\quad(r_j\in\mb{R}_{\geq 0},\;1\leq j\leq n).
\end{align*} 
Then, subject to the above conjecture, $K_t^{(\kappa)}$ is bounded above by
\begin{align*}
   K_t^{(\kappa)}(2R)\leq c_n\frac{\exp{\big(\!-\!\sum_{j=1}^n j^2t/4\big)}}{t^{n^2+n/2}}\!\! \int\limits_{k\in K}\!\! \frac{\ve(\vr(r,k))\exp\big(\!-\!\sum_{j=1}^n \big(\vr_j(r,k)^2/t-\kappa\vert\vr_j(r,k)\vert\big)\big)}{\delta(\vr(r,k))\prod_{j=1}^n \ch^{\kappa}(r_j)} \dif\mu(k),
\end{align*}
where $\vr(r,k)$ is the diagonal matrix $\vr(r,k)=\big(\begin{smallmatrix}P(r,k) & 0\\0 & -P(r,k)\end{smallmatrix}\big)$ with 
\begin{align*}
P(r,k)= \begin{pmatrix}
{\vr_{1}(r,k)} &{} & {0} \\
{} & {\ddots}& \\
{0} &{}& {\vr_{n}(r,k)}
\end{pmatrix}\quad(\vr_j(r,k)\in\mb{R},\;1\leq j\leq n)
\end{align*}
related to $r=\big(\begin{smallmatrix}R & 0\\0 & -R\end{smallmatrix}\big)$ via the matrix equality $k \exp(r)\overline{k}^t=u\exp(\vr)\overline{u}^t$ with $k\in K$ and $u\in U$. The functions $\ve$ and $\delta$ are given by
\begin{align*}
    \varepsilon(\vr)&=\prod_{1\leq j\leq n} \vr_j\,
    \prod_{1\leq j<k\leq n} (\vr_j+\vr_k)\,\prod_{1\leq j<k\leq n} (\vr_j-\vr_k) ,\\
    \delta(\vr)&=\prod_{1\leq j\leq n} \sh(\vr_j) \prod_{1\leq j<k\leq n} \sh\Big(\frac{\vr_j+\vr_k}{2}\Big)\prod_{1\leq j<k\leq n} \sh\Big(\frac{\vr_j-\vr_k}{2}\Big),
\end{align*}
while $c_n$ is a positive real constant depending only on $n$. 
\end{theorem}
\begin{proof}
Follows immediately from Theorems \ref{k-spherical theorem}, \ref{heat kernel theorem} and \ref{heat kernel correction theorem}
\end{proof}


\section{Sup-norm bounds on average}

Let $K_t^{(\kappa)}(R(Z,W))\,(Z,W\in\mb{H}_n)$ denote the heat kernel on $\mb{H}_n$, where $R(Z,W)$ is the matrix 
\begin{align*}
R(Z,W)= \begin{pmatrix}
{r_{1}}(Z,W) &{} & {0} \\
{} & {\ddots}& \\
{0} &{}& {r_{n}(Z,W)}
\end{pmatrix}\quad(r_j(Z,W)\in\mb{R},1\leq j\leq n),
\end{align*}
with the entries $r_j(Z,W)\,(1\leq j\leq n)$ of $R(Z,W)$  related to the eigenvalues $\rho_j(Z,W)\,(1\leq j\leq n)$ of the cross-ratio matrix (see \autoref{SEC Siegel geometry}) 
\begin{align*}
\rho(Z, W)=(Z-W)(\overline{Z}-W)^{-1}(\overline{Z}-\overline{W})(Z-\overline{W})^{-1}\quad(Z,W\in\mb{H}_n)
\end{align*}
by the relation
\begin{align*}
    \exp(2r_j(Z,W))=\frac{1+\sqrt{\rho_{j}(Z,W)}}{1-\sqrt{\rho_{j}(Z,W)}}\quad\quad(1\leq j\leq n).
\end{align*}
Then the heat kernel $K_t^{(\kappa,\Gamma)}$ on $\Gamma\backslash\mb{H}_n$ is given by the $\Gamma$-periodization
\begin{align*}
    K_t^{(\kappa,\Gamma)}(Z,W):=\sum_{\gamma\in\Gamma} \det\bigg(\frac{Z-\gamma\overline{W}}{\gamma W-\overline{Z}}\bigg)^{\kappa/2}\det\bigg(\frac{C\overline{W}+D}{CW+D}\bigg)^{\kappa/2}K_t^{(\kappa)}(2R(Z,\gamma W)).
\end{align*}
We write $K_t^{(\kappa,\Gamma)}(Z):=K_t^{(\kappa,\Gamma)}(Z,Z)$ and $R^{\gamma}(Z):=R(Z,\gamma Z)$ with entries
\begin{align*}
R^{\gamma}(Z)= \begin{pmatrix}
{r^{\gamma}_{1}}(Z) &{} & {0} \\
{} & {\ddots}& \\
{0} &{}& {r^{\gamma}_{n}(Z)}
\end{pmatrix}\quad(r^{\gamma}_j(Z)\in\mb{R},1\leq j\leq n).
\end{align*}

Since $\Delta^{(\kappa)}$ is symmetric, it extends to an essentially self-adjoint linear operator acting on a dense subspace of $\mathcal{H}_{\kappa}^n(\Gamma)$. Therefore the heat kernel $K_t^{(\kappa,\Gamma)}(Z,W)$ has a spectral decomposition
\begin{align}\label{spectral decomposition}
K_t^{(\kappa,\Gamma)}(Z,W)=&\sum_{j=1}^{\infty}\exp(-\lambda_j t) \varphi_{\lambda_j}(Z)\overline{\varphi}_{\lambda_j}(W) \notag \\ 
&\quad+\sum_{\mathcal{P}\in \mathcal{C}}c_{\mathcal{P}}\int\limits_{\lambda\in\mf{a}_{\mathcal{P}}\spcheck} \exp(-(\langle \rho_{\mathcal{P}},\rho_{\mathcal{P}}\rangle+\langle\lambda,\lambda\rangle)t)\; E_{\mathcal{P}}(Z,\rho_{\mathcal{P}}+i\lambda)\;
\overline{E}_{\mathcal{P}}(W,\rho_{\mathcal{P}}+i \lambda)\;\dif \lambda
\end{align}
converging absolutely and uniformly on compacta for $t>0$. The discrete part of the spectrum given by the first sum runs over the eigenvalues $\lambda_j$ of the Siegel--Maa\ss{} Laplacian $\Delta^{(\kappa)}$ with eigenfunctions $\varphi_{\lambda_j}$. The continuous part of the spectrum given by the second sum runs over the set $\mathcal{C}$ of inequivalent chains of rational boundary components of $M$ with $c_{\mathcal{P}}$ denoting a positive constant depending on the cusp $\mathcal{P}\in\mathcal{C}$, $\mf{a}_{\mathcal{P}}$ the Lie algebra of the diagonal component $A_{\mathcal{P}}$ of $\mathcal{P}$, $\rho_{\mathcal{P}}$ the half-sum of positive roots with multiplicity in  $\mf{a}_{\mathcal{P}}$ and $E_{\mathcal{P}}$ the Eisenstein series attached to the cusp $\mathcal{P}$. Setting $Z=W$ in equation $\eqref{spectral decomposition}$, we obtain
\begin{align*}
    K_t^{(\kappa,\Gamma)}(Z)=\sum_{j=1}^{\infty}\exp(-\lambda_j t)\; \vert \varphi_{\lambda_j}(Z)\vert^2+\sum_{\mathcal{P}\in \mathcal{C}}c_{\mathcal{P}}\int\limits_{\lambda\in\mf{a}_{\mathcal{P}}\spcheck}\exp(-(\langle \rho_{\mathcal{P}},\rho_{\mathcal{P}}\rangle+\langle\lambda,\lambda\rangle)t)\;\vert E_{\mathcal{P}}(Z,\rho_{\mathcal{P}}+i\lambda)\vert^2 \; \dif \lambda
\end{align*}

Now, let $\kappa\geq n+1$ and multiply both sides of the above equation by $\exp((n\kappa/4)((n+1)-\kappa)t)$. Then
\begin{align*}
    \frac{n\kappa}{4}((n+1)-\kappa)-\langle \rho_F,\rho_F\rangle-\langle\lambda,\lambda\rangle<0.
\end{align*}
Also, since $\lambda_j\geq (n\kappa/4)((n+1)-\kappa)$ by \autoref{kernel connection}, we have 
\begin{align*}
\frac{n\kappa}{4}((n+1)-\kappa)-\lambda_j\leq 0. 
\end{align*}
Therefore, on taking limit $t\rightarrow\infty$ on both sides of the above equation, on the right-hand side only the $\varphi_{\lambda_j}$'s corresponding to
$\lambda_j= (n\kappa/4)((n+1)-\kappa)$ survive. By \autoref{kernel connection}, these are of the form $\varphi_{\lambda_j}(Z)=\det(Y)^{\kappa/2}f_j(Z)$. Therefore, we have
\begin{align}\label{heat equivalence}
\lim\limits_{t\rightarrow\infty}\exp\big(-\frac{n\kappa}{4}(\kappa-(n+1))\,t\big) \,K_{t}^{(\kappa,\Gamma)}(Z)=\sum_{j=1}^{d}(\det Y)^{\kappa}\vert{f_j(Z)}\vert^2 \quad(\kappa>(n+1)),
\end{align}    
where $d=\dim(\mathcal{S}_{\kappa}^n(\Gamma))$ and $\{f_j\}_{1\leq j\leq d}$ is an orthonormal  basis of $\mathcal{S}_{\kappa}^n(\Gamma)$ with respect to the Petersson inner product. We denote
\begin{align*}
    S_{\kappa}^{\Gamma}(Z):=\sum_{j=1}^{d}\det (Y)^{\kappa}\vert{f_j(Z)}\vert^2\quad(Z\in\mb{H}_n).
\end{align*}
Thus, we have 
\begin{align}\label{limit connection}
    S_{\kappa}^{\Gamma}(Z)=\lim\limits_{t\rightarrow\infty}\exp\big(-\frac{n\kappa}{4}(\kappa-(n+1))\,t\big) \,K_{t}^{(\kappa,\Gamma)}(Z).
\end{align}
Since the function $\exp(-n\kappa(\kappa-(n+1))\,t/4) K_{t}^{(\kappa,\Gamma)}(Z)$ is monotonically decreasing for any $t>0$, we also have
\begin{align*}
    S_{\kappa}^{\Gamma}(Z)\leq \exp\big(-\frac{n\kappa}{4}(\kappa-(n+1))\,t\big) \,K_{t}^{(\kappa,\Gamma)}(Z).
\end{align*}
Further, as 
\begin{align*}
\Bigg\vert\det\bigg(\frac{Z-\gamma\overline{Z}}{\gamma Z-\overline{Z}}\bigg)^{\kappa/2}\det\bigg(\frac{C\overline{Z}+D}{CZ+D}\bigg)^{\kappa/2}\Bigg\vert=1, 
\end{align*}
this also implies that for any $t>0$ and $Z\in\mb{H}_n$, we have
\begin{align}\label{heat connection}
    S_{\kappa}^{\Gamma}(Z)\leq \exp\big(-\frac{n\kappa}{4}(\kappa-(n+1))\,t\big)\sum_{\gamma\in\Gamma}K_t^{(\kappa)}(2R^{\gamma}(Z)).
\end{align}


\subsection{Sup-norm bounds in the cocompact setting}

Note that, to make the calculations less cumbersome, we continue clubbing all positive real constants depending only on $n$ under the generic symbol $c_n$. 

\begin{lemma}\label{coordinate change lemma 1}
Let $G$ denote the complex symplectic group $\mathrm{Sp}_n(\mb{C})$, $K\subsetneq G$ denote the complex orthogonal group $K=\{k\in G\;\vert\;kk^t=\mathbbm{1}_{2n}\}$ and $X=G/U$ denote the symmetric space $X=\{x=g\overline{g}^t\;\vert\;g\in G\}$. Then the invariant volume form $\dif\mu(x)$ on $X$ in the coordinates $x=k\exp(r)\overline{k}^t\,(x\in X)$, where $r$ is given by the diagonal matrix $r=\big(\begin{smallmatrix}R & 0\\0 & -R\end{smallmatrix}\big)$ with 
\begin{align*}
R=\begin{pmatrix}
{r_{1}} &{} & {0} \\
{} & {\ddots}& \\
{0} &{}& {r_{n}}
\end{pmatrix}\quad\quad(r_j\in\mb{R},\;1\leq j\leq n),
\end{align*}
is given by
\begin{align*}
    \dif\mu(x)=c_n\, \vert\delta(2r)\vert\,\bigwedge\limits_{j=1}^n\dif r_j\wedge \dif\mu(k),
\end{align*}
where $\dif\mu(k)$ denotes the Haar measure on $K$, $\delta(r)$ denotes the function 
\begin{align*}
\delta(r)&=\prod_{1\leq j\leq n} \sh(r_j) \prod_{1\leq j<k\leq n} \sh\Big(\frac{r_j+r_k}{2}\Big) \prod_{1\leq j<k\leq n} \sh\Big(\frac{r_j-r_k}{2}\Big)
\end{align*}
on $R$ and $c_n$ is a constant depending only on $n$. 
\end{lemma}
\begin{proof}
The tangent space of $X=G/U$ at identity is given by the space $\mf{p}$ of real dimension $n(2n+1)$. 
Therefore, to calculate the invariant volume form $\dif\mu(x)$ at $x\in X=G/U$, we first calculate the invariant matrix differential form $x^{-1}\dif x\in\mf{p}\spcheck$. Then, for a choice of dual basis $e_1,\ldots,e_j,\ldots,e_{n(2n+1)}$ of $\mf{p}\spcheck$, we have
\begin{align*}
x^{-1}\dif x=\omega_1(x) e_1+\ldots+ \omega_j(x) e_j+\ldots+\omega_{n(2n+1)}(x) e_{n(2n+1)},
\end{align*}
where each $\omega_j(x)\,(1\leq j \leq n(2n+1))$ is a real $1$-form. The volume form $\dif\mu(x)$, denoted by $[x^{-1}\dif x]$ is then obtained by taking the wedge product
\begin{align*}
 [x^{-1}\dif x]=\omega_1(x) \wedge \cdots \omega_j(x) \cdots+\omega_{n(2n+1)}(x). 
\end{align*}

From $x=k\exp(r)\overline{k}^t\,(x\in X)$, one obtains
\begin{align*}
    x^{-1}\dif x &=\overline{k}e^{-r} k^t\dif k e^r\overline{k}^t+\overline{k}\dif r\overline{k}^t+\overline{k}\dif\overline{k}^t\\
    &=\overline{k}e^{-r/2}\big(e^{-r/2} (k^t\dif k) e^{r/2}+e^{r/2}(\dif r) e^{-r/2}+e^{r/2}(\dif\overline{k}^t\overline{k})e^{-r/2}\big)e^{r/2}\overline{k}^t.
\end{align*}
Then taking volume form, denoted by the parentheses $[\;\cdot\;]$, on both sides, we have
\begin{align}\label{localdif1}
   \dif\mu(x)= [x^{-1}\dif x]=[e^{-r/2} (k^t\dif k) e^{r/2}+e^{r/2}(\dif r) e^{-r/2}+e^{r/2}(\dif\overline{k}^t\overline{k})e^{-r/2}]. 
\end{align}
Now, as the invariant differential form $k^t\dif k$ has the structure of the elements of the complex orthogonal Lie algebra $\mf{k}$, we take
\begin{align*}
    k^t\dif k=\begin{pmatrix} A & B \\ -B & A\end{pmatrix}\quad\quad (A,B\in\mathbb{C}^{n\times n},\,A=-A^t,\,B=B^t). 
\end{align*}
Then the form $\dif\overline{k}^t\overline{k}=\overline{k^t\dif k}^t$ is given by
\begin{align*}
    \dif\overline{k}^t\overline{k}=\begin{pmatrix} \overline{A}^t & -\overline{B}^t \\ \overline{B}^t & \overline{A}^t\end{pmatrix}=\begin{pmatrix} -\overline{A} & -\overline{B} \\ \overline{B} & -\overline{A}\end{pmatrix}.
\end{align*}
Then writing the right-hand side of \autoref{localdif1} in block decomposed form, we have
\begin{align*}
    \dif\mu(x)=&\Bigg[\begin{pmatrix}e^{-R/2}&0\\0&e^{R/2}\end{pmatrix}\begin{pmatrix} A & B \\ -B & A\end{pmatrix}\begin{pmatrix}e^{R/2}&0\\0&e^{-R/2}\end{pmatrix}\\
    &\quad+\begin{pmatrix}e^{R/2}&0\\0&e^{-R/2}\end{pmatrix}
    \begin{pmatrix} -\overline{A} & -\overline{B} \\ \overline{B} &-\overline{A}\end{pmatrix}\begin{pmatrix}e^{-R/2}&0\\0&e^{R/2}\end{pmatrix}
    +\begin{pmatrix}\dif R&0\\0&-\dif R\end{pmatrix}\Bigg]\\
    =&\Bigg[ \begin{pmatrix}
    e^{-R/2}Ae^{R/2}-e^{R/2}\overline{A}e^{-R/2} & e^{-R/2}Be^{-R/2}-e^{R/2}\overline{B}e^{R/2}\\
    -e^{R/2}Be^{R/2}+e^{-R/2}\overline{B}e^{-R/2}&
    e^{R/2}Ae^{-R/2}-e^{-R/2}\overline{A}e^{R/2}
    \end{pmatrix}
    +\begin{pmatrix}\dif R&0\\0&-\dif R\end{pmatrix}\Bigg].
\end{align*}
Now, writing the matrices $A$ and $B$ as
\begin{align*}
    A&=\big(\alpha_{j,k}=\xi_{j,k}+i\eta_{j,k}\big)_{1\leq j,k\leq n}\quad\quad\big(\alpha_{j,j}=0,\,\alpha_{k,j}=-\alpha_{j,k}\;(1\leq j<k\leq n)\big),\\
    B&=\big(\beta_{j,k}=\omega_{j,k}+i\tau_{j,k}\big)_{1\leq j,k\leq n}\quad\quad\big(\beta_{k,j}=\beta_{j,k}\;(1\leq j\leq k \leq n)\big),
\end{align*}
where $\xi_{j,k},\,\eta_{j,k},\,\omega_{j,k},\,\tau_{j,k}$ are real $1$-forms, one obtains
\begin{align*}
    \big(e^{-R/2}Ae^{R/2}-e^{R/2}\overline{A}e^{-R/2}\big)_{j,k}
    \!\!=&e^{(r_k-r_j)/2}(\xi_{j,k}+i\eta_{j,k})-e^{(r_j-r_k)/2}(\xi_{j,k}-i\eta_{j,k})\\
    =&-2 \sh\Big(\frac{r_j-r_k}{2}\Big)\xi_{j,k}+2i \ch\Big(\frac{r_j-r_k}{2}\Big)\eta_{j,k}\,(1\leq j< k\leq n)
\end{align*}
and similarly
\begin{align*}
    &\big(e^{R/2}Ae^{-R/2}-e^{-R/2}\overline{A}e^{R/2}\big)_{j,k}\!\!=2 \sh\Big(\frac{r_j-r_k}{2}\Big)\xi_{j,k}+2i \ch\Big(\frac{r_j-r_k}{2}\Big)\eta_{j,k}\,\,(1\leq j< k\leq n),\\
    &\big(e^{-R/2}Be^{-R/2}-e^{R/2}\overline{B}e^{R/2}\big)_{j,k}\!\!=
    \!-2 \sh\Big(\frac{r_j+r_k}{2}\Big)\omega_{j,k}
    +2i \ch\Big(\frac{r_j+r_k}{2}\Big)\tau_{j,k}\,(1\leq j\leq k\leq n),\\
    &\big(e^{R/2}Be^{R/2}-e^{-R/2}\overline{B}e^{-R/2}\big)_{j,k}\!\!=2 \sh\Big(\frac{r_j+r_k}{2}\Big)\omega_{j,k}
    +2i \ch\Big(\frac{r_j+r_k}{2}\Big)\tau_{j,k}\,\,(1\leq j\leq k\leq n).
\end{align*}
Now taking wedge product of the above entries, it easily follows that
\begin{align*}
\dif\mu(x)= c_n \,\vert\delta(2r)\vert\,\bigwedge\limits_{j=1}^n\dif r_j\bigwedge\limits_{1\leq j<k\leq n}(\xi_{j,k}\wedge\eta_{j,k})\bigwedge\limits_{1\leq j\leq k\leq n}(\omega_{j,k}\wedge\tau_{j,k}),
\end{align*}
whence, identifying 
\begin{align*}
    \dif\mu(k)=\bigwedge\limits_{1\leq j<k\leq n}(\xi_{j,k}\wedge\eta_{j,k})\bigwedge\limits_{1\leq j\leq k\leq n}(\omega_{j,k}\wedge\tau_{j,k})
\end{align*}
we have the result stated in the lemma.
\end{proof}

\begin{lemma}\label{coordinate change lemma 2}
Let $G$ denote the complex symplectic group $\mathrm{Sp}_n(\mb{C})$, $U\subsetneq G$ denote the symplectic unitary group $U=\{u\in G\;\vert\;u\overline{u}^t=\mathbbm{1}_{2n}\}$ and $X=G/U$ denote the symmetric space $X=\{x=g\overline{g}^t\;\vert\;g\in G\}$. Then the invariant volume form $\dif\mu(x)$ on $X$ in the coordinates $x=u\exp(\vr)\overline{u}^t\,(x\in X)$, where $\vr$ is given by the diagonal matrix $\vr=\big(\begin{smallmatrix}P & 0\\0 & -P\end{smallmatrix}\big)$ with 
\begin{align*}
P=\begin{pmatrix}
{\vr_{1}} &{} & {0} \\
{} & {\ddots}& \\
{0} &{}& {\vr_{n}}
\end{pmatrix}\quad\quad(\vr_j\in\mb{R},\;1\leq j\leq n),
\end{align*}
is given by
\begin{align*}
    \dif\mu(x)=c_n\, \delta(\vr)^2\,\bigwedge\limits_{j=1}^n\dif \vr_j\wedge \dif\mu(u),
\end{align*}
where $\dif\mu(u)$ denotes the Haar measure on $U$, $\delta(\vr)$ denotes the function 
\begin{align*}
\delta(\vr)&=\prod_{1\leq j\leq n} \sh(\vr_j) \prod_{1\leq j<k\leq n} \sh\Big(\frac{\vr_j+\vr_k}{2}\Big) \prod_{1\leq j<k\leq n} \sh\Big(\frac{\vr_j-\vr_k}{2}\Big)
\end{align*}
on $P$ and $c_n$ is a constant depending only on $n$. 
\end{lemma}
\begin{proof}
Proceeding as in the proof of \autoref{coordinate change lemma 1}, from $x=u\exp(\vr)\overline{u}^t$, one obtains
\begin{align*}
   x^{-1}\dif x&=ue^{-\vr}\overline{u}^t\dif u e^{\vr}\overline{u}^t+u\dif\vr\overline{u}^t+u\dif\overline{u}^t \\
   &=ue^{-\vr/2}\big(e^{-\vr/2}(\overline{u}^t\dif u)e^{\vr/2}
   +e^{\vr/2}(\dif \rho)e^{-\vr/2}+
   e^{\vr/2}(\dif \overline{u}^t u)e^{-\vr/2}\big)e^{\vr/2}\overline{u}^t.
\end{align*}
Then taking volume form, denoted by the parentheses $[\;\cdot\;]$, on both sides, we have
\begin{align}\label{localdif2}
   \dif\mu(x)= [x^{-1}\dif x]=[e^{-\vr/2}(\overline{u}^t\dif u)e^{\vr/2}
   +e^{\vr/2}(\dif \rho)e^{-\vr/2}+
   e^{\vr/2}(\dif \overline{u}^t u)e^{-\vr/2}]. 
\end{align}
Now, as the invariant differential form $\overline{u}^t\dif u$ has the structure of the elements of the unitary symplectic Lie algebra $\mf{u}$, we take
\begin{align*}
   \overline{u}^t\dif u= \begin{pmatrix} A & B \\ -\overline{B} & \overline{A}\end{pmatrix} \quad\quad (A,B\in\mathbb{C}^{n\times n},\, \overline{A}^t=-A,\,B=B^t).
\end{align*}
Also, from $\overline{u}^t u=\mathbbm{1}_{2n}$ it follows that that $\dif \overline{u}^t u=-\overline{u}^t \dif u$. Then writing the right-hand side of \autoref{localdif2} in block decomposed form, we have
\begin{align*}
    \dif\mu(x)=&\Bigg[\begin{pmatrix}e^{-P/2}&0\\0&e^{P/2}\end{pmatrix}\begin{pmatrix} A & B \\ -\overline{B} & \overline{A}\end{pmatrix}\begin{pmatrix}e^{P/2}&0\\0&e^{-P/2}\end{pmatrix}\\
    &\quad-\begin{pmatrix}e^{P/2}&0\\0&e^{-P/2}\end{pmatrix}
    \begin{pmatrix} A & B \\ -\overline{B} & \overline{A}\end{pmatrix}\begin{pmatrix}e^{-P/2}&0\\0&e^{P/2}\end{pmatrix}
    +\begin{pmatrix}\dif P&0\\0&-\dif P\end{pmatrix}\Bigg]\\
    =&\Bigg[ \begin{pmatrix}
    e^{-P/2}Ae^{P/2}-e^{P/2}Ae^{-P/2} & e^{-P/2}Be^{-P/2}-e^{P/2}Be^{P/2}\\
    -e^{P/2}\overline{B}e^{P/2}+e^{-P/2}\overline{B}e^{-P/2}&
    e^{P/2}\overline{A}e^{-P/2}-e^{-P/2}\overline{A}e^{P/2}
    \end{pmatrix}
    +\begin{pmatrix}\dif P&0\\0&-\dif P\end{pmatrix}\Bigg].
\end{align*}
Now, writing the matrices $A$ and $B$ as
\begin{align*}
    A&=\big(\alpha_{j,k}\big)_{1\leq j,k\leq n}\quad\quad\big(\re(\alpha_{j,j})=0,\,\alpha_{k,j}=-\overline{\alpha_{j,k}}\;(1\leq j<k\leq n)\big),\\
    B&=\big(\beta_{j,k}\big)_{1\leq j,k\leq n}\quad\quad\big(\beta_{k,j}=\beta_{j,k}\;(1\leq j\leq k \leq n)\big),
\end{align*}
where $\alpha_{j,k},\,\beta_{j,k}$ are complex $1$-forms, one obtains
\begin{align*}
    &\big(e^{-P/2}Ae^{P/2}-e^{P/2}Ae^{-P/2}\big)_{j,k}=-2\sh\Big(\frac{\vr_j-\vr_k}{2}\Big) \,\alpha_{j,k}\quad\quad&(1\leq j<k\leq n),\\
    &\big(e^{P/2}\overline{A}e^{-P/2}-e^{-P/2}\overline{A}e^{P/2}\big)_{j,k}=2\sh\Big(\frac{\vr_j-\vr_k}{2}\Big) \,\overline{\alpha}_{j,k}\quad\quad&(1\leq j<k\leq n),\\
    &\big(e^{-P/2}Be^{-P/2}-e^{P/2}Be^{P/2}\big)_{j,k}=-2\sh\Big(\frac{\vr_j+\vr_k}{2}\Big) \,\beta_{j,k}\quad\quad&(1\leq j\leq k\leq n),\\
    &\big(e^{P/2}\overline{B}e^{P/2}+e^{-P/2}\overline{B}e^{-P/2}\big)_{j,k}=2\sh\Big(\frac{\vr_j+\vr_k}{2}\Big) \,\overline{\beta}_{j,k}\quad\quad&(1\leq j\leq k\leq n).
\end{align*}
Now taking wedge product of the above entries, it easily follows that
\begin{align*}
\dif\mu(x)= c_n \,\delta(\vr)^2\,\bigwedge\limits_{j=1}^n\dif \vr_j\bigwedge\limits_{1\leq j<k\leq n}(\alpha_{j,k}\wedge\overline{\alpha}_{j,k})\bigwedge\limits_{1\leq j\leq k\leq n}(\beta_{j,k}\wedge\overline{\beta}_{j,k}),
\end{align*}
whence, identifying 
\begin{align*}
    \dif\mu(u)=\bigwedge\limits_{1\leq j<k\leq n}(\alpha_{j,k}\wedge\overline{\alpha}_{j,k})\bigwedge\limits_{1\leq j\leq k\leq n}(\beta_{j,k}\wedge\overline{\beta}_{j,k})
\end{align*}
we have the result stated in the lemma.
\end{proof}

\begin{theorem}\label{compact case}
For any arithmetic subgroup $\Gamma\subsetneq \mathrm{Sp}_n(\mb{R})$ such that $M:=\Gamma\backslash\mb{H}_n$ is compact, we have
\begin{align*}
   \sup_{Z\in\mb{H}_n} S_{\kappa}^{\Gamma}(Z)\leq c_{n,\Gamma}\, \kappa^{n(n+1)/2}\quad\quad( \kappa\geq n+1),
\end{align*}
where $c_{n,\Gamma}$ is a positive real constant depending only on $n$ and $\Gamma$.
\end{theorem}
\begin{proof}
For $Z,\,W\in\mb{H}_n$, let $R(Z,W)$ denote the matrix 
\begin{align*}
R(Z,W)= \begin{pmatrix}
{r_{1}}(Z,W) &{} & {0} \\
{} & {\ddots}& \\
{0} &{}& {r_{n}(Z,W)}
\end{pmatrix}\quad(r_j(Z,W)\in\mb{R},1\leq j\leq n),
\end{align*}
with the entries $r_j(Z,W)\,(1\leq j\leq n)$ of $R(Z,W)$  related to the eigenvalues $\rho_j(Z,W)\,(1\leq j\leq n)$ of the cross-ratio matrix  
\begin{align*}
\rho(Z, W)=(Z-W)(\overline{Z}-W)^{-1}(\overline{Z}-\overline{W})(Z-\overline{W})^{-1}\quad(Z,W\in\mb{H}_n)
\end{align*}
by the relation
\begin{align*}
    \exp(2r_j(Z,W))=\frac{1+\sqrt{\rho_{j}(Z,W)}}{1-\sqrt{\rho_{j}(Z,W)}}\quad\quad(1\leq j\leq n).
\end{align*}
Let $R(Z)$ denote the matrix $R(Z,i\mathbbm{1}_n)$ with corresponding diagonal entries $r_j(Z)\,(1\leq j\leq n)$ and for $\gamma\in\Gamma$, let $R^{\gamma}(Z)$ denote the matrix $R(Z,\gamma Z)$ with corresponding diagonal entries $r_j^{\gamma}(Z)\,(1\leq j\leq n)$.

Now, since $M=\Gamma\backslash\mb{H}_n$ is compact, there are only finitely many elements $\gamma\in\Gamma$, namely, the torsion elements of $\Gamma$, for which the point $\gamma Z$ can get arbitrarily close to $Z$. Then, denoting the set of torsion elements of $\Gamma$ by $\Gamma_T$, there is a positive real constant $c_{n,\Gamma}$, such that $R^{\gamma}(Z)\geq c_{n,\Gamma}\,\mathbbm{1}_n$ for all $\gamma\in\Gamma\setminus\Gamma_T$ and $Z\in\mb{H}_n$. Therefore, given $n$ positive real numbers $r_j\,(1\leq j\leq n)$, we have
\begin{align*}
    \#\{\gamma\in\Gamma\;\vert\;r_j^{\gamma}(Z)\leq r_j,\,1\leq j\leq n\}\leq c_{n,\Gamma}\, \mathrm{vol}_n(\{Z\in\mb{H}_n\;\vert\;r_j(Z)\leq r_j\,\,1\leq j\leq n\}),
\end{align*}
for some positive real positive constant $c_{n,\Gamma}$ depending only on $n$ and $\Gamma$. The dependence on $\Gamma$ here is given by the maximal order of the torsion elements of $\Gamma$. 

As the volume form on $\mb{H}_n$ in polar coordinates is given by 
\begin{align*}
    \dif\mu_n(k_0\,\exp(2R)i)
    =\vert\delta(2r)\vert\bigwedge\limits_{j=1}^n\dif r_j\wedge \dif\mu(k_0)
\end{align*}
with $R=\mathrm{diag}(r_1,\,r_2,\ldots,\,r_n)$, $r=\big(\begin{smallmatrix}R & 0\\0 & -R\end{smallmatrix}\big)$ and $k_0\in K_0=\mathrm{Sp}_n(\mb{R})\cap \mathrm{O}(2n,\mb{R})$, we have
\begin{align*}
    \dif\mathrm{vol}_n(\{Z\in\mb{H}_n\;\vert\;r_j(Z)\leq r_j\,,\,1\leq j\leq n\})=\vert\delta(2r)\vert\bigwedge\limits_{j=1}^n\dif r_j.
\end{align*}
Therefore, as the heat kernel $K_t^{(\kappa)}(2R)$ is non-negative, continuous, and monotonically decreasing in each $r_j\,(1\leq j\leq n)$, we have 
\begin{align*}
    \sum_{\gamma\in\Gamma}K_t^{(\kappa)}(R^{\gamma}(Z))\leq c_{n,\Gamma}
    \int\limits_{r_1=0}^{\infty}\ldots\int\limits_{r_n=0}^{\infty}K_t^{(\kappa)}(2R)
    \,\vert\delta(2r)\vert \bigwedge\limits_{j=1}^n\dif r_j.
\end{align*}
Hence, from \autoref{heat connection} and \autoref{heat corollary}, we have
\begin{align}\label{basic inequality}
        S_{\kappa}^{\Gamma}(Z)\leq&\, c_{n,\Gamma}\, I_n(\kappa,t),
\end{align}
where the function $I_n(\kappa,t)$ is given by the integral  
\begin{align*} 
   &I_n(\kappa,t):=\frac{\exp{\big(\big(-n\kappa(\kappa-(n+1))-\sum_{j=1}^n j^2\big)t/4\big)}}{t^{n^2+n/2}} \\
   &\int\limits_{r_1=0}^{\infty}\ldots\int\limits_{r_n=0}^{\infty}\; \int\limits_{k\in K} \frac{\ve(\vr(r,k))\exp\big(-\sum_{j=1}^n\,\big(\vr_j(r,k)^2/t-\kappa\vert\vr_j(r,k)\vert\big)\big)}{\delta(\vr(r,k))\prod_{j=1}^n \ch^\kappa(r_j)}\vert\delta(2r)\vert \bigwedge\limits_{j=1}^n\dif r_j\bigwedge \dif\mu(k),
\end{align*}
    
Now, using Lemmas \ref{coordinate change lemma 1} and \ref{coordinate change lemma 2}, we switch the above integral from over $r,\,k\,(k\in K)$ coordinates on $X=G/U$ to $\vr,\,u\,(u\in U)$ coordinates on $X$. As we have
\begin{align*}
    \vert\delta(2r)\vert\,\bigwedge\limits_{j=1}^n\dif r_j\wedge \dif\mu(k)=c_n\, \delta(\vr)^2\,\bigwedge\limits_{j=1}^n\dif \vr_j\wedge \dif\mu(u), 
\end{align*}
with this change of variables, the above integral becomes
\begin{align*}
    I_n(\kappa,t)&=c_n\frac{\exp{\big(\big(-n\kappa(\kappa-(n+1))-\sum_{j=1}^n j^2\big)t/4\big)}}{t^{n^2+n/2}} \\
    &\quad\int\limits_{\vr_1=-\infty}^{\infty}\ldots\int\limits_{\vr_n=-\infty}^{\infty}\; \int\limits_{u\in U} \frac{\ve(\vr)\exp\big(-\sum_{j=1}^n\,\big(\vr_j^2/t-\kappa\vert\vr_j\vert\big)\big)}{\prod_{j=1}^n \ch^\kappa(r_j)}\delta(\vr) \bigwedge\limits_{j=1}^n\dif \vr_j\bigwedge \dif\mu(u)\\
    &=c_n\frac{\exp{\big(\big(-n\kappa(\kappa-(n+1))-\sum_{j=1}^n j^2\big)t/4\big)}}{t^{n^2+n/2}} \\
    &\quad\int\limits_{\vr_1=-\infty}^{\infty}\ldots\int\limits_{\vr_n=-\infty}^{\infty} \ve(\vr)\delta(\vr)\exp\big(-\sum_{j=1}^n\,\big(\vr_j^2/t-\kappa\vert\vr_j\vert\big)\big) J_n(\vr,\kappa) \bigwedge\limits_{j=1}^n\dif \vr_j,
\end{align*}
 where $J_n(\vr,\kappa)$ is the integral given by  
 \begin{align*}
 J_n(\vr,\kappa):=\int\limits_{u\in U}\frac{\dif\mu(u)}{\prod_{j=1}^n \ch^\kappa(r_j)}.
 \end{align*}
 Now, as $1/\ch(r_j)\leq 2\exp(-r_j)\,(1\leq j\leq n)$ and $\kappa\geq n+1$ we have
 \begin{align*}
    J_n(\vr,\kappa)\leq c_n \int\limits_{u\in U} \exp\Big(-\kappa\sum_{j=1}^n r_j\Big)\dif\mu(u)\leq  c_n \int\limits_{u\in U} \exp\Big(-(n+1)\sum_{j=1}^n r_j\Big)\dif\mu(u). 
 \end{align*}
 Then, as $r_j\in\mb{R}_{\geq 0}\,(1\leq j\leq n)$, from \autoref{helgason lemma}, we have
 \begin{align*}
     \phi_0(\exp(r))\geq \exp(-\rho(r))= \exp(-nr_1-(n-1)r_2-\ldots-r_n)\geq \exp\Big(-(n+1)\sum_{j=1}^n r_j\Big),
 \end{align*}
 where $\phi_0$ is the real spherical function on $\mb{H}_n$ corresponding to $\lambda=0\in\mf{a}\spcheck$. Then, we have
 \begin{align*}
     J_n(\vr,\kappa)\leq c_n \int\limits_{u\in U} \phi_0(\exp(r))\dif\mu(u)= c_n \,\Phi_0(\exp(\vr)),
 \end{align*}
 where $\Phi_0(\exp(\vr))$ is the complex spherical function on $X=G/U$ corresponding to $\lambda=0$. Now, from \autoref{0-complex spherical}, we have
 \begin{align*}
     J_n(\vr,\kappa)\leq c_n \frac{\varepsilon(\vr)}{\delta(2\vr)}.
 \end{align*}
 Putting $J_n(\vr,\kappa)$ back in $I_n(\kappa,t)$, we have
 \begin{align*}
     I_n(\kappa,t)&\leq c_n\frac{\exp{\big(\big(-n\kappa(\kappa-(n+1))-\sum_{j=1}^n j^2\big)t/4\big)}}{t^{n^2+n/2}} \\
        &\quad\int\limits_{\vr_1=-\infty}^{\infty}\ldots\int\limits_{\vr_n=-\infty}^{\infty} \frac{\ve(\vr)^2\exp\big(-\sum_{j=1}^n\,\big(\vr_j^2/t-\kappa\vert\vr_j\vert\big)\big)}{\nu(\vr)} \bigwedge\limits_{j=1}^n\dif \vr_j,
 \end{align*}
 where $\nu(\vr)$ is the function given by
 \begin{align*}
     \nu(\vr)=\prod_{1\leq j\leq n} \ch(\vr_j) \prod_{1\leq j< k\leq n} \ch\Big(\frac{\vr_j+\vr_k}{2}\Big) \prod_{1\leq j< k \leq n} \ch\Big(\frac{\vr_j-\vr_k}{2}\Big).
 \end{align*}
Now, since the integrand in the above integral is an even function of each $\vr_j\,(1\leq j\leq n)$, we can integrate each $\vr_j$ in the limit $\vr_j\in[0,\infty]$, thereby giving
 \begin{align*}
     I_n(\kappa,t)&\leq c_n\frac{\exp{\big(\big(-n\kappa(\kappa-(n+1))-\sum_{j=1}^n j^2\big)t/4\big)}}{t^{n^2+n/2}} \\
        &\quad\int\limits_{\vr_1=0}^{\infty}\ldots\int\limits_{\vr_n=0}^{\infty}\; \frac{\ve(\vr)^2\exp\big(-\sum_{j=1}^n\,\big(\vr_j^2/t-\kappa\vr_j\big)\big)}{\nu(\vr)} \bigwedge\limits_{j=1}^n\dif \vr_j
 \end{align*}
 Now, as $1/\nu(\vr)\leq c_n \exp(-n\vr_1-(n-1)\vr_2-\ldots-\vr_n) $ and
 $n\kappa(\kappa-(n+1))+\sum_{j=1}^n j^2=(\kappa-n)^2+\ldots+(\kappa-1)^2$,
 we have
 \begin{align*}
     I_n(\kappa,t)
     &\leq c_n\int\limits_{\vr_1=0}^{\infty}\ldots\int\limits_{\vr_n=0}^{\infty}\;\frac{\ve(\vr)^2 \exp\big(-\sum_{j=1}^n(\vr_j/\sqrt{t}-(\kappa-(n-j+1)\sqrt{t}/2)^2\big)}{t^{n^2+n/2}}
     \bigwedge\limits_{j=1}^n\dif \vr_j\\
     & \leq c_n\int\limits_{\vr_1=-\infty}^{\infty}\ldots\int\limits_{\vr_n=-\infty}^{\infty}\;\frac{\ve(\vr)^2 \exp\big(-\sum_{j=1}^n(\vr_j/\sqrt{t}-(\kappa-(n-j+1)\sqrt{t}/2)^2\big)}{t^{n^2+n/2}}
     \bigwedge\limits_{j=1}^n\dif \vr_j;
 \end{align*}
 in the sequel we denote the latter integral by $H_n(\kappa,t)$.
 Then setting $\xi_j=\vr_j/\sqrt{t}-(\kappa-(n-j+1))\sqrt{t}/2$, we have 
 \begin{align*}
 \vr_j=\xi_j\sqrt{t}+(\kappa-(n-j+1))t/2,
 \end{align*}
 whence one obtains
 \begin{align*}
     \bigwedge\limits_{j=1}^n\dif \vr_j=t^{n/2}\bigwedge\limits_{j=1}^n\dif \xi_j.
 \end{align*}
 The quantity $\varepsilon(\vr)^2$ now becomes
 \begin{align*}
    \varepsilon(\vr)^2=&\prod_{1\leq j\leq n} \vr_j^2 \prod_{1\leq l<m\leq n} (\vr_l-\vr_m)^2 
    \prod_{1\leq l<m\leq n} (\vr_l+\vr_m)^2 \\ 
    =&\;t^{n^2} \prod_{1\leq j\leq n} \bigg(\xi_j+\frac{(\kappa-(n-j+1))\sqrt{t}}{2}\bigg)^2 \prod_{1\leq l<m\leq n}\bigg( (\xi_l-\xi_m)+\frac{l-m}{2}\sqrt{t}\bigg)^2 \times \\ 
    &\quad\times \prod_{1\leq l<m\leq n}\bigg( (\xi_l+\xi_m)+\Big(\kappa-\big(n-\frac{l+m}{2}+1\big)\Big)\sqrt{t}\bigg)^2,
\end{align*}
which is a polynomial in $\xi=(\xi_1,\ldots,\xi_n),\,\kappa$ and $t$. 
Then, putting $\ve(\vr)^2$ back in \autoref{local integral3}, we have
\begin{align*}
    H_n(\kappa,t)=&\int\limits_{\xi_1=-\infty}^{\infty}\ldots \int\limits_{\xi_n=-\infty}^{\infty}\exp\big(-\sum_{j=1}^n\xi_j^2\big)  \prod_{1\leq l<m\leq n}\bigg( (\xi_l+\xi_m)+\Big(\kappa-\big(n-\frac{l+m}{2}+1\big)\Big)\sqrt{t}\bigg)^2 \times \\
    &\;\times \prod_{1\leq j\leq n} \bigg(\xi_j+\frac{(\kappa-(n-j+1))\sqrt{t}}{2}\bigg)^2 \prod_{1\leq l<m\leq n}\bigg( (\xi_l-\xi_m)+\frac{l-m}{2}\sqrt{t}\bigg)^2 
    \bigwedge\limits_{j=1}^n\dif \xi_j.
\end{align*}
After evaluating integrals of the form
\begin{align*}
    \int\limits_{\xi=-\infty}^{\infty} \xi^m \exp(-\xi^2) \dif\xi=\frac{1+(-1)^m}{2}\Gamma\Big(\frac{m+1}{2}\Big)\quad\quad(m\in\mb{N}_{\geq 0}),
\end{align*}
$H_n(\kappa,t)$ becomes a polynomial in $\kappa$ and $t$. As $\kappa\geq n+1>0$ and $t>0$, for an upper bound, we need only consider the highest powers of $\kappa$ and $t$ in  $H_n(\kappa,t)$. From the above formula, it follows that
\begin{align*}
    H_n(\kappa,t)\leq c_n\, \kappa^{n(n+1)}\,t^{n(n+1)/2} \mu(\sqrt{t}),
\end{align*}
where $\mu$ is a polynomial of order $n(n-1)$. Thus, from equations (\ref{basic inequality}) and (\ref{local integral2}), one obtains
\begin{align*}
    S_{\kappa}^{\Gamma}(Z)\leq c_{n,\Gamma} \,\kappa^{n(n+1)}\,t^{n(n+1)/2} \mu(\sqrt{t})\quad\quad(Z\in\mb{H}_n).
\end{align*}
Now multiplying both sides of the above inequality by $\exp(-\kappa t)$ and integrating over $t\in[0,\infty]$, we have
\begin{align*}
    \int\limits_{t=0}^{\infty}\exp(-\kappa t)S_{\kappa}^{\Gamma}(Z)\dif t=\frac{S_{\kappa}^{\Gamma}(Z)}{\kappa}\leq c_{n,\Gamma}\,\kappa^{n(n+1)} \int\limits_{t=0}^{\infty} \exp(-\kappa t)\,t^{n(n+1)/2} \mu(\sqrt{t})\dif t\leq   c_{n,\Gamma}\,\frac{\kappa^{n(n+1)}}{\kappa^{n(n+1)/2+1}},
\end{align*}
whence it easily follows that
\begin{align*}
    S_{\kappa}^{\Gamma}(Z)\leq c_{n,\Gamma}\, \kappa^{n(n+1)/2}\quad\quad(Z\in\mb{H}_n)
\end{align*}
thereby proving the result stated in the theorem. 
\end{proof}


\subsection{Sup-norm bounds in the cofinite setting}

\begin{theorem}\label{cofinite 1}
For any arithmetic subgroup $\Gamma\subsetneq \mathrm{Sp}_n(\mb{R})$ such that $M:=\Gamma\backslash\mb{H}_n$ is of finite volume, we have
\begin{align*}
   \sup_{Z\in\mb{H}_n} S_{\kappa}^{\Gamma}(Z)\leq c_{n}\, \kappa^{n(n+1)/2}\sum_{\gamma\in\Gamma}\,\frac{1}{\prod_{j=1}^n \ch^{\kappa}(r_j^{\gamma}(Z))}\quad\quad( \kappa\geq n+1),
\end{align*}
where $r_j^{\gamma}(Z)$ denotes the diagonal entries of the diagonal matrix $R^{\gamma}(Z)=R(Z,\gamma Z)$ and $c_{n}$ is a positive real constant depending only on $n$ and $\Gamma$.
\end{theorem}
\begin{proof}
From \autoref{limit connection}, we have
\begin{align}\label{local Sk-inequality}
    S_{\kappa}^{\Gamma}(Z)= \lim_{t\rightarrow\infty}\exp\big(-\frac{n\kappa}{4}(\kappa-(n+1))\,t\big)\sum_{\gamma\in\Gamma}K_t^{(\kappa)}(2R^{\gamma}(Z))
\end{align}
and from \autoref{heat corollary}, we have
\begin{align*}
   K_t^{(\kappa)}(2R^{\gamma}(Z))\leq \frac{I_n(\kappa,t,R^{\gamma}(Z))}{\prod_{j=1}^n \ch^{\kappa}(r^{\gamma}_j(Z))},
\end{align*}
where $I_n(\kappa,t,R^{\gamma}(Z))$ is the integral given by
\begin{align*}
    I_n(\kappa,t,R^{\gamma}(Z))=c_n\frac{\exp{\big(-\sum_{j=1}^n j^2t/4\big)}}{t^{n^2+n/2}} \int\limits_{k\in K} \frac{\ve(\vr)\exp\big(-\sum_{j=1}^n\,\big(\vr_j^2/t-\kappa\vert\vr_j\vert\big)\big)}{\delta(\vr)} \dif\mu(k).
\end{align*}
Here $\vr=\vr(r^{\gamma}(Z),k)$ is the diagonal matrix $\vr(r^{\gamma}(Z),k)=\big(\begin{smallmatrix}P(r^{\gamma}(Z),k) & 0\\0 & -P(r^{\gamma}(Z),k)\end{smallmatrix}\big)$ with 
\begin{align*}
P(r^{\gamma}(Z),k)= \begin{pmatrix}
{\vr_{1}(r^{\gamma}(Z),k)} &{} & {0} \\
{} & {\ddots}& \\
{0} &{}& {\vr_{n}(r^{\gamma}(Z),k)}
\end{pmatrix}\quad(\vr_j(r^{\gamma}(Z),k)\in\mb{R},\;1\leq j\leq n)
\end{align*}
related to 
\begin{align*}
r^{\gamma}(Z)=\begin{pmatrix}R^{\gamma}(Z) & 0\\0 & -R^{\gamma}(Z)\end{pmatrix}
\end{align*}
via the matrix equality 
\begin{align}\label{local matrix equality}
k\exp(r^{\gamma}(Z))\overline{k}^t=u\exp(\vr)\overline{u}^t\quad\quad(k\in K,\,u\in U).
\end{align}
Since heat kernels decrease rapidly with increasing distance, the integral $I_n(\kappa,t,R^{\gamma}(Z))$ also decreases rapidly with increasing distance and hence we have 
\begin{align}\label{local zero reduction}
    I_n(\kappa,t,R^{\gamma}(Z))\leq I_n(\kappa,t,0_n). 
\end{align}
Then, for $R^{\gamma}(Z)=0_n$, the matrix equality (\ref{local matrix equality}) becomes
\begin{align}\label{local eigendecomposition}
    k\overline{k}^t=u\exp(\vr)\overline{u}^t\quad\quad(k\in K,\,u\in U).
\end{align}

For the eigendecomposition \eqref{local eigendecomposition} of $k\overline{k}^t$, we next determine $u$ and $\exp(\vr)$ in terms of $k\in K$. For this, we use the matrix $l$ from \eqref{ldef} and calculate
\begin{align}\label{local l-decomposition}
    k\overline{k}^t=\begin{pmatrix} A & B\\ -B & A\end{pmatrix}=
    l^{-1}\begin{pmatrix} A+iB & 0\\ 0 & A-iB\end{pmatrix}l,
\end{align}
where the matrix $h:=A+iB$ is Hermitian, as $k\overline{k}^t$ is Hermitian; note that $A-iB=h^{-t}$. Since $h$ is Hermitian, we have
\begin{align*}
h=vD\overline{v}^t,
\end{align*}
where $v\in U_n$ and $D$ is a real diagonal $(n\times n)$-matrix. Substituting this into \eqref{local l-decomposition} yields
\begin{align*}
 k\overline{k}^t=l^{-1}\begin{pmatrix} v & 0\\ 0 & \overline{v}\end{pmatrix} \begin{pmatrix} D & 0\\ 0 & D^{-1}\end{pmatrix} \begin{pmatrix} \overline{v}^t & 0\\ 0 & v^t \end{pmatrix} l.
\end{align*}
In this way, the factors on the right-hand side of \eqref{local eigendecomposition} become
\begin{align}\label{local u-structure}
    u=l^{-1}\begin{pmatrix} v & 0\\ 0 & \overline{v}\end{pmatrix}\quad \text{and} \quad \exp(\vr)= \begin{pmatrix} D & 0\\ 0 & D^{-1}\end{pmatrix}. 
\end{align}
Note that the eigendecomposition \eqref{local eigendecomposition} is unique only up to the ordering of the eigenvalues $\exp(\pm\vr_j)$ ($1\leq j\leq n$), i.e., we can always choose $u\in U$ in such a way so that $\vr_j\in\mb{R}_{\geq 0}$ ($1\leq j\leq n$). Therefore, without loss of generality, for the rest of the calculation, we assume $\vr_j\in\mb{R}_{\geq 0}$ ($1\leq j\leq n$).
 
Next, we determine the invariant volume form $\dif\mu(k)$ in terms of $\vr$ and $v$ by proceeding as in the proof of \autoref{coordinate change lemma 1}. From $x=k\overline{k}^t=u\exp(\vr)\overline{u}^t$, one obtains
\begin{align*}
    \dif x=\dif k\;\overline{k}^t+k\dif\overline{k}^t=\dif u\exp(\vr)\overline{u}^t+u\exp(\vr)\dif\vr\;\overline{u}^t+u\exp(\vr)\dif\overline{u}^t.
\end{align*}
Now as $x^{-1}=(k\overline{k}^t)^{-1}=\overline{k}k^t=u\exp(-\vr)\overline{u}^t$, we have
\begin{align*}
    x^{-1}\dif x=\overline{k}\;(k^t\dif k+\dif\overline{k}^t\overline{k})\;\overline{k}^t
    =ue^{-\vr/2}\;(e^{-\vr/2}(\overline{u}^t \dif u)e^{\vr/2}+e^{-\vr/2}(\dif\vr)e^{\vr/2}+e^{\vr/2}(\dif\overline{u}^t u)e^{-\vr/2})\;e^{\vr/2}\overline{u}^t. 
\end{align*}
Noting that $\dif\overline{u}^t u=-\overline{u}^t \dif u$, we take the volume form on both sides, denoted by the square brackets $[\,\cdot\,]$, to obtain
\begin{align*}
    [k^t\dif k+\overline{(k^t\dif k)}^t]=[e^{-\vr/2}(\overline{u}^t \dif u)e^{\vr/2}+e^{-\vr/2}(\dif\vr)e^{\vr/2}-e^{\vr/2}(\overline{u}^t \dif u)e^{-\vr/2}].
\end{align*}
From the structure of $u$ obtained in \eqref{local u-structure}, it is easy to see that the invariant matrix differential form $\overline{u}^t \dif u$ is of the form
\begin{align*}
    \overline{u}^t \dif u=\begin{pmatrix} \overline{v}^t \dif v & 0\\ 0 & v^t\dif\overline{v}\end{pmatrix}. 
\end{align*}
Now, writing $e^{-\vr/2}(\overline{u}^t \dif u)e^{\vr/2}+e^{-\vr/2}(\dif\vr)e^{\vr/2}-e^{\vr/2}(\overline{u}^t \dif u)e^{-\vr/2}$ in the familiar block decomposed form, we have
\begin{align*}
    [k^t\dif k+\overline{(k^t\dif k)}^t]=&\Bigg[\begin{pmatrix}e^{-P/2}&0\\0&e^{P/2}\end{pmatrix}\begin{pmatrix} \overline{v}^t \dif v & 0\\ 0 & v^t\dif\overline{v}\end{pmatrix}\begin{pmatrix}e^{P/2}&0\\0&e^{-P/2}\end{pmatrix}\\
    &\quad-\begin{pmatrix}e^{P/2}&0\\0&e^{-P/2}\end{pmatrix}
    \begin{pmatrix} \overline{v}^t \dif v & 0\\ 0 & v^t\dif\overline{v}\end{pmatrix}\begin{pmatrix}e^{-P/2}&0\\0&e^{P/2}\end{pmatrix}
    +\begin{pmatrix}\dif P&0\\0&-\dif P\end{pmatrix}\Bigg]\\
\end{align*}
The right-hand side of the above equation gives
\begin{align*}
    \sbr{\begin{pmatrix}\dif P+ e^{-P/2}(\overline{v}^t \dif v)e^{P/2}-e^{P/2}(\overline{v}^t \dif v)e^{-P/2} & 0\\
    0 & -\dif P+ e^{P/2}(v^t\dif\overline{v})e^{-P/2}-e^{-P/2}(v^t\dif\overline{v})e^{P/2}
    \end{pmatrix}}
\end{align*}
Now, taking $\overline{v}^t \dif v=(\omega_{j,k})_{1\leq j,k\leq n}$, we have 
\begin{align*}
    (\dif P+e^{-P/2}(\overline{v}^t \dif v)e^{P/2}-e^{P/2}(\overline{v}^t \dif v)e^{-P/2})_{j,k}=\delta_{j,k}\,\dif\vr_j+2\sh\big(\frac{\vr_j-\vr_k}{2}\big)\,\omega_{j,k}.
\end{align*}
Therefore, we have
\begin{align*}
    [k^t\dif k+\overline{(k^t\dif k)}^t]=c_n\prod\limits_{1\leq j< k\leq n}\sh^2\big(\frac{\vr_j-\vr_k}{2}\big)\bigwedge\limits_{j=1}^n\dif \vr_j\bigwedge\limits_{1\leq j<k\leq n}(\omega_{j,k}\wedge\overline{\omega}_{j,k}).
\end{align*}
Since $k^t\dif k\in\mf{k}$ is of the form
\begin{align*}
    k^t\dif k=\begin{pmatrix} A & B \\ -B & A\end{pmatrix}\quad\quad
    ( A,B\in\mathbb{C}^{n\times n}, B=B^t, A=-A^t),
\end{align*}
we have
\begin{align*}
    k^t\dif k+\overline{(k^t\dif k)}^t=2i\begin{pmatrix} \im(A) & \im(B) \\ -\im(B) & \im(A)\end{pmatrix}= 2i\im(k^t\dif k). 
\end{align*}
Then,identifying $\dif\mu(v)=[\overline{v}^t \dif v]=\wedge_{1\leq j<k\leq n}(\omega_{j,k}\wedge\overline{\omega}_{j,k})$, we have
\begin{align*}
    \dif\mu(k)=c_n\prod\limits_{1\leq j< k\leq n}\sh^2\big(\frac{\vr_j-\vr_k}{2}\big)\bigwedge\limits_{j=1}^n\dif \vr_j
    \bigwedge \dif\mu(v) \bigwedge \dif\mu(k_0)\quad\quad(k_0\in K_0). 
\end{align*}
This allows us to write
\begin{align*}
    I_n(\kappa,t,0_n)=& c_n\frac{\exp{\big(-\sum_{j=1}^n j^2t/4\big)}}{t^{n^2+n/2}} \int\limits_{\vr_1=0}^{\infty}\ldots\int\limits_{\vr_n=0}^{\infty}\frac{\ve(\vr)\exp\big(-\sum_{j=1}^n\,(\vr_j^2/t-\kappa\vr_j)\big)}{\delta(\vr)}\times \\
    &\quad \quad \times \prod\limits_{1\leq l< m\leq n}\sh^2\big(\frac{\vr_l-\vr_m}{2}\big)\bigwedge\limits_{j=1}^n\dif \vr_j.
\end{align*}
Therefore, we have
\begin{align*}
    \exp\big(-\frac{n\kappa}{4}(\kappa-(n+1))\,t\big)I_n(\kappa,t,0_n)
     = & c_n\int\limits_{\vr_1=0}^{\infty}\ldots\int\limits_{\vr_n=0}^{\infty}\;
    \frac{\ve(\vr) \prod\limits_{j=1}^n \exp\big(-(\vr_j/\sqrt{t}-(\kappa-(n-j+1))\sqrt{t}/2)^2\big)}{t^{n^2+n/2}} \times \\
    & \quad \quad \times \prod\limits_{j=1}^n\frac{\exp(\vr_j)}{\sh(\vr_j)}
    \prod\limits_{1\leq l< m\leq n}\frac{\exp(\vr_l)\sh\big((\vr_l-\vr_m)/2\big)}{\sh\big((\vr_l+\vr_m)/2\big)}
     \bigwedge\limits_{j=1}^n\dif \vr_j
\end{align*}
Now setting $\xi_j=\vr_j/\sqrt{t}-(\kappa-(n-j+1))\sqrt{t}/2$, we have 
 \begin{align}\label{local substitution}
 \vr_j=\xi_j\sqrt{t}+(\kappa-(n-j+1))t/2,
 \end{align}
 whence one obtains
 \begin{align*}
     \bigwedge\limits_{j=1}^n\dif \vr_j=t^{n/2}\bigwedge\limits_{j=1}^n\dif \xi_j.
 \end{align*}
 Now see that
 \begin{align*}
     \lim_{t\rightarrow\infty}\frac{\ve(\vr)}{t^{n^2}}
     =&\lim_{t\rightarrow\infty} \; \prod_{1\leq j\leq n} \frac{\vr_j}{t} \prod_{1\leq l<m\leq n} \Big(\frac{\vr_l}{t}+\frac{\vr_m}{t}\Big) 
     \prod_{1\leq l<m\leq n} \Big(\frac{\vr_l}{t}-\frac{\vr_m}{t}\Big)\\
     =& \prod_{1\leq j\leq n} \frac{\kappa-(n-j+1)}{2} \prod_{1\leq l<m\leq n}
    \Big(\kappa-\big(n-\frac{l+m}{2}+1\big)\Big) \prod_{1\leq l<m\leq n} \frac{l-m}{2}.
 \end{align*}
 Also, as $t\rightarrow \infty$, by the substitution (\ref{local substitution}), we have $\vr_j\rightarrow\infty\,(1\leq j\leq n)$. Therefore, taking the limit at $t\rightarrow\infty$, we obtain
 \begin{align*}
     \lim_{t\rightarrow\infty}\frac{\exp(\vr_j)}{\sh(\vr_j)}
     =\lim_{\vr_j\rightarrow\infty}\frac{\exp(\vr_j)}{\sh(\vr_j)}=2.
 \end{align*}
 Next, as 
 \begin{align*}
     \vr_l-\vr_m=(\xi_l-\xi_m)\sqrt{t}+(l-m)t
 \end{align*}
 and for $l<m$ we have $l-m<0$, the quantity $\exp(\vr_l-\vr_m)\rightarrow 0$ as $t\rightarrow\infty$.  Therefore, we have
 \begin{align*}
 \lim_{t\rightarrow\infty}\frac{\exp(\vr_l)\sh\big((\vr_l-\vr_m)/2\big)}{\sh\big((\vr_l+\vr_m)/2\big)}
 &=\lim_{t\rightarrow\infty}\frac{\exp(\vr_l)\;
 \big(\exp\big((\vr_l-\vr_m)/2\big)-\exp\big(-(\vr_l-\vr_m)/2\big)\big)}{\exp\big((\vr_l+\vr_m)/2\big)-\exp\big(-(\vr_l+\vr_m)/2\big)}\\
 &=\lim_{t\rightarrow\infty}\frac{\exp(\vr_l-\vr_m)-1}{1-\exp(-(\vr_l+\vr_m))}
 =-1
 \quad \quad (1\leq l<m\leq n).
 \end{align*}
 Combining all the above limits, we have
 \begin{align*}
      &\lim_{t\rightarrow\infty}   \exp\big(-\frac{n\kappa}{4}(\kappa-(n+1))\,t\big)I_n(\kappa,t,0_n) \\
      &=c_n \prod_{1\leq j\leq n} \frac{\kappa-(n-j+1)}{2} \prod_{1\leq l<m\leq n}
    \Big(\kappa-\big(n-\frac{l+m}{2}+1\big)\Big) \prod_{1\leq l<m\leq n}  \frac{m-l}{2} \\
    &\leq c_n\; \kappa^{n(n+1)/2}. 
 \end{align*}
 Then, from \autoref{local Sk-inequality} and \autoref{local zero reduction}, the statement of the theorem easily follows. 
\end{proof}
\begin{theorem}\label{cofinite 2}
For any arithmetic subgroup $\Gamma\subsetneq \mathrm{Sp}_n(\mb{R})$ such that $M:=\Gamma\backslash\mb{H}_n$ is of finite volume, we have
\begin{align*}
   \sup_{Z\in\mb{H}_n} S_{\kappa}^{\Gamma}(Z)\leq c_{n,\Gamma}\, \kappa^{3n(n+1)/4}\quad\quad( \kappa\geq n+1),
\end{align*}
where $c_{n,\Gamma}$ is a positive real constant depending only on $n$ and $\Gamma$.
\end{theorem}
\begin{proof}
By \autoref{arithmetic boundary}, we know that the boundary $M^{\star}\setminus M$ of $M$ consists of finite union of subspaces $M_j:=(\Gamma\cap P(\mb{P}_j))\backslash \mb{P}_j$, where $\mb{P}_j$ runs through a set of representatives of equivalence classes modulo $\Gamma$ of rational boundary components of $\mb{H}_n$, and its subspaces of strictly smaller degree. We denote by $\mathcal{C}$ the set of all such inequivalent chains of boundary components of $M$. Then, for $\mathcal{P}\in\mathcal{C}$, we can define boundary neighbourhoods $U_{\ve}(\mathcal{P})$ containing the entire chain $\mathcal{P}$, such that the complement of their union in $M$, i.e.,
\begin{align*}
    K_{\ve}:=M\setminus \bigcup_{P\in\mathcal{C}} U_{\ve}(\mathcal{P})
\end{align*}
is a compact subset of $M$. We shall now estimate $S_{\kappa}^{\Gamma}(Z)$ for $Z$ ranging through $K_{\ve}$ and $U_{\ve}(\mathcal{P})\,(\mathcal{P}\in\mathcal{C})$, respectively.

In case of the compact set $K_{\ve}$, using \autoref{compact case}, we have already determined that
\begin{align*}
   \sup_{Z\in K_{\ve}} S_{\kappa}^{\Gamma}(Z)\leq c_{n,\Gamma}\, \kappa^{n(n+1)/2}\quad\quad( \kappa\geq n+1),
\end{align*}
where the constant $c_{n,\Gamma}>0$ depends only on $n$ and $\Gamma$.. 

Next, in case of $U_{\ve}(\mathcal{P})\,(\mathcal{P}\in\mathcal{C})$,
by \autoref{boundary equivalence}, without loss of generality, we can assume $\mathcal{P}$ to be the chain 
\begin{align*}
\Gamma_0\backslash\mb{H}_0< \Gamma_1\backslash\mb{H}_1<\ldots<\Gamma_j\backslash\mb{H}_j<\ldots<\Gamma_{n-1}\backslash\mb{H}_{n-1}
\end{align*}
of standard boundary components of $\Gamma_n\backslash\mb{H}_n$. 

Let $\mathscr{F}_n$ denote the standard fundamental domain of the Siegel modular group $\Gamma_n$. For $Z\in\mathscr{F}_n$, there exists a constant $ c_3(n)>0$ depending only on $n$, such that $Y\geq  c_3(n)\mathbbm{1}_n$ (see \autoref{SEC arithmetic subgroups}). Let $\lambda_j(Y)\,(1\leq j\leq n)$ denote the ordered set of eigenvalues
\begin{align*}
     c_3(n)\leq\lambda_1(Y)\leq \ldots\leq \lambda_j(Y)\leq\ldots\leq\lambda_n(Y)
\end{align*}
of the positive definite matrix $Y$. Then $U_{\ve}(\mathcal{P})$ can be taken as the neighbourhood
\begin{align*}
    S_{\ve}:=\{Z=X+iY\in\mathscr{F}_n\;\vert\;\lambda_n(Y)>\ve\}
\end{align*}
of the standard boundary components of $\mathscr{F}_n$. As $\lambda_n(Y)$ denotes the highest eigenvalue of $Y$, the complement of $S_{\ve}$ in $\mathscr{F}_n$ is then given by the compact subset 
\begin{align*}
   K_{\ve}=\{Z=X+iY\in\mathscr{F}_n\;\vert\; c_3(n)\mathbbm{1}_n\leq Y\leq \ve\mathbbm{1}_n\}.
\end{align*}

Let $f\in \mathcal{S}_{\kappa}^n(\Gamma)$ be a cusp form of weight $\kappa$. For $Z\in\mb{H}_n$ such that $Y$ is Minkowski reduced and $Y>c\mathbbm{1}_n$ for some $c>0$, there exist positive numbers $c_1(n,c)>0$ and $c_2(n,c)>0$ depending only on $n$ and $c$, such that
\begin{align*}
    \vert f(Z)\vert \leq c_1(n,c)\, \exp(-c_2(n,c)\tr(Y)). 
\end{align*}
(see \cite[page 57]{Klingen}). Since here we consider $Z\in\mathscr{F}_n$, we can take $c= c_3(n)$. In that case, the positive numbers $c_1(n,c)>0$ and $c_2(n,c)>0$ depend only on $n$ and we have
\begin{align}\label{klingen bound}
    \vert f(Z)\vert \leq c_1(n)\, \exp(-c_2(n)\tr(Y))
    \quad\quad(Z\in\mathscr{F}_n). 
\end{align}
This shows that the function $f(Z)/\exp(ic_2(n)\tr(Z))$ is a bounded holomorphic function on $S_{\ve}$ and hence, by maximum modulus principle, its absolute value
\begin{align*}
\bigg\vert\frac{ f(Z)}{\exp\big(ic_2(n)\tr(Z)\big)}\bigg\vert^2= exp(2c_2(n)\tr(Y))\vert f(Z)\vert^2
\end{align*}
takes its maximum value at the boundary 
\begin{align*}
   \partial S_{\ve}=\{Z=X+iY\in\mathscr{F}_n\;\vert\;\lambda_n(Y)=\ve\}
\end{align*}
of $S_{\ve}$. Now, write $\det(Y)^{\kappa} \vert f(Z)\vert^2$ as
\begin{align*}
    \det(Y)^{\kappa} \vert f(Z)\vert^2=\exp(2c_2(n)\tr(Y))\vert f(Z)\vert^2\,\frac{\det(Y)^{\kappa}}{\exp(2c_2(n)\tr(Y))}. 
\end{align*}
Then writing the eigenvalues of $Y$ as $\lambda_j(Y)\,(1\leq j\leq n)$, we have
\begin{align*}
    \frac{\det(Y)^{\kappa}}{\exp(2c_2(n)\tr(Y))}=\prod_{j=1}^n \frac{\lambda_j(Y)^{\kappa}}{\exp(2c_2(n)\lambda_j(Y))}. 
\end{align*}
The functions ${\lambda_j^{\kappa}}/{\exp(2c_2(n)\lambda_j})$ attain maxima at $\lambda_j=\kappa/(2c_2(n))$ and decreases monotonically for  $\lambda_j>\kappa/(2c_2(n))$. Therefore, if we choose $\ve>\kappa/(2c_2(n))$, then we have
\begin{align*}
    \sup_{Z\in M} S_{\kappa}^{\Gamma}(Z)=\sup_{Z\in K_{\ve}} S_{\kappa}^{\Gamma}(Z)\leq c_{n,\Gamma}\, \kappa^{n(n+1)/2}\quad\quad\bigg( \kappa\geq n+1,\,\ve>\frac{\kappa}{2c_2(n)}\bigg). 
\end{align*}
Now, in case $\ve\leq \kappa/(2c_2(n))$, we need to determine $\sup_{Z\in M} S_{\kappa}^{\Gamma}(Z)$ in the annulus
\begin{align*}
    S_{\ve}\setminus S_{\kappa/(2c_2(n))}&= \bigg\{Z=X+iY\in\mathscr{F}_n\;\bigg\vert\;\ve<\lambda_n(Y)\leq\frac{\kappa}{2c_2(n)}\bigg\}\\
    &\subsetneq  \bigg\{Z=X+iY\in\mathscr{F}_n\;\bigg\vert\;Y\leq\frac{\kappa}{2c_2(n)}\mathbbm{1}_n\bigg\}.
\end{align*}
We do this using Theorems \ref{compact case} and \ref{cofinite 1}.

From \autoref{local Sk-inequality}, we have
\begin{align*}
   S_{\kappa}^{\Gamma}(Z)= \lim_{t\rightarrow\infty}\exp\big(-\frac{n\kappa}{4}(\kappa-(n+1))\,t\big)\sum_{\gamma\in\Gamma}K_t^{(\kappa)}(2R^{\gamma}(Z)).
\end{align*}
We split the sum over $\Gamma$ according as whether there is a minimum distance between the point $\gamma Z$ and $Z$ or they can get arbitrarily close. Let $\Gamma_{\infty}$ denote the set of elements of $\Gamma$ for which $\gamma Z$ and $Z$ can get arbitrarily close. Then we split the above sum as
\begin{align*}
    S_{\kappa}^{\Gamma}(Z)=&\lim_{t\rightarrow\infty} \exp\big(-\frac{n\kappa}{4}(\kappa-(n+1))\,t\big)\sum_{\gamma\in\Gamma\setminus\Gamma_{\infty}}K_t^{(\kappa)}(2R^{\gamma}(Z))\\
    &\quad +\lim_{t\rightarrow\infty} \exp\big(-\frac{n\kappa}{4}(\kappa-(n+1))\,t\big)\sum_{\gamma\in\Gamma_{\infty}}K_t^{(\kappa)}(2R^{\gamma}(Z)).
\end{align*}
As the function $\exp\big(-\frac{n\kappa}{4}(\kappa-(n+1))\,t\big)\sum_{\gamma\in\Gamma\setminus\Gamma_{\infty}}K_t^{(\kappa)}(2R^{\gamma}(Z))$ is monotonically decreasing in $t$, we have
\begin{align}\label{infinity split}
\begin{aligned}
    S_{\kappa}^{\Gamma}(Z)\leq& \exp\big(-\frac{n\kappa}{4}(\kappa-(n+1))\,t\big)\sum_{\gamma\in\Gamma\setminus\Gamma_{\infty}}K_t^{(\kappa)}(2R^{\gamma}(Z)).\\
    &\quad +\lim_{t\rightarrow\infty} \exp\big(-\frac{n\kappa}{4}(\kappa-(n+1))\,t\big)\sum_{\gamma\in\Gamma_{\infty}}K_t^{(\kappa)}(2R^{\gamma}(Z)).
    \end{aligned}
\end{align}
As for $\gamma\in\Gamma\setminus\Gamma_{\infty}$ the points $\gamma Z$ and $Z$ cannot be arbitrarily close, the first sum can be handled exactly as in \autoref{compact case} using the counting function to estimate the sum by an integral to give
\begin{align}\label{gamma minus gamma infinity sum}
 \exp\big(-\frac{n\kappa}{4}(\kappa-(n+1))\,t\big)\sum_{\gamma\in\Gamma\setminus\Gamma_{\infty}}K_t^{(\kappa)}(2R^{\gamma}(Z))\leq   c_{n,\Gamma}\, \kappa^{n(n+1)/2}. 
\end{align}
The second sum was estimated in \autoref{cofinite 1} to be
\begin{align}\label{local pretheorem}
\lim_{t\rightarrow\infty} \exp\big(-\frac{n\kappa}{4}(\kappa-(n+1))\,t\big)\sum_{\gamma\in\Gamma_{\infty}}K_t^{(\kappa)}(2R^{\gamma}(Z))
\leq 
c_{n}\, \kappa^{n(n+1)/2}\sum_{\gamma\in\Gamma_{\infty}}\,\frac{1}{\prod_{j=1}^n \ch^{\kappa}(r_j^{\gamma}(Z))}.
\end{align}
Thus, it only remains to estimate the sum
\begin{align}\label{gamma-sum}
    \sum_{\gamma\in\Gamma_{\infty}}\,\frac{1}{\prod_{j=1}^n \ch^{\kappa}(r_j^{\gamma}(Z))}. 
\end{align}

Since $\Gamma_{\infty}$ is defined as the set of elements of $\Gamma$ for which $\gamma Z$ and $Z$ can get arbitrarily close, by \autoref{arbitrary closeness}, we have
\begin{align*}
    \Gamma_{\infty}=\bigcup_{j=0}^{n-1}\Gamma_{\infty}^j, 
\end{align*}
where $\Gamma_{\infty}^j:=\Gamma\cap W_j$. Thus, by \eqref{W_j definition}, these groups are explicitly given by
\begin{align*}
    \Gamma_{\infty}^0&=\bigg\{\begin{pmatrix}\mathbbm{1}_n & S\\ 0 & \mathbbm{1}_n\end{pmatrix}\;\bigg\vert\; S\in\mr{Sym}_n(\mb{Z})\bigg\},\\
    \Gamma_{\infty}^j&=\bigg\{\begin{pmatrix}A & A S\\ 0 & A^{-t}\end{pmatrix}\;\bigg\vert\; 
    A=\begin{pmatrix}\mathbbm{1}_j & 0\\ L & \mathbbm{1}_{n-j}\end{pmatrix},\,
    S=\begin{pmatrix}0 & H^t\\ H & S_2\end{pmatrix}\bigg\}\quad\quad(1\leq j\leq (n-1)), 
\end{align*}
where $\,L,H\in\mb{Z}^{(n-j)\times j}$ and $S_2\in\mr{Sym}_{n-j}(\mb{Z})$.

Next, we need an effective way of calculating the quantity 
$1/\prod_{j=1}^n \ch^{\kappa}(r_j^{\gamma}(Z))$. Here we derive a more general formula for the quantity $1/{\prod_{j=1}^n \ch^2(r_j(Z,W))}\,(Z,W\in\mb{H}_n)$, where setting $W=\gamma Z\,(\gamma\in\Gamma_{\infty})$, we can easily obtain the sum in (\ref{gamma-sum}) above.

Recall from \autoref{SEC Siegel geometry} the cross ratio
\begin{align*}
    \rho(W, Z)=(W-Z)(\overline{W}-Z)^{-1}(\overline{W}-\overline{Z})(W-\overline{Z})^{-1}\quad(Z,W\in\mb{H}_n).
\end{align*}
Let $\rho_j(Z,W)\,(1\leq j\leq n)$ denote the eigenvalues of $\rho(W, Z)$. The point $Z=X+iY\in\mb{H}_n$, where $X,Y\in\mb{R}^{n\times n}$ with $Y>0$ can be written as
\begin{align*}
    Z=\begin{pmatrix}\mathbbm{1}_n & X\\ 0 & \mathbbm{1}_n\end{pmatrix}
    \begin{pmatrix}Y^{1/2} & 0\\ 0 & Y^{-1/2}\end{pmatrix}\cdot i\mathbbm{1}_n.
\end{align*}
Now, as the matrices $\rho(Z,W)$ and $\rho(gZ,gW)$ have the same set of eigenvalues for all $g\in \mathrm{Sp}_n(\mb{R})$, setting
\begin{align}\label{V-definition}
    V=\begin{pmatrix}Y^{-1/2} & 0\\ 0 & Y^{1/2}\end{pmatrix}\begin{pmatrix}\mathbbm{1}_n & -X\\ 0 & \mathbbm{1}_n\end{pmatrix}\cdot W=Y^{-1/2}(W-X)Y^{-1/2},
\end{align}
the cross ratio $\rho(V, i\mathbbm{1}_n)$ has the same eigenvalues as $\rho(Z,W)$, i.e., $\rho_j(Z,W)\,(1\leq j\leq n)$. Therefore, we have
\begin{align*}
    \det(\mathbbm{1}_n-\rho(Z,W))=\det\big(\mathbbm{1}_n-(V-i\mathbbm{1}_n)(V+i\mathbbm{1}_n)^{-1}(\overline{V}+i\mathbbm{1}_n)
    (\overline{V}-i\mathbbm{1}_n)^{-1}\big). 
\end{align*}
Since these eigenvalues are of the form 
\begin{align*}
    \rho_{j}(Z,W)=\Th^2(r_j(Z,W))\quad\quad(1\leq j\leq n),
\end{align*}
from the above equations, using the fact $(V-i\mathbbm{1}_n)(V+i\mathbbm{1}_n)^{-1}=(V+i\mathbbm{1}_n)^{-1}(V-i\mathbbm{1}_n)$, one obtains
\begin{align*}
    \frac{1}{\prod_{j=1}^n \ch^2(r_j(Z,W))}&=
    \frac{\det\big((V+i\mathbbm{1}_n)(\overline{V}-i\mathbbm{1}_n)
    -(V-i\mathbbm{1}_n)(\overline{V}+i\mathbbm{1}_n)\big)}{\det(V+i\mathbbm{1}_n)\det(\overline{V}-i\mathbbm{1}_n)}\\
    &=\frac{\det(2i(\overline{V}-V))}{\det(V+i\mathbbm{1}_n)\det(\overline{V}-i\mathbbm{1}_n)}.
\end{align*}
Then, using the definition of $V$ in \autoref{V-definition}, one obtains
\begin{align}\label{cosh product formula}
    \frac{1}{\prod_{j=1}^n \ch^2(r_j(Z,W))}=\frac{4^n \det(\im(Z))\det(\im(W))}{\vert\det(W-\overline{Z})\vert^2}. 
\end{align}

Next we estimate the sum in (\ref{gamma-sum}) by breaking the sum over $\Gamma_{\infty}$ into sums over $\Gamma^j_{\infty}\,(0\leq j\leq (n-1))$. 

We begin with $\Gamma^0_{\infty}$. For $\gamma\in \Gamma^0_{\infty}$, i.e.,
\begin{align*}
    \gamma=\begin{pmatrix}\mathbbm{1}_n & S\\ 0 & \mathbbm{1}_n\end{pmatrix}\quad\quad (S\in\mr{Sym}_n(\mb{Z}))
\end{align*}
we have $\gamma Z=Z+S$. Therefore, putting $W=Z+S$ in \autoref{cosh product formula}, we obtain
\begin{align*}
    \frac{1}{\prod_{j=1}^n \ch^2(r_j^{\gamma}(Z))}
    =\frac{4^n \det(Y)^2}{\det(S-2iY)\det(S+2iY)}
    =\frac{1}{\det(\mathbbm{1}_n+(\frac{1}{2}Y^{-1/2}SY^{-1/2})^2)}.
\end{align*}
Then we estimate the sum over $\Gamma^0_{\infty}$ by the matrix beta integral
\begin{align}\label{beta integral connection}
    \sum_{\gamma\in\Gamma^0_{\infty}}\,\frac{1}{\prod_{j=1}^n \ch^{\kappa}(r_j^{\gamma}(Z))}\leq \int\limits_{S\in\mr{Sym}_n(\mb{R})}
    \frac{[\dif S]}{\det(\mathbbm{1}_n+(\frac{1}{2}Y^{-1/2}SY^{-1/2})^2)^{\kappa/2}}
\end{align}
Now, setting $T=\frac{1}{2}Y^{-1/2}SY^{-1/2}$, we have
\begin{align*}
    [\dif T]=c_n \det(Y)^{-(n+1)/2}[\dif S]. 
\end{align*}
This gives us
\begin{align*}
     \sum_{\gamma\in\Gamma^0_{\infty}}\,\frac{1}{\prod_{j=1}^n \ch^{\kappa}(r_j^{\gamma}(Z))}\leq c_n \det(Y)^{(n+1)/2}\int\limits_{T\in\mr{Sym}_n(\mb{R})}
     \frac{[\dif T]}{\det(\mathbbm{1}_n+T^2)^{\kappa/2}}. 
\end{align*}
Then, using Hua's matrix beta integral (see \cite[page 33]{Hua}) 
\begin{align}\label{hua integral}
\int\limits_{T\in\mr{Sym}_n(\mb{R})} \frac{[\dif T]}{\left(\operatorname{det}\left(I+T^{2}\right)\right)^{\alpha}}
= \pi^{n(n+1)/4}\, \frac{\Gamma\left(\alpha-n/2\right)}{\Gamma(\alpha)} \prod_{\nu=1}^{n-1} \frac{\Gamma\left(2 \alpha-(n+\nu)/2\right)}{\Gamma(2 \alpha-\nu)}\quad\quad(\alpha>n/2),
\end{align}
and $\det(Y)<(\kappa/(2c_2(n)))^n$, from the above calculations, it easily follows that
\begin{align}\label{gamma 0-sum}
    \sum_{\gamma\in\Gamma^0_{\infty}}\,\frac{1}{\prod_{j=1}^n \ch^{\kappa}(r_j^{\gamma}(Z))}\leq c_n \kappa^{n(n+1)/4}.
\end{align}

Next we consider the sum
\begin{align*}
    \sum_{\gamma\in\Gamma^j_{\infty}}\,\frac{1}{\prod_{j=1}^n \ch^{\kappa}(r_j^{\gamma}(Z))}\quad\quad(1\leq j\leq (n-1)). 
\end{align*}
For $\gamma\in \Gamma^j_{\infty}\,(1\leq j\leq (n-1))$, i.e.,
\begin{align*}
    \gamma=\begin{pmatrix}A & A S\\ 0 & A^{-t}\end{pmatrix}\quad\quad 
    \bigg(A=\begin{pmatrix}\mathbbm{1}_j & 0\\ L & \mathbbm{1}_{n-j}\end{pmatrix},\,
    S=\begin{pmatrix}0 & H^t\\ H & S_2\end{pmatrix}\bigg),
\end{align*}
where $\,L,H\in\mb{Z}^{(n-j)\times j}$ and $S_2\in\mb{Z}^{(n-j)\times(n-j)},\,S_2=S_2^t$, we have 
$\gamma Z=A(Z+S)A^t$. Therefore, putting $W=A(Z+S)A^t$ in \autoref{cosh product formula}, we obtain
\begin{align*}
    \frac{1}{\prod_{j=1}^n \ch^2(r_j^{\gamma}(Z))}
    &=\frac{4^n \det(Y)^2}{\vert\det(A(Z+S)A^t-\overline{Z})\vert^2}\\
    &=\frac{4^n \det(Y)^2}{\big\vert\det\big((A(X+S)A^t-X)+i(AYA^t+Y)\big)\big\vert^2}
\end{align*}
Now, just as in the $j=0$ case above, we estimate the sum over $\Gamma^j_{\infty}$ by a matrix integral $I_{n,\kappa}(Z)$, i.e.,
\begin{align*}
\sum_{\gamma\in\Gamma^j_{\infty}}\,\frac{1}{\prod_{j=1}^n \ch^{\kappa}(r_j^{\gamma}(Z))}\leq I^j_{n,\kappa}(Z),
\end{align*}
where the integral $I_{n,\kappa}(Z)$ is given by
\begin{align*}
I^j_{n,\kappa}(Z)=\int\limits_{L}\int\limits_{H}\int\limits_{S_2} 
\frac{2^{n\kappa} \det(Y)^{\kappa}\;\; [\dif S_2]\wedge[\dif H]\wedge[\dif L]}{\big\vert\det\big((A(X+S)A^t-X)+i(AYA^t+Y)\big)\big\vert^{\kappa}}.
\end{align*}
Next consider the the block decomposition  
\begin{align*}
    X= \begin{pmatrix}X_1 & X_{12}^t\\ X_{12} & X_2\end{pmatrix}\quad\quad(X_1\in\mb{R}^{j\times j},\,X_2\in\mb{R}^{(n-j)\times (n-j)},\,X_{12}\in\mb{R}^{(n-j)\times j}).
\end{align*}
of the matrix $X\in\mb{R}^{n\times n}$. Then, we have
\begin{align*}
    A(X+S)A^t-X= \begin{pmatrix}0 & H^t+X_1L^t\\ H+LX_1 & S_2+(LH^t+HL^t)+(LX_{12}^t+X_{12}L^t+LX_1L^t)\end{pmatrix}.
\end{align*}
Now, since
\begin{align*}
    &[\dif\,(S_2+(LH^t+HL^t)+(LX_{12}^t+X_{12}L^t+LX_1L^t))]\wedge [\dif\,(H+LX_1)]\wedge [\dif L]\\
    &\quad =[\dif S_2]\wedge[\dif H]\wedge[\dif L],
\end{align*}
we can simply replace the term $A(X+S)A^t-X$ in $I_{n,\kappa}(Z)$ with $S$, to write
\begin{align*}
    I^j_{n,\kappa}(Z)=\int\limits_{L}\int\limits_{H}\int\limits_{S_2} 
\frac{2^{n\kappa} \det(Y)^{\kappa}\;\; [\dif S_2]\wedge[\dif H]\wedge[\dif L]}{\vert\det((AYA^t+Y)+iS)\vert^{\kappa}}.
\end{align*}
Next we write the positive definite matrix $Y>0$ in the Cholesky decomposed form $Y=BB^t$, where 
\begin{align*}
    B= \begin{pmatrix}P_1 & 0\\ P & P_2\end{pmatrix}\quad\quad(P_1\in\mb{R}^{j\times j},\,P_2\in\mb{R}^{(n-j)\times (n-j)},\,P\in\mb{R}^{(n-j)\times j})
\end{align*}
with $P_1,\,P_2$ non-singular lower triangular. 
Then we have
\begin{align}\label{intermediate integral 1}
    I^j_{n,\kappa}(Z)=\int\limits_{L}\int\limits_{H}\int\limits_{S_2} 
\frac{\;\; [\dif S_2]\wedge[\dif H]\wedge[\dif L]}{\big\vert\det\big(1/2(\mathbbm{1}_n+(B^{-1}AB)(B^{-1}AB)^t+iB^{-1}SB^{-t})\big)\big\vert^{\kappa}}.
\end{align}
The matrices $B^{-1}AB$ and $B^{-1}SB^{-t}$, in block decomposed form, are given by
\begin{align*}
    B^{-1}AB&=\begin{pmatrix}\mathbbm{1}_j & 0\\ P_2^{-1}LP_1 & \mathbbm{1}_{n-j}\end{pmatrix},\\
    B^{-1}SB^{-t}&=\begin{pmatrix}0 & P_1^{-1}H^tP_2^{-t}
    \\ P_2^{-1}H P_1^{-t}  & P_2^{-1}S_2P_2^{-t}-P_2^{-1}(HP_1^{-t}P^t+PP_1^{-1}H^t)P_2^{-t}\end{pmatrix},
\end{align*}
respectively. We set
\begin{align}
    2\ti{L}&=P_2^{-1}LP_1, \label{tilde setting 1}\\
    2\ti{H}&=P_2^{-1}H P_1^{-t},\label{tilde setting 2}\\
    2\ti{S_2}&=P_2^{-1}S_2P_2^{-t}-P_2^{-1}(HP_1^{-t}P^t+PP_1^{-1}H^t)P_2^{-t} \label{tilde setting 3}
\end{align}
Then the matrix $1/2(\mathbbm{1}_n+(B^{-1}AB)(B^{-1}AB)^t+iB^{-1}SB^{-t})$ in the denominator of the integrand in \autoref{intermediate integral 1} is given by
\begin{align*}
    \frac{1}{2}(\mathbbm{1}_n+(B^{-1}AB)(B^{-1}AB)^t+iB^{-1}SB^{-t})=
    \begin{pmatrix}
    \mathbbm{1}_j & \ti{L}^t+i\ti{H}^t\\ 
    \ti{L}+i\ti{H} & \mathbbm{1}_{n-j}+2\ti{L}\ti{L}^t+i\ti{S_2}
    \end{pmatrix}
\end{align*}
and the corresponding determinant is given by
\begin{align*}
    &\det\Big(\frac{1}{2}(\mathbbm{1}_n+(B^{-1}AB)(B^{-1}AB)^t+iB^{-1}SB^{-t})\Big)\\
    &\quad=\det((\mathbbm{1}_{n-j}+\ti{L}\ti{L}^t+
    \ti{H}\ti{H}^t)+i(\ti{S_2}-\ti{H}\ti{L}^t-\ti{L}\ti{H}^t)).
\end{align*}
Next we set
\begin{align*}
    Q&=\mathbbm{1}_{n-j}+\ti{L}\ti{L}^t+\ti{H}\ti{H}^t,\\
    T&=\ti{S_2}-\ti{H}\ti{L}^t-\ti{L}\ti{H}^t. 
\end{align*}
Then the integral $I^j_{n,\kappa}(Z)$ in \autoref{intermediate integral 1} is given by
\begin{align*}
    I^j_{n,\kappa}(Z)=\int\limits_{S_2}\int\limits_{H}\int\limits_{L} 
\frac{\;\; [\dif S_2]\wedge[\dif H]\wedge[\dif L]}{\vert\det(Q+iT)\vert^{\kappa}}.
\end{align*}

Next we need to calculate the volume form $[\dif S_2]\wedge[\dif H]\wedge[\dif L]$ in terms of $T,\ti{L}$ and $\ti{H}$. From equations (\ref{tilde setting 1}) and (\ref{tilde setting 2}), we obtain
\begin{align*}
    2^{j(n-j)}[\dif\ti{L}]&=\frac{\det(P_1)^{n-j}}{\det(P_2)^j}[\dif L],\\
    2^{j(n-j)}[\dif\ti{H}]&=\frac{1}{\det(P_1)^{n-j}\det(P_2)^j}[\dif H]. 
\end{align*}
From \autoref{tilde setting 3}, one obtains
\begin{align*}
    2^{(n-j)(n-j+1)/2}[\dif\ti{S_2}]\wedge[\dif H]\wedge[\dif L]=\frac{[\dif S_2]\wedge[\dif H]\wedge[\dif L]}{\det(P_2)^{n-j+1}}. 
\end{align*}
Now, since $[\dif\ti{S_2}]\wedge[\dif \ti{H}]\wedge[\dif \ti{L}]=[\dif T]\wedge[\dif \ti{H}]\wedge[\dif \ti{L}]$, we conclude that
\begin{align*}
    [\dif S_2]\wedge[\dif H]\wedge[\dif L]=2^{(n-j)(n-j+1)/2+2j(n-j)}
    \det(P_2)^{n-j+1}\det(P_2)^{2j} [\dif T]\wedge[\dif \ti{H}]\wedge[\dif \ti{L}]. 
\end{align*}
Therefore, we have
\begin{align}\label{first estimate}
    I^j_{n,\kappa}(Z)\leq c_n \det(P_2)^{n-j+1}\det(P_2)^{2j} \int\limits_{\ti{L}}\int\limits_{\ti{H}}\int\limits_{T} 
\frac{\;\; [\dif T]\wedge[\dif \ti{H}]\wedge[\dif\ti{L}]}{\vert\det(Q+iT)\vert^{\kappa}},
\end{align}
where $c_n$, as usual, stands for a generic constant depending only on $n$. 
As the matrix $Q=\mathbbm{1}_{n-j}+\ti{L}\ti{L}^t+\ti{H}\ti{H}^t$ is positive definite, the integral 
\begin{align*}
    \int\limits_{T\in\mr{Sym}_{n-j}(\mb{R})}
\frac{ [\dif T]}{\vert\det(Q+iT)\vert^{\kappa}}
\end{align*}
can be written as
\begin{align*}
    \int\limits_{T\in\mr{Sym}_{n-j}(\mb{R})} 
\frac{ [\dif T]}{\vert\det(Q+iT)\vert^{\kappa}}
=\frac{1}{\det(Q)^{\kappa}}\int\limits_{T\in\mr{Sym}_{n-j}(\mb{R})} \frac{ [\dif T]}{\vert\det(\mathbbm{1}_{n-j}+iQ^{-1/2}TQ^{-1/2})\vert^{\kappa}}.
\end{align*}
Then setting $\ti{T}=Q^{-1/2}TQ^{-1/2}$, we have
\begin{align*}
    [\dif T]=\det(Q)^{(n-j+1)/2}[\dif \ti{T}].
\end{align*}
Then, using the Hua integral in (\ref{hua integral}), we obtain
\begin{align}\label{second estimate}
    \int\limits_{T\in\mr{Sym}_{n-j}(\mb{R})} 
\frac{ [\dif T]}{\vert\det(Q+iT)\vert^{\kappa}}
&=\frac{1}{\det(Q)^{\kappa-(n-j+1)/2}}\int\limits_{\ti{T}=\ti{T}^t}
     \frac{[\dif \ti{T}]}{\det(\mathbbm{1}_n+\ti{T}^2)^{\kappa/2}}\nonumber\\
    & \leq c_n \frac{\kappa^{-(n-j)(n-j+1)/4}}{\det(Q)^{\kappa-(n-j+1)/2}}.
\end{align}
Also, from $Y\leq(\kappa/2c_2(n))\mathbbm{1}_n$, i.e.,
\begin{align*}
    Y=BB^t=\begin{pmatrix}P_1 & 0\\ P & P_2\end{pmatrix}
    \begin{pmatrix}P_1^t & P^t\\ 0 & P_2^t\end{pmatrix}
    =\begin{pmatrix}P_1 P_1^t & P_1P^t\\ PP_1^t & P_2P_2^t+PP^t\end{pmatrix}\leq \frac{\kappa}{2c_2(n)}\mathbbm{1}_n,
\end{align*}
one obtains that 
\begin{align}\label{schur technique}
    P_2P_2^t+PP^t-PP_1^t(P_1P_1^t)^{-1}P_1P^t=P_2P_2^t\leq \frac{\kappa}{2c_2(n)}\mathbbm{1}_{n-j}. 
\end{align}
Thus we have $\det(P_2)\leq c_n \kappa^{(n-j)/2}$. Hence, by \autoref{first estimate} and \autoref{second estimate}, 
we have the estimate
\begin{align*}
    I^j_{n,\kappa}(Z)\leq& c_n\, \kappa^{-(n-j)(n-j+1)/4}\,\kappa^{(n-j)(n-j+1)/2}\,\kappa^{j(n-j)}\\
    &\quad \cdot \int\limits_{\ti{L}}\int\limits_{\ti{H}} 
\frac{[\dif \ti{H}]\wedge[\dif\ti{L}]}{\det(\mathbbm{1}_{n-j}+\ti{L}\ti{L}^t+\ti{H}\ti{H}^t)^{\kappa-(n-j+1)/2}}. 
\end{align*}
Now, to estimate the the integral
\begin{align*}
    \int\limits_{\ti{L}}\int\limits_{\ti{H}} 
\frac{[\dif \ti{H}]\wedge[\dif\ti{L}]}{\det(\mathbbm{1}_{n-j}+\ti{L}\ti{L}^t+\ti{H}\ti{H}^t)^{\kappa-(n-j+1)/2}},
\end{align*}
we set
\begin{align*}
    \mathbbm{1}_{n-j}+\ti{L}\ti{L}^t=EE^t,\quad E^{-1}\ti{H}=U. 
\end{align*}
Then, the above integral splits as
\begin{align*}
    &\int\limits_{\ti{L}}\int\limits_{\ti{H}} 
\frac{[\dif \ti{H}]\wedge[\dif\ti{L}]}{\det(\mathbbm{1}_{n-j}+\ti{L}\ti{L}^t+\ti{H}\ti{H}^t)^{\kappa-(n-j+1)/2}}\\
&\quad=
\int\limits_{\ti{L}\in\mb{R}^{(n-j)\times j}}
\frac{[\dif\ti{L}]}{\det(\mathbbm{1}_{n-j}+\ti{L}\ti{L}^t)^{\kappa-(n-j)-1/2}}
\int\limits_{U\in\mb{R}^{(n-j)\times j}}
\frac{[\dif U]}{\det(\mathbbm{1}_{n-j}+U U^t)^{\kappa-(n-j+1)/2}}.
\end{align*}
Proceeding as in \cite[Theorem 2.2.1]{Hua}, for matrices $X\in\mb{R}^{p\times q}\,(p,q\in\mb{N}_{\geq 1})$ one obtains the formula 
\begin{align*}
    \int\limits_{X\in \mb{R}^{p\times q}} 
    \frac{[\dif X]}{\det(\mathbbm{1}_{p}+X X^t)^{\mu}}
    =\pi^{pq/2}\prod\limits_{l=1}^q \frac{\Gamma(\mu-(l-1)/2-p/2)}{\Gamma(\mu-(l-1)/2)}\quad\quad(\mu>(p+q-1)/2). 
\end{align*}
Using this formula, it immediately follows that
\begin{align*}
    &\int\limits_{U\in\mb{R}^{(n-j)\times j}}
\frac{[\dif U]}{\det(\mathbbm{1}_{n-j}+U U^t)^{\kappa-(n-j+1)/2}} 
\leq c_n \kappa^{-j(n-j)/2},\\
 & \int\limits_{\ti{L}\in\mb{R}^{(n-j)\times j}}
\frac{[\dif\ti{L}]}{\det(\mathbbm{1}_{n-j}+\ti{L}\ti{L}^t)^{\kappa-(n-j)-1/2}}\leq c_n \kappa^{-j(n-j)/2}. 
\end{align*}
Thus, we have the estimate for the integral
\begin{align*}
    \int\limits_{\ti{L}}\int\limits_{\ti{H}} 
\frac{[\dif \ti{H}]\wedge[\dif\ti{L}]}{\det(\mathbbm{1}_{n-j}+\ti{L}\ti{L}^t+\ti{H}\ti{H}^t)^{\kappa-(n-j+1)/2}}\leq c_n \kappa^{-j(n-j)}, 
\end{align*}
thereby giving
\begin{align*}
    I^j_{n,\kappa}(Z)\leq& c_n\, \kappa^{(n-j)(n-j+1)/4}. 
\end{align*}
Hence, we have that
\begin{align*}
\sum_{\gamma\in\Gamma^j_{\infty}}\,\frac{1}{\prod_{j=1}^n \ch^{\kappa}(r_j^{\gamma}(Z))}\leq c_n\, \kappa^{(n-j)(n-j+1)/4} \quad\quad(0\leq j\leq (n-1))
\end{align*}
and consequently,
\begin{align*}
\sum_{\gamma\in\Gamma_{\infty}}\,\frac{1}{\prod_{j=1}^n \ch^{\kappa}(r_j^{\gamma}(Z))}\leq c_n\, \kappa^{n(n+1)/4}. 
\end{align*}

Thus, from \autoref{local pretheorem}, we have
\begin{align*}
\lim_{t\rightarrow\infty} \exp\big(-\frac{n\kappa}{4}(\kappa-(n+1))\,t\big)\sum_{\gamma\in\Gamma_{\infty}}K_t^{(\kappa)}(2R^{\gamma}(Z))
\leq 
c_{n}\, \kappa^{3n(n+1)/4},
\end{align*}
whence the theorem follows. 
\end{proof}


\subsection{Uniform sup-norm bounds}

\begin{theorem}\label{cover case}
Let $\Gamma_0\subsetneq \mathrm{Sp}_n(\mathbb{R})$ be a fixed arithmetic subgroup of $\mathrm{Sp}_n(\mathbb{R})$ such that $M_0:=\Gamma_0\backslash \mb{H}_n$ is of finite volume. Let $\Gamma\subseteq \Gamma_0$ a subgroup of finite index. Then, for $\kappa\geq n+1$, we have 
\begin{align*}
   \sup_{Z\in\mb{H}_n} S_{\kappa}^{\Gamma}(Z)\leq c_{n,\Gamma_0}\, \kappa^{3n(n+1)/4}\quad\quad( \kappa\geq n+1),
\end{align*}
where $c_{n,\Gamma_0}$ is a positive real constant depending only on $n$ and $\Gamma_0$.
\end{theorem}
\begin{proof}
As in the proof of \autoref{cofinite 2}, we denote by $\mathcal{C}_0$ the set of all inequivalent chains of boundary components of $M_0$ and choose boundary neighbourhoods $U_{\ve}(\mathcal{P}_0)\,(\mathcal{P}_0\in\mathcal{C}_0)$ containing the entire chain $\mathcal{P}_0$ such that the complement of their union in $M_0$, i.e.,
\begin{align*}
    K_{0,\ve}:=M_0\setminus \bigcup_{P_0\in\mathcal{C}_0} U_{\ve}(\mathcal{P}_0)
\end{align*}
is a compact subset of $M_0$. 
 
Let $M:=\Gamma\backslash\mb{H}_n$ and $\pi\colon M\rightarrow M_0$ denote the covering map. Then by means of $K_{0,\ve}$, we obtain the compact subset $K_{\ve}:=\pi^{-1}(K_{0,\ve})$ of $M$. Since $\Gamma\subseteq\Gamma_0$, from \autoref{heat connection}, by expanding the sum over $\Gamma$ to that over the larger group $\Gamma_0$, we have
 \begin{align*}
    S_{\kappa}^{\Gamma}(Z)&\leq \exp\big(-\frac{n\kappa}{4}(\kappa-(n+1))\,t\big)\sum_{\gamma\in\Gamma}K_t^{(\kappa)}(2R^{\gamma}(Z))\\
    &\leq \exp\big(-\frac{n\kappa}{4}(\kappa-(n+1))\,t\big)\sum_{\gamma\in\Gamma_0}K_t^{(\kappa)}(2R^{\gamma}(Z)),
\end{align*}
which, by \autoref{compact case}, gives the uniform bound
\begin{align}\label{compact gamma0}
   \sup_{Z\in K_{\ve}} S_{\kappa}^{\Gamma}(Z)\leq c_{n,\Gamma_0}\, \kappa^{n(n+1)/2}\quad\quad( \kappa\geq n+1).
\end{align}

We are thus left to bound the quantity $S_{\kappa}^{\Gamma}(Z)$ in the neighbourhoods of $M$ obtained by pulling back the neighbourhoods $U_{\ve}(\mathcal{P}_0)\,(\mathcal{P}_0\in\mathcal{C}_0)$ of $M_0$ to $M$. In order to do this, as in the proof of \autoref{cofinite 2}, we can again assume without loss of generality that $\mathcal{P}_0$ is the chain 
\begin{align*}
\Gamma_0\backslash\mb{H}_0< \Gamma_1\backslash\mb{H}_1<\ldots<\Gamma_j\backslash\mb{H}_j<\ldots<\Gamma_{n-1}\backslash\mb{H}_{n-1}
\end{align*}
of standard boundary components of $\Gamma_n\backslash\mb{H}_n$. Furthermore, we may also assume that the chain $\mathcal{P}\in\mathcal{C}$ of boundary components of $M$ lying over $\mathcal{P}_0$ is also the chain of standard boundary components of $\Gamma_n\backslash\mb{H}_n$ of ramification index $\ell$, say. Then a cusp form $f\in \mathcal{S}_{\kappa}^n(\Gamma)$ of weight $\kappa$ has a Fourier expansion (see \autoref{eq:fourier})
\begin{align*}
f(Z)=\sum_{\substack{T\in\mr{Sym}_n(\mb{Q}),\,T>0\\T\, \text{half-integral}}} a(T)\exp\bigg(\frac{2\pi i}{\ell}\:\tr(TZ)\bigg)
\end{align*}
at $\mathcal{P}$. Then, just like in \autoref{klingen bound} in \autoref{cofinite 2}, we obtain positive numbers $c_1(n)>0$ and $c_2(n)>0$ depending only on $n$ such that
\begin{align*}
    \vert f(Z)\vert \leq c_1(n)\, \exp(-c_2(n)\tr(Y)/\ell)
    \quad\quad(Z\in\mathscr{F}_n). 
\end{align*}
Then proceeding as in \autoref{cofinite 2} with $c_2(n)$ replaced by $c_2(n)/\ell$, one sees that for $\ve>\kappa\ell/(2c_2(n))$, we have
\begin{align*}
    \sup_{Z\in M} S_{\kappa}^{\Gamma}(Z)=\sup_{Z\in K_{\ve}} S_{\kappa}^{\Gamma}(Z)\quad\quad\bigg(\ve>\frac{\kappa\ell}{2c_2(n)}\bigg),
\end{align*}
which, by \autoref{compact gamma0} gives the uniform estimate
\begin{align*}
   \sup_{Z\in M} S_{\kappa}^{\Gamma}(Z)\leq c_{n,\Gamma_0}\, \kappa^{n(n+1)/2}\quad\quad( \kappa\geq n+1).
\end{align*}

Thus, we are left only to bound $S_{\kappa}^{\Gamma}(Z)$ in the range
$Y\leq (\kappa\ell/2c_2(n))\mathbbm{1}_n$. Again, as in \autoref{infinity split} in \autoref{cofinite 2}, we split the sum 
\begin{align*}
    S_{\kappa}^{\Gamma}(Z)= \lim_{t\rightarrow\infty}\exp\big(-\frac{n\kappa}{4}(\kappa-(n+1))\,t\big)\sum_{\gamma\in\Gamma}K_t^{(\kappa)}(2R^{\gamma}(Z))
\end{align*}
in \autoref{local Sk-inequality} into sums over $\Gamma\setminus\Gamma_{\infty}$ and $\Gamma_{\infty}$, with $\Gamma_{\infty}:=\Gamma\cap\Gamma_{0,\infty}$, to obtain
\begin{align}\label{infinity split 2}
\begin{aligned}
    S_{\kappa}^{\Gamma}(Z)\leq& \exp\big(-\frac{n\kappa}{4}(\kappa-(n+1))\,t\big)\sum_{\gamma\in\Gamma\setminus\Gamma_{\infty}}K_t^{(\kappa)}(2R^{\gamma}(Z)).\\
    &\quad +\lim_{t\rightarrow\infty} \exp\big(-\frac{n\kappa}{4}
    (\kappa-(n+1))\,t\big)
    \sum_{\gamma\in\Gamma_{\infty}}K_t^{(\kappa)}(2R^{\gamma}(Z)).
    \end{aligned}
\end{align}
Now as $\Gamma\setminus\Gamma_{\infty}\subseteq \Gamma_0\setminus\Gamma_{0,\infty}$, expanding the first sum to
$\Gamma_0\setminus\Gamma_{0,\infty}$ and using \autoref{gamma minus gamma infinity sum}, we obtain
\begin{align*}
 \exp\big(-\frac{n\kappa}{4}(\kappa-(n+1))\,t\big)
 \sum_{\gamma\in\Gamma\setminus\Gamma_{\infty}}K_t^{(\kappa)}(2R^{\gamma}(Z))\leq   c_{n,\Gamma_0}\, \kappa^{n(n+1)/2}. 
\end{align*}
Thus, it only remains to estimate the sum
\begin{align*}
    \sum_{\gamma\in\Gamma_{\infty}}\,\frac{1}{\prod_{j=1}^n \ch^{\kappa}(r_j^{\gamma}(Z))}
\end{align*}
in \autoref{local pretheorem}. Note that here we now have 
\begin{align*}
    \Gamma_{\infty}=\bigcup_{j=0}^{n-1}\Gamma_{\infty}^j,
\end{align*}
with $\Gamma_{\infty}^j=\Gamma\cap\Gamma_{0,\infty}^j$, i.e.,
\begin{align*}
    \Gamma_{\infty}^0&=\bigg\{\begin{pmatrix}\mathbbm{1}_n & \ell S\\ 0 & \mathbbm{1}_n\end{pmatrix}\;\bigg\vert\; S\in\mr{Sym}_n(\mb{Z})\bigg\},\\
    \Gamma_{\infty}^j&=\bigg\{\begin{pmatrix}A & A S\\ 0 & A^{-t}\end{pmatrix}\;\bigg\vert\; 
    A=\begin{pmatrix}\mathbbm{1}_j & 0\\ \ell L & \mathbbm{1}_{n-j}\end{pmatrix},\,
    S=\begin{pmatrix}0 & \ell H^t\\ \ell H & \ell S_2\end{pmatrix}\bigg\}\quad\quad(1\leq j\leq (n-1)), 
\end{align*}
where $\,L,H\in\mb{Z}^{(n-j)\times j}$ and $S_2\in\mr{Sym}_{n-j}(\mb{Z})$. 

Then, for $j=0$, proceeding as in \autoref{beta integral connection}, we have
\begin{align*}
    \sum_{\gamma\in\Gamma^0_{\infty}}\,\frac{1}{\prod_{j=1}^n \ch^{\kappa}(r_j^{\gamma}(Z))}\leq \int\limits_{S\in\mr{Sym}_n(\mb{R})}
    \frac{[\dif S]}{\det(\mathbbm{1}_n+(\frac{1}{2}Y^{-1/2}\ell SY^{-1/2})^2)^{\kappa/2}}. 
\end{align*}
Now, setting $T=\frac{1}{2}Y^{-1/2}\ell SY^{-1/2}$, we have
\begin{align*}
    [\dif T]=c_n \ell^{n(n+1)/2}\det(Y)^{-(n+1)/2}[\dif S],  
\end{align*}
thereby giving 
\begin{align*}
     \sum_{\gamma\in\Gamma^0_{\infty}}\,\frac{1}{\prod_{j=1}^n \ch^{\kappa}(r_j^{\gamma}(Z))}\leq c_n \ell^{-n(n+1)/2}\det(Y)^{(n+1)/2}\int\limits_{T\in\mr{Sym}_n(\mb{R})}
     \frac{[\dif T]}{\det(\mathbbm{1}_n+T^2)^{\kappa/2}}. 
\end{align*}
Now, from $\det(Y)\leq (\kappa\ell/(2c_2(n)))^n$, it easily follows that we have an uniform estimate
\begin{align*}
    \sum_{\gamma\in\Gamma^0_{\infty}}\,\frac{1}{\prod_{j=1}^n \ch^{\kappa}(r_j^{\gamma}(Z))}\leq c_n \kappa^{n(n+1)/4}
\end{align*}
independent of the ramification index $\ell$. 

Similarly, for the $j>0$ case, proceeding as in \autoref{cofinite 2}, in place of substitution equations \eqref{tilde setting 1}--\eqref{tilde setting 3}, we set
\begin{align*}
    2\ti{L}&=P_2^{-1}\ell LP_1, \\
    2\ti{H}&=P_2^{-1}\ell H P_1^{-t},\\
    2\ti{S_2}&=P_2^{-1}\ell S_2P_2^{-t}-P_2^{-1}(\ell HP_1^{-t}P^t+PP_1^{-1}\ell H^t)P_2^{-t}, 
\end{align*}
which results in
\begin{align*}
    I^j_{n,\kappa}(Z)\leq c_n\ell^{-(n-j+1)(n-j)/2-2j(n-j)} \det(P_2)^{n-j+1}\det(P_2)^{2j} \int\limits_{\ti{L}}\int\limits_{\ti{H}}\int\limits_{T} 
\frac{\;\; [\dif T]\wedge[\dif \ti{H}]\wedge[\dif\ti{L}]}{\vert\det(Q+iT)\vert^{\kappa}}
\end{align*}
in place of \eqref{first estimate}. Now, with $\det(P_2)\leq c_n (\ell \kappa)^{(n-j)/2}$ coming from $Y\leq (\kappa\ell/2c_2(n))\mathbbm{1}_n$ via \autoref{schur technique}, it follows that 
\begin{align*}
  \sum_{\gamma\in\Gamma^j_{\infty}}\,\frac{1}{\prod_{j=1}^n \ch^{\kappa}(r_j^{\gamma}(Z))}  I^j_{n,\kappa}(Z)\leq& c_n\,\ell^{-j(n-j)} \kappa^{(n-j)(n-j+1)/4}\leq c_n \kappa^{(n-j)(n-j+1)/4}.
\end{align*}
Thus, combined, we get an uniform estimate 
\begin{align*}
\sum_{\gamma\in\Gamma_{\infty}}\,\frac{1}{\prod_{j=1}^n \ch^{\kappa}(r_j^{\gamma}(Z))}\leq c_n\, \kappa^{n(n+1)/4}
\end{align*}
resulting in the uniform estimate 
\begin{align*}
\lim_{t\rightarrow\infty} \exp\big(-\frac{n\kappa}{4}(\kappa-(n+1))\,t\big)\sum_{\gamma\in\Gamma_{\infty}}K_t^{(\kappa)}(2R^{\gamma}(Z))
\leq 
c_{n}\, \kappa^{3n(n+1)/4}
\end{align*}
in the second sum in \autoref{infinity split 2}, thereby proving the theorem. 
\end{proof}


\bibliographystyle{amsplain}
\bibliography{bibliogen}
\end{document}